\newif\ifjota\jotafalse

\PassOptionsToPackage{usenames,dvipsnames,table}{xcolor}
\ifjota
  \documentclass[smallextended,envcountsect,]{svjour3}
  \smartqed
  \usepackage{graphicx}
  \journalname{JOTA}
\else
  \documentclass[USenglish,10pt,a4paper]{article}
\fi
  \PassOptionsToPackage{noend}{algpseudocode}
  \PassOptionsToPackage{capitalize,nameinlink}{cleveref}
  \PassOptionsToPackage{breakspaces}{clevethm}
  \PassOptionsToPackage{dvipsnames}{xcolor}

  \usepackage[cal=boondoxo, scr=euler]{mathalfa}
  \ifjota
    \usepackage{myPreamble}
  \else
    \usepackage[bottom=1.18in,left=1.1in,right=1.1in]{geometry}
    \usepackage[
      alglinenumber,
      eqreset=section,
      thmreset=section,
    ]{myPreamble}
  \fi

\makeatletter
  \foreach \x in {gray, blue, orange, red, black} {%
    \edef\temp{%
      \noexpand\expandafter\noexpand\gdef\noexpand\csname\x\noexpand\endcsname{\noexpand\color{\x}}
    }%
    \temp
  }

  \let\OLDand\and
  \def\and{\texorpdfstring{\OLDand}{, }}%
  \ifjota
    \newcommand{\Sep}{\,\(\cdot\)\,}
  \else
    \newenvironment{keywords}{\par\noindent{\bf Keywords. }}{}
    \newenvironment{AMS}{\par\noindent{\bf AMS subject classifications. }}{}
    \newcommand{\Sep}{ \(\cdot\)\ }
  \fi

  \ifjota
    
    \CreateTheoremLists{assumption}
    \CreateTheoremLists{corollary}
    \CreateTheoremLists{definition}
    \CreateTheoremLists{example}
    \CreateTheoremLists{lemma}
    \CreateTheoremLists{proposition}
    \CreateTheoremLists{remark}
    \CreateTheoremLists{theorem}
    \CreateTheoremLists{stepsize}
    \CreateTheoremLists{relaxation}

    \makeatletter
      \let\oldproposition\proposition
      \let\oldendproposition\endproposition
      \renewenvironment{proposition}[1][]{%
        \def\@currentlabelname{#1}\ifstrempty{#1}{\oldproposition}{\oldproposition[#1]}%
      }{\oldendproposition}
      \let\oldlemma\lemma
      \let\oldendlemma\endlemma
      \renewenvironment{lemma}[1][]{%
        \def\@currentlabelname{#1}\ifstrempty{#1}{\oldlemma}{\oldlemma[#1]}%
      }{\oldendlemma}
      \let\oldtheorem\theorem
      \let\oldendtheorem\endtheorem
      \renewenvironment{theorem}[1][]{%
        \def\@currentlabelname{#1}\ifstrempty{#1}{\oldtheorem}{\oldtheorem[#1]}%
      }{\oldendtheorem}
      \let\oldcorollary\corollary
      \let\oldendcorollary\endcorollary
      \renewenvironment{corollary}[1][]{%
        \def\@currentlabelname{#1}\ifstrempty{#1}{\oldcorollary}{\oldcorollary[#1]}%
      }{\oldendcorollary}
    \makeatother

    \foreach \x in {assumption, corollary, definition, example, fact, lemma, proposition, remark, theorem, stepsize, relaxation} {
      \edef\y{\expandafter\euppercase\x}
      \RegisterTheoremName{\x}{\y}
    }
    \newenvironment{assumption}{\begin{ass}}{\end{ass}}
    \newenvironment{stepsize}{\begin{step}}{\end{step}}
    \newenvironment{relaxation}{\begin{relaxat}}{\end{relaxat}}
  \else
    \newenvironment{assumption}{\begin{ass}}{\end{ass}}
    \newenvironment{corollary}{\begin{cor}}{\end{cor}}
    \newenvironment{definition}{\begin{defin}}{\end{defin}}
    \newenvironment{example}{\begin{es}}{\end{es}}
    \newenvironment{lemma}{\begin{lem}}{\end{lem}}
    \newenvironment{proposition}{\begin{prop}}{\end{prop}}
    \newenvironment{remark}{\begin{rem}}{\end{rem}}
    \newenvironment{theorem}{\begin{thm}}{\end{thm}}
    \newenvironment{stepsize}{\begin{step}}{\end{step}}
    \newenvironment{relaxation}{\begin{relaxat}}{\end{relaxat}}
  \fi
  
  \newlist{claims}{enumerate}{1}
    \setlist[claims,1]{
      label={\it Claim \oldstylenums{\arabic*}:},
      ref={\oldstylenums{\arabic*}},
      partopsep=0pt,
      parsep=0pt,
      itemsep=3pt,
      wide,
      labelindent=0pt,
    }
    \newlist{claims*}{enumerate}{1}
    \setlist[claims*,1]{
      label={\it Contradiction claim \oldstylenums{\arabic*}$^*$:},
      ref={\oldstylenums{\arabic*}$^*$},
      partopsep=0pt,
      parsep=0pt,
      itemsep=3pt,
      wide,
      labelindent=0pt,
    }
    \Crefname{claimsi}{Claim}{Claims}

  \allowdisplaybreaks

  \makeatletter
    \renewcommand{\clevethm@proofsectiontitle}{Omitted proofs of }
  \makeatother
  \Crefname{assumption}{Assumption}{Assumptions}
  \Crefname{example}{Example}{Examples}
  \Crefname{fact}{Fact}{Facts}
  \Crefname{remark}{Remark}{Remarks}
  \Crefname{cor}{Corollary}{Corollaries}
  \crefname{ALG@line}{step}{steps}
  \crefname{enumeratpropi}{property}{properties}
  \crefname{enumeratpropii}{property}{properties}
  \counterwithout{algorithm}{section}
  \counterwithout{table}{section}
  \renewcommand{\thealgorithm}{\arabic{algorithm}}
  
  \crefformat{equation}{(#2#1#3)}

  \crefname{table}{table}{tables}
  \Crefname{table}{Table}{Tables}

\usepackage{hyphenat}
\usepackage{cases} 
\usepackage{diagbox}
\usepackage{tabularray}
\UseTblrLibrary{diagbox}
\NewColumnType{M}[1]{Q[m,c,#1]}  
\usepackage{makecell}
\usepackage{makebox}
\usepackage{calrsfs} 
\DeclareMathAlphabet{\pazocal}{OMS}{zplm}{m}{n}

\makeatletter
\newsavebox{\@brx}
\newcommand{\llangle}[1][]{\savebox{\@brx}{\(\m@th{#1\langle}\)}%
  \mathopen{\copy\@brx\mkern2mu\kern-0.9\wd\@brx\usebox{\@brx}}}
\newcommand{\rrangle}[1][]{\savebox{\@brx}{\(\m@th{#1\rangle}\)}%
  \mathclose{\copy\@brx\mkern2mu\kern-0.9\wd\@brx\usebox{\@brx}}}
\makeatother

\makeatletter
\newcommand{\qindef}{\@ifstar{\@qindefb}{\@qindefa}}
\newcommand{\@qindefb}[2]{\inner{#1,#1}_{#2}}
\newcommand{\@qindefa}[2]{\operatorname{q}_{#2}(#1)}
\makeatother

\newcommand{\qsemi}[2]{\operatorname{q}_{#2}(#1)}

\def\A{A}
\def\B{B}
\def\mon{\mu}
\def\com{\rho}
\def\monA{\mon_\A}
\def\monB{\mon_\B}
\def\comA{\com_\A}
\def\comB{\com_\B}

\def\Mon{M}
\def\Com{R}
\def\MonA{\Mon_\A}
\def\MonB{\Mon_\B}
\def\ComA{\Com_\A}
\def\ComB{\Com_\B}
\def\parsumset{\dom_{\Box}}
\def\etamin{\bar\eta}

\def\M{P}                           
\def\projQ{\proj_{\range{\M}}}      

\def\DRSRho{V}
\def\DRSrhocom{\beta_{\rm P}}
\def\DRSrhomon{\beta_{\rm D}}

\def\s1{\|L\|}
\def\UL{Y}
\def\VL{X}
\def\SL{\Sigma}
\def\basisP{U}

\newcommand\optional[1]{[#1]}

\newcommand\tp[1]{#1^\top}
\newcommand\inv[1]{#1^{-1}}

\def\Tpd{T_{\!\,\mathrm{PD}}}
\def\Tp{T_{\!\,\mathrm{P}}}
\def\Td{T_{\!\,\mathrm{D}}}

\newcommand\other[1]{{#1}^\prime}

\newcommand{\x}{\bar x}
\newcommand{\y}{\bar y}
\newcommand{\z}{\bar z}

\newcommand{\D}{\mathcal{D}}

\renewcommand\thealgorithm{\arabic{algorithm}}

\makeatletter
  \renewcommand{\appendixproof@toc}{subsubsection}

\makeatother

\newcommand{\dfn}{\coloneqq}
\newcommand{\gph}{\graph}
\newcommand{\rge}{\range}
\DeclarePairedDelimiter{\nrm}{\lVert}{\rVert}
\DeclarePairedDelimiter{\inner}{\langle}{\rangle}

\let\oldBox\Box
\renewcommand{\Box}{\mathbin{\oldBox}}
\newcommand\sym[1]{\mathop{\mathbb{S}}^{#1}}

\pgfplotsset{
  compat=1.17,
  colormap={paper}{rgb255(0cm)=(230,230,230); rgb255(1cm)=(200,200,200); rgb255(2cm)=(185,185,185)}
  }
\usepgfplotslibrary{fillbetween, colormaps}
\usetikzlibrary{arrows.meta, calc, shapes.geometric, perspective}
\usetikzlibrary{decorations.pathreplacing,decorations.markings,arrows,shapes.misc}
\pgfplotscolorbarCMYKworkaroundfalse
\showgrayfalse
\usepackage{caption}

\definecolor{PaperBlue}{HTML}{185477}
\definecolor{PaperGreen}{HTML}{228B22}
\definecolor{PaperOrange}{HTML}{E69C24}
\definecolor{PaperRed}{HTML}{902A3C}

  \newif\ifshowold\showoldfalse
  \newif\ifshownew\shownewfalse
  \AtBeginDocument{%
    \ifshownew\else
      \renewenvironment{claims}{\let\item\relax}{}%
      \ifshowold
        \PackageWarning{DioBestia}{WARNING: Showing only OLD version}%
      \else
        \PackageWarning{DioBestia}{ERROR: One at least between \protect\ifshow old and \protect\ifshow new should be True}%
        \PackageError{DioBestia}{One at least between \protect\ifshow old and \protect\ifshow new should be True}{}%
      \fi
    \fi
  }

    \colorlet{newcolor}{orange!70!red}
    \colorlet{oldcolor}{black!30}
    
    \newcommand{\disablecolorlinks}{\def\HyColor@UseColor##1{}}

    \newcommand{\citenum}[1]{{%
      \def\@cite##1##2{{##1\if@tempswa , ##2\fi}}%
      \cite{#1}%
    }}%

  \newcounter{saveTheorem}
  \newcounter{saveEquation}

  \makeatletter%
  \newcommand{\mytag}[2]{%
    \text{#1}%
    \@bsphack
    \begingroup
      \@onelevel@sanitize\@currentlabelname
      \edef\@currentlabelname{%
        \expandafter\strip@period\@currentlabelname\relax.\relax\@@@%
      }%
      \protected@write\@auxout{}{%
        \string\newlabel{#2}{%
          {#1}%
          {\thepage}%
          {\@currentlabelname}%
          {\@currentHref}{}%
        }%
      }%
    \endgroup
    \@esphack
  }
  \makeatother
\pgfplotsset{compat=1.16}
\usepgfplotslibrary{fillbetween}
\showgrayfalse

\shownewtrue               
\newcommand{\TheShortTitle}{%
  Chambolle--Pock Algorithm without Monotonicity
}
\newcommand{\TheTitle}{
  Convergence of the Chambolle--Pock Algorithm in the Absence of Monotonicity
  }

\newcommand{\TheFunding}{%
  This work was supported by:
  the Research Foundation Flanders (FWO) PhD grant No. 1183822N;
  postdoctoral grant No. 12Y7622N and research projects G081222N, G033822N, G0A0920N;
  European Union's Horizon 2020 research and innovation programme under the Marie Skłodowska-Curie grant agreement No. 953348.%
}  
\newcommand{\TheKeywords}{%
convex/nonconvex optimization\Sep
\ifjota
  monotone/nonmonotone inclusions\Sep
\else
  monotone/nonmonotone inclusion problems\Sep
\fi
Chambolle--Pock algorithm\Sep
primal-dual hybrid gradient\Sep
semimonotone operators
}
\newcommand{\TheAMSsubj}{%
  47H04\Sep 
  49J52\Sep 
  49J53\Sep
  65K15\Sep 
  90C26. 
}
\newcommand{\TheAbstract}{%
The Chambolle--Pock algorithm (CPA), also known as the primal-dual hybrid gradient method, has gained popularity over the last decade due to its success in solving large-scale convex structured problems. 
This work extends its convergence analysis for problems with varying degrees of (non)monotonicity, quantified through a so-called oblique weak Minty condition on the associated primal-dual operator. %
Our results reveal novel stepsize and relaxation parameter ranges 
which do not only depend on the norm of the linear mapping, but also on its other singular values.
In particular, in nonmonotone settings, in addition to the classical stepsize conditions, extra bounds on the stepsizes and relaxation parameters are required. 
On the other hand, in the strongly monotone setting, the relaxation parameter is allowed to exceed the classical upper bound of two. 
Moreover, we build upon the recently introduced class of semimonotone operators,
providing sufficient convergence conditions for CPA when the individual operators are semimonotone. 
Since this class of operators encompasses traditional operator classes including (hypo)- and co(hypo)-monotone operators, this analysis recovers and extends existing results for CPA.
Tightness of the proposed stepsize ranges is demonstrated through several examples.
}

\author{%
    Brecht Evens\textsuperscript{1}
    \and
    Puya Latafat\textsuperscript{2}
    \and
    Panagiotis Patrinos\textsuperscript{1}
    \thanks{
      \TheFunding\newline
      \indent
      \textsuperscript{1}\TheAddressKU,
      {\tt
        \{%
          \href{mailto:brecht.evens@kuleuven.be}{brecht.evens},%
          \href{mailto:panos.patrinos@kuleuven.be}{panos.patrinos}%
          \}\texttt{@kuleuven.be}.
      }\newline
      \indent
      \textsuperscript{2}\TheAddressIMT, 
      {\tt
          \href{mailto:puya.latafat@imtlucca.it}{puya.latafat@imtlucca.it}.
      }
    }
  }
\renewcommand\footnotemark{}

\ifjota
  
  \title{%
    \TheTitle%
  }%

  \author{
    Brecht Evens\textsuperscript{1}
    \and
    Puya Latafat\textsuperscript{2}
    \and
    Panagiotis Patrinos\textsuperscript{1}
  }

  \institute{
    \TheFunding\newline
      \indent
      \textsuperscript{1}\TheAddressKU,
      {\tt
        \{%
          \href{mailto:brecht.evens@kuleuven.be}{brecht.evens},%
          \href{mailto:panos.patrinos@kuleuven.be}{panos.patrinos}%
          \}\texttt{@kuleuven.be}.
      }\newline
      \indent
      \textsuperscript{2}\TheAddressIMT, 
      {\tt
          \href{mailto:puya.latafat@imtlucca.it}{puya.latafat@imtlucca.it}.
      }
  }

  \date{Received: date / Accepted: date}
\else
  \title{\texorpdfstring{\TheTitle}{\TheShortTitle}}
  \date{}
\fi

\begin{document}
    \maketitle

    \begin{abstract}
      \TheAbstract
    \end{abstract}
    
    \ifjota
\else
    \fi

\section{Introduction}

This paper considers composite inclusion problems of the form 
\begin{equation}\label{prob:composite}\tag{P-I}
    \text{find} \quad x \in \R^n \quad \text{such that} \quad 0 \in \Tp x \dfn Ax + L^\top B Lx, 
\end{equation}
where $A : \R^n \rightrightarrows \R^n$, $B : \R^m \rightrightarrows \R^m$ are two (possibly nonmonotone) operators and \(L \in \R^{m\times n}\) is a nonzero matrix.  
Problems of this form emerge naturally in a wide variety of applications in optimization and variational analysis.
For instance, in the framework of convex optimization, the first-order optimality condition for minimizing $g(x) + h(Lx)$ is of the form \cref{prob:composite},
with \(A = \partial g\) and \(B = \partial h\) representing the subdifferentials of proper lsc convex functions \(g\) and \(h\). 

One of the central algorithms for solving \eqref{prob:composite} is the Chambolle--Pock algorithm (CPA) \cite{chambolle2011firstorder} (also known as the primal-dual hybrid gradient (PDHG) method \cite{zhu2008efficient,esser2010general,he2012Convergence}).
Given strictly positive stepsizes 
\(\gamma,\tau > 0\),
a sequence of strictly positive relaxation parameters \(\seq{\lambda_k}\) and an initial guess \((x^0,y^0)\in\R^{n+m}\), this algorithm consists of the following iterates.
\begin{equation}
    \begin{cases}\tag{CPA}\label{eq:CP}
        \displaystyle
        \bar{x}^k
            {}\in{}&
        J_{\gamma \A}\big(x^k - \gamma L^\top y^k\big)\\
        \bar y^k 
            {}\in{}&
        J_{\tau \inv\B}\big(y^k + \tau L(2\bar x^k - x^k)\big)\\
        x^{k+1}&=x^k+\lambda_k(\bar{x}^k-x^k)\\
        y^{k+1}&=y^k+\lambda_k(\bar{y}^k-y^k)
    \end{cases}
\end{equation}
The convergence analysis of \ref{eq:CP} in literature largely relies upon an underlying monotonicity assumption. In this work, we identify classes of nonmonotone problems along with corresponding stepsize and relaxation parameter conditions for which \ref{eq:CP} remains convergent.
To this end, we rely on casting \ref{eq:CP} as an instance of the preconditioned proximal point algorithm (PPPA).
This connection was previously exploited in \cite{esser2010general,he2012Convergence,condat2013primal,latafat2017Asymmetric,bredies2022degenerate} in the monotone setting.
Many other widely used numerical methods can also be interpreted as special cases of \ref{eq:PPPA-intro}, see e.g. \cite{eckstein1988LionsMercier,rockafellar1976augmented,rockafellar2019solving,eckstein1992DouglasRachford,condat2013primal}. 
In particular, consider the inclusion problem of finding a zero of a set-valued operator \(T:\R^n \rightrightarrows \R^n\), \ie,
\begin{equation} \label{prob:P1} \tag{G-I} 
    \text{find} \quad z \in \R^n \quad \text{such that} \quad 0 \in T z. 
\end{equation}
Then, given a symmetric positive \emph{semidefinite} preconditioning matrix \(\M \in \R^{n\times n}\) and a sequence of strictly positive relaxation parameters \(\seq{\lambda_k}\), the (relaxed) preconditioned proximal point algorithm applied to \eqref{prob:P1} consists of the fixed point iterations
\begin{equation}\label{eq:PPPA-intro}\tag{PPPA}
    \begin{cases}
        \bar{z}^k&\in(\M+T)^{-1}\M z^k\\
        z^{k+1}&=z^k+\lambda_k(\bar{z}^k-z^k).
    \end{cases}
\end{equation}
The Chambolle--Pock algorithm can be retrieved for certain choice of the preconditioner $\M$ and the operator $T$. 
Specifically, consider the so-called \emph{primal-dual inclusion}
\begin{equation}\label{eq:primaldual}\tag{PD-I}
    \text{find} \quad z = (x,y) \in \R^{n+m} \quad \text{such that} \quad 
    0\in \Tpd (z)
        {}\coloneqq{}
    \begin{bmatrix}
        Ax\\ B^{-1}y\end{bmatrix}+\begin{bmatrix}L^\top y\\ - L x
    \end{bmatrix}    
    .
\end{equation}
Then, letting \(z^k = (x^k, y^k)\) and \(\bar z^k = (\bar x^k, \bar y^k)\), \ref{eq:CP} is equivalent to applying \ref{eq:PPPA-intro} to \eqref{eq:primaldual}, with preconditioner
\begin{equation}\label{eq:CP:P}
    \M
        {}={}
    \begin{bmatrix}
        \tfrac1\gamma\I_n&-L^\top\\
        -L&\tfrac1\tau \I_m
    \end{bmatrix}.
\end{equation} 
As a result of this equivalence, the convergence properties of \ref{eq:CP} can be inferred from those of \ref{eq:PPPA-intro}.
In the monotone setting, convergence of \ref{eq:PPPA-intro} is well understood, not only for positive definite preconditioners \cite{martinet1970Regularisation,rockafellar1976augmented,rockafellar1976Monotone,rockafellar2019Progressive} but also for positive semidefinite ones \cite[Thm. 3.4]{latafat2017Asymmetric}, \cite[\S 2.1]{bredies2022degenerate}.
Analogously, the convergence of \ref{eq:CP} for monotone inclusions is relatively well-understood, provided that the stepsizes $\gamma$ and $\tau$ satisfy a certain stepsize condition.
The standard assumption in the first works on \ref{eq:CP} such as \cite{chambolle2011firstorder,esser2010general,he2012Convergence} was that the stepsizes $\gamma$ and $\tau$ satisfy $\gamma\tau\nrm{L}^2 < 1$.
This assumption was later relaxed to $\gamma\tau\nrm{L}^2 \leq 1$ in \cite{condat2013primal,latafat2017Asymmetric,oconnor2020equivalence}, broadening the scope of the analysis to Douglas-Rachford splitting (DRS), for which $\tau = \nicefrac1\gamma$ and $L = \I$. 
Interestingly, when interpreting \ref{eq:CP} as a particular instance of \ref{eq:PPPA-intro},
the stepsize condition discussed in these works is directly linked to the positive definiteness of the preconditioning matrix $\M$ in \ref{eq:PPPA-intro}. 
This connection becomes evident by observing that, owing to the Schur complement lemma, \(\M\) is positive definite under the traditional stepsize condition $\gamma\tau\nrm{L}^2 < 1$
and positive semidefinite under the relaxed stepsize condition $\gamma\tau\nrm{L}^2 \leq 1$.
Note that the extended stepsize range $\gamma\tau\nrm{L}^2 \leq \nicefrac43$ was examined in \cite{yan2024improved,banert2023chambollepock}, although they rely upon a more involved Lyapunov analysis.
Recently, convergence of \ref{eq:PPPA-intro} in the nonmonotone setting has been considered in \cite{evens2023convergence} under the assumption that $T$ admits a set of oblique weak Minty solutions, defined as follows.

\begin{definition}[$\DRSRho$\hyp{}oblique weak Minty solutions \cite{evens2023convergence}]\label{def:SWMVI}
    Let $\DRSRho\in\R^{n \times n}$ be symmetric.
    An operator \(T : \R^n \rightrightarrows \R^n\) is said to have $\DRSRho$\hyp{}oblique weak Minty solutions at (a nonempty set) \(\pazocal{S}^\star\subseteq\zer T\) 
    on a set $\pazocal U \subseteq \R^n \times \R^n$ if
    \begin{equation}\label{def:WMVI}
         \langle v,z-z^\star\rangle\geq \qindef{v}{\DRSRho},\qquad \text{for all $z^\star\in\pazocal{S}^\star, (z,v)\in\graph T \cap \pazocal{U}
        $}, 
    \end{equation}
    where the quadratic form \(\qindef{v}{\DRSRho} \dfn \langle v, \DRSRho v \rangle\).
    Whenever $\pazocal{U} = \R^n \times \R^n$ we refrain from explicitly mentioning the set $\pazocal{U}$, and whenever $\DRSRho = \rho \I$ for some $\rho \in \R$, we refer to them as $\rho$\hyp{}weak Minty solutions.
\end{definition}
One key aspect of this assumption is its generality, as $\DRSRho$ is allowed to be any (possibly indefinite) symmetric matrix.
For instance, if $\DRSRho$ is equal to the zero matrix, \eqref{def:WMVI} reduces to the classic Minty variational inequality (MVI) \cite{minty1962,giannessi1998Minty},
while if $\DRSRho = \rho \I$ the so-called weak MVI is retrieved.
In literature, weak MVI and the closely related notion of cohypomonotonicity have been employed in the context of the extragradient and the forward-backward-forward method \cite{diakonikolas2021Efficient,pethick2021Escaping,pethick2023solving,bohm2023solving,gorbunov2023convergence}, as well as the classic PPA method \cite{pennanen2002Local,iusem2003Inexact,combettes2004proximal,gorbunov2023convergence}.

Leveraging the results from \cite{evens2023convergence} and the primal-dual connection between \ref{eq:CP} and \ref{eq:PPPA-intro}, the first part of this work will focus on establishing convergence of \ref{eq:CP} under the assumption that the primal-dual operator $\Tpd$ admits a set of $\DRSRho$\hyp{}oblique weak Minty solutions. 
To account for the inherent structure present within $\Tpd$, we impose a specific block diagonal form for $\DRSRho = \blkdiag(\DRSRho_{\rm P}, \DRSRho_{\rm D})$, which depends on the fundamental subspaces of $L$ (see \eqref{eq:WMSDRS} and the discussion thereafter).
Furthermore, we demonstrate that by restricting our obtained results to the case where $L = \I$ and $\tau = \nicefrac1\gamma$, the convergence results for nonmonotone 
DRS from \cite[Sec. 3]{evens2023convergence} are retrieved.

In contrast to the DRS setting, where the convergence results follow in a straightforward manner from those of \ref{eq:PPPA-intro} (see proof of \cite[Thm. 3.3]{evens2023convergence}), convergence results for \ref{eq:CP} are more challenging to obtain, not only due to additional stepsize parameter $\tau$, but mainly due to the additional complexity in the algorithm introduced by the matrix $L$.
This difficulty is overcome by considering the singular value decomposition of \(L\) and using the corresponding orthonormal basis to carefully decompose the preconditioner $\M$ and the oblique weak Minty matrix $\DRSRho$ (see proof of \Cref{thm:CP:full}).

In practice, it might be difficult to determine whether the associated primal-dual operator of a given inclusion problem admits $\DRSRho$-oblique weak Minty solutions.  
This issue will be addressed in the second part of this work, where we introduce the class of $(\Mon, \Com)$\hyp{}semimonotone operators and provide several calculus rules for this class, allowing to verify the existence of $\DRSRho$-oblique weak Minty solutions 
based on the semimonotonicity properties of the underlying operators $\A$ and $\B$.
The class of semimonotone operators is defined as follows.

\begin{definition}[semimonotonicity]\label{def:semimonotonicity}
    Let $\Mon,\Com\in\R^{n \times n}$ be symmetric (possibly indefinite) matrices. 
    An operator 
    $A : \R^n \rightrightarrows \R^n$ is said to be $(\Mon,\Com)$\hyp{}semimonotone at $(\other{x}, \other{y}) \in \gph A$ if
    \begin{equation}\label{eq:obliquequasisemimonotonicity:matrix}
        \inner*{x - \other{x}, y - \other{y}} \geq \qindef{x - \other{x}}{\Mon} + \qindef{y - \other{y}}{\Com}, \qquad \text{for all $(x, y) \in \gph A$},
    \end{equation}
    where ${\displaystyle \qindef{\cdot}{X} \dfn \inner{\cdot,\cdot}_X}$ for any symmetric matrix $X\in\R^{n \times n}$.
    An operator \(A\) is said to be $(\Mon,\Com)$\hyp{}semimonotone if it is $(\Mon,\Com)$\hyp{}semimonotone at all $(\other{x}, \other{y})\in \gph A$.
    It is said to be maximally \((\Mon,\Com)\)-semimonotone if its graph is not strictly contained in the graph of another \((\Mon,\Com)\)-semimonotone operator. 

    Throughout, whenever
    \(\Mon = \mon \I_n\) and \(\Com = \com \I_n\) where $\mon, \com \in \R$, 
    the prefix $(\Mon, \Com)$ is replaced by $(\mon,\com)$ and condition \eqref{eq:obliquequasisemimonotonicity:matrix} reduces to
    \begin{equation}\label{eq:quasisemimonotonicity}
        \inner*{x - \other{x}, y - \other{y}} \geq \mon \nrm{x - \other{x}}^2 + \com \nrm{y - \other{y}}^2, \qquad \text{for all $(x, y) \in \gph A$}.
    \end{equation}
\end{definition}
The class of $(\mon,\com)$-semimonotone operators was introduced in \cite[Sec. 4]{evens2023convergence} and enjoys a lot of additional freedom compared to more traditional operators classes.
For instance, it encompasses the classes of (hypo)monotone, co(hypomonotone), $\rho$\hyp{}semimonotone \cite[Def. 2]{otero2011Regularity}, averaged and firmly nonexpansive operators (see \cite[Rem. 4.2 \& Fig. 5]{evens2023convergence}).

In this work, this notion is generalized by characterizing the operator class with matrices $(\Mon, \Com)$ instead of scalars $(\mon, \com)$. This generalization is crucial to capture and exploit the specific structure emerging in \ref{eq:CP}.
To illustrate this, the next theorem provides a simplified version of our main result (see \Cref{cor:CP:semi} for the full statement).
For instance, if $\monA$ is positive, $(\monA L^\top L, \comA \I_n)$\hyp{}semimonotonicity of $A$ in \cref{cor:CP:semi:simple} could be replaced by $(\monA \nrm{L}^2, \comA)$\hyp{}semimonotonicity, which is in general a much more restrictive assumption.

\begin{theorem}[convergence of \ref{eq:CP} under semimonotonicity (simplified)]\label{cor:CP:semi:simple}
    Let operators $A:\R^n \rightrightarrows \R^n$ and $B:\R^m \rightrightarrows \R^m$ be outer semicontinuous. 
    Suppose that there exists $(x^\star, y^\star) \in \zer \Tpd$ such that $\A$ is $(\monA \tp L L, \comA \I_n)$\hyp{}semimonotone 
    at $(x^\star, -L^\top y^\star) \in \gph A$, $\B$ is
    $(\monB \I_m, \comB L \tp L)$\hyp{}semimonotone
    at $(L x^\star, y^\star) \in \gph B$
    and the
    \ifjota\else
        semimonotonicity
    \fi
    moduli $(\monA, \monB, \comA, \comB) \in \R^4$ satisfy one of the following conditions.
    \begin{enumerate}
        \item \label{it:cor:CP:semi:simple:1} (either) \(\monA = \monB = 0\) and $\comA = \comB = 0$.
        \item \label{it:cor:CP:semi:simple:2} (or) $\monA + \monB > 0$ and $\comA = \comB = 0$.
        \item \label{it:cor:CP:semi:simple:3} (or) $\comA + \comB > 0$ and \(\monA = \monB = 0\).
        \item \label{it:cor:CP:semi:simple:4} (or) $\monA + \monB > 0$, $\comA + \comB > 0$ and 
        $
            \min\{0, \frac{\monA\monB}{\monA + \monB}\}\min\{0, \frac{\comA\comB}{\comA + \comB}\} < \tfrac{1}{4\nrm{L}^2}
        $.
    \end{enumerate}
    Then, there exist positive stepsizes $\gamma, \tau$ and relaxation sequences \(\seq{\lambda_k}\) such that if the resolvents $J_{\gamma A}, J_{\tau B^{-1}}$ have full domain\footnote{The full domain assumption is imposed to ensure that the iterates of \ref{eq:CP} are well-defined.}, any sequence $\seq{\z^k} = \seq{\x^k,\y^k}$ generated by \ref{eq:CP} starting from $z^0 \in \R^n$ either reaches a point $\bar z^k \in \zer \Tpd$ in a finite number of iterations or every limit point of $(\z^k)_{k\in\N}$ belongs to $\zer \Tpd$.
\end{theorem}
This convergence result possesses two primary attributes that deserve attention.
First of all, it only requires semimonotonicity of the involved operators at a single point, as opposed to the traditional global assumptions of (hypo)monotonicity and co(hypo)monotonicity.
Secondly, by considering the more general class of semimonotone operators, we obtain fundamentally new convergence results (see case \ref{it:cor:CP:semi:simple:4}), not covered by any existing theory for \ref{eq:CP}. Most notably, this includes examples where $\monA\monB < 0$ and $\comA\comB < 0$, for which neither the primal nor the dual nor the primal-dual inclusion are monotone (see e.g. \Cref{ex:CP:example5:revisited}).

As first observed in \cite{shefi2014rate}, \ref{eq:CP} can be viewed as a particular instance of proximal ADMM. Exploiting this connection it is possible to obtain convergence results for \ref{eq:CP} based on those for nonconvex proximal ADMM, see \cite{li2015Global,bot2020proximal}. 
This approach leads to requirements for \(L\) such as a full row rank assumption and restrictions on its condition number (see \cite[Ass. 1, Rem. 2(c)]{bot2020proximal}).
Our convergence results for \ref{eq:CP} also do not depend on any explicit rank conditions on $L$, allowing to cover rank-deficient cases without introducing a switching mechanism.

  \subsection{Contributions}

The main contribution of the paper is to establish convergence of \ref{eq:CP} under the assumption that the primal-dual operator $\Tpd$ admits a set of $\DRSRho$\hyp{}oblique weak Minty solutions, which leads to novel stepsize and relaxation parameter ranges in both strongly monotone and nonmonotone settings (see \Cref{thm:CP:full} and the preceding discussion). Interestingly, in contrast to the classical stepsize condition $\gamma\tau\nrm{L}^2 \leq 1$ in the monotone setting, the conditions obtained through our analysis not only depend on the norm of $L$ but also on its other singular values. The tightness of our main convergence theorem is demonstrated through \Cref{ex:CP:example5,ex:CP:second}.

As our second main contribution, convergence results are provided for the class of semimonotone operators \cite[Sec. 4]{evens2023convergence}, which can be viewed as a natural extension of the (hypo)- and co(hypo)monotone operators.
We show that the stepsize requirements reduce to a look-up table depending on the level of (hypo)- and co(hypo)monotonicity (see \cref{cor:CP:semi}).
These results are made possible by establishing a link between the oblique weak Minty assumption for the primal-dual operator and semimonotonicity of the underlying operators $A$ and $B$, relying on the extended calculus rules developed in \cref{sec:semi} (see also \Cref{thm:calculus:primaldual}).

  \subsection{Organization}

The paper is structured in the following manner.
In \Cref{subsec:notation}, some notation and standard definitions are provided.
\Cref{sec:preliminaries} recalls the main convergence results from \cite{evens2023convergence} for \ref{eq:PPPA-intro} in the nonmonotone setting.
In \Cref{sec:CP}, the primal-dual equivalence between \ref{eq:CP} and \ref{eq:PPPA-intro} is established, which lead to convergence of \ref{eq:CP} under an oblique weak Minty assumption on the associated primal-dual operator.
\Cref{sec:semi} discusses and introduces various calculus rules for the class of $(\Mon, \Com)$\hyp{}semimonotone operators. Leveraging these calculus rules, \Cref{sec:CP:semi} presents a set of sufficient conditions for the convergence of \ref{eq:CP}, based on the semimonotonicity of the underlying operators, along with several examples.
Finally, \Cref{sec:conclusion} concludes the paper.
Several proofs and auxiliary results are deferred to the Appendix.

  \subsection{Notation}\label{subsec:notation}

The set of natural numbers including zero is denoted by \(\N\coloneqq\set{0,1,\hdots}\).
The set of real and extended-real numbers are denoted by \(\R\coloneqq(-\infty,\infty)\) and \(\Rinf\coloneqq\R\cup\set\infty\), while the positive and strictly positive reals are \(\R_+\coloneqq[0,\infty)\) and \(\R_{++}\coloneqq(0,\infty)\).
We use the notation $\seq{w^k}[k\in I]$ to denote a sequence with indices in the set $I\subseteq \N$. When dealing with scalar sequences we use the subscript notation $\seq{\gamma_k}[k\in I]$. 
We denote the positive part of a real number by $[\cdot]_{+} \dfn \max\{0, \cdot\}$ and the negative part by 
$
    [\cdot]_{-} \dfn -\min\{0, \cdot\}
$.
With \(\id\) we indicate the identity function \(x\mapsto x\) defined on a suitable space.
The identity matrix is denoted by $\I_n\in\R^{n\times n}$ and the zero matrix by $0_{m\times n}\in\R^{m\times n}$; we write respectively $\I$ and $0$ when no ambiguity occurs. 
Adopting the notation from \cite{bernstein2009matrix}, we say a matrix $P \in \R^{m \times n}$ is empty if $\min(m,n) = 0$ and use the conventions
$P \mkern 1mu 0_{n\times 0} = 0_{m\times 0}$,
$0_{0\times m} P = 0_{0\times n}$
and
$0_{m\times0}0_{0\times n} = 0_{m\times n}$.
Moreover, we use the convention that $\tfrac10 = \infty$ and that $0 \cdot \infty = 0$.
Given a matrix $P \in \R^{m \times n}$, we denote the range of $P$ by $\range{P}$ and the kernel of $P$ by $\ker{P}$.
The trace of a square matrix $P \in \R^{n\times n}$ is denoted by $\trace P$.

We denote by $\R^n$ the standard $n$-dimensional Euclidean space with inner product $\langle\cdot,\cdot\rangle$ and induced norm $\|\cdot\|$. 
The set of symmetric $n$-by-$n$ matrices is denoted by $\sym{n}$. Given a symmetric matrix $P \in \sym{n}$, we write $P \succeq 0$ and $P \succ 0$ to denote that $P$ is positive semidefinite and positive definite, respectively.
Furthermore, for any $P \in \sym{n}$ we define the 
quadratic form \(\qsemi{x}{P} \dfn \inner{x, Px}\).
Let $\diag(\cdot)$ denote the diagonal
matrix whose arguments constitute its diagonal elements.
For arbitrary matrices $A$ and $B$, we define the direct sum $A \oplus B = \blkdiag(A,B)$,
where
$\blkdiag(\cdot)$ denotes the block diagonal
matrix whose arguments constitute its diagonal blocks.
We denote the Kronecker product between two matrices of arbitrary size by \(\otimes\).

Two vectors $u, v \in \R^n$ are said to be orthogonal if $\inner{u,v} = 0$, and orthonormal if they are orthogonal and $\nrm{u} = \nrm{v} = 1$. Two linear subspaces $\mathbf{U} \subseteq \R^n$ and $\mathbf{V} \subseteq \R^n$ are said to be orthogonal if any $u \in \mathbf{U}$ and any $v \in \mathbf{V}$ are orthogonal. We say that $U \in \R^{n \times m}$ is an orthonormal basis for a linear subspace $\mathbf{U} \subseteq \R^n$ if $U$ has orthonormal columns and $\range{U} = \mathbf{U}$.
 

The effective domain of an extended-real-valued function $f : \R^n \rightarrow \Rinf$ is given by the set $\dom f \dfn \set{x \in \R^n}[f(x) < \infty]$. We say that $f$ is proper if $\dom f \neq \emptyset$ and that $f$ is lower semicontinuous (lsc) if the epigraph $\epi f \dfn \set{(x,\alpha) \in \R^n \times \R}[f(x) \leq \alpha]$ is a closed subset of $\R^{n+1}$. We denote the limiting subdifferential of $f$ by $\partial f$.
We denote the normal cone of a set $E\subseteq\R^n$ by \(N_E\) and 
the projection onto $E$ is denoted by 
    \(
        \proj_E(x)\coloneqq{}\argmin_{z\in E}\|z-x\|.
    \)
An operator or set-valued mapping $A:\R^n\rightrightarrows\R^d$ maps each point $x\in\R^n$ to a subset $A(x)$ of $\R^d$. We will use the notation $A(x)$ and $Ax$ interchangeably. 
 We denote the domain of $A$ by $\dom A\coloneqq\{x\in\R^n\mid Ax\neq\emptyset\}$,
its graph by $\graph A\coloneqq\{(x,y)\in\R^n\times \R^d\mid y\in Ax\}$, and 
the set of its zeros by $\zer A\coloneqq\{x\in\R^n \mid 0\in Ax\}$. 
The inverse of $A$ is defined through its graph: $\graph A^{-1}\coloneqq\{(y,x)\mid (x,y)\in\graph A\}$.
The \emph{resolvent} of $A$ is defined by $J_A\coloneqq(\id+A)^{-1}$.
We say that \(A\) is \emph{outer semicontinuous (osc)} at $\other{x}\in\dom A$ if 
\ifjota
\else
    \begin{equation*}
        \limsup_{x\to \other{x}} Ax \coloneqq \{y\mid \exists x^k \to \other{x}, \exists y^k \to y \textrm{ with } y^k\in Ax^k\} \subseteq A\other{x}.
    \end{equation*}
\fi
\ifjota
\else
    Outer semicontinuity of $A$ everywhere is equivalent to its graph being a closed subset of $\R^n\times \R^d$.
\fi
\ifjota
\else
    \begin{definition}[(co)monotonicity]
        An operator $A:\R^n\rightrightarrows\R^n$ is said to be $\mon$-monotone for some $\mon\in \R$ if
        \[
            \langle x-\other{x},y-\other{y}\rangle \geq \mon\|x-\other{x}\|^2, \qquad \text{for all $(x,y),(\other{x},\other{y})\in\graph A$},
        \]
        and it is said to be $\com$-comonotone for some $\com\in \R$ if
        \[
            \langle x-\other{x},y-\other{y}\rangle \geq \com\|y-\other{y}\|^2, \qquad \text{for all $(x,y),(\other{x},\other{y})\in\graph A$}.
        \]
        $A$ is said to be maximally 
        (co-)monotone if its graph is not strictly contained in the graph of another 
        (co-)monotone operator. We say that $A$ is monotone if it is $0$-monotone.
    \end{definition}
\fi

\section{Preliminaries on the preconditioned proximal point method}
\label{sec:preliminaries}
Departing from the classical monotone setting, convergence of relaxed \ref{eq:PPPA-intro} was established in \cite{evens2023convergence} for a class of nonmonotone operators that admit a set of \emph{oblique weak Minty solutions} (see \cref{def:SWMVI}). This result will serve as our primary tool for establishing convergence of \ref{eq:CP} in the nonmonotone setting, which is why we will reiterate it here.
In particular, its analysis involves the following assumptions.

\begin{assumption} \label{ass:PPPA}
The operator \(T\) in \eqref{prob:P1} and the symmetric positive semidefinite preconditioner \(\M\) in \eqref{eq:PPPA-intro} satisfy the following properties. 
	\begin{enumeratass}
		\item \label{ass:PPPA:0} \(T:\R^n\rightrightarrows\R^n\) is outer semicontinuous.
		\item \label{ass:PPPA:0.5} The preconditioned resolvent \((\M + T)^{-1}\circ\M\) has full domain. 
		\item \label{ass:PPPA:1} 
		There exists a nonempty set \(\pazocal{S}^\star\subseteq \zer T\) and a symmetric, possibly indefinite matrix \(V \in \sym{n}\) such that \(T\) has $\DRSRho$\hyp{}oblique weak Minty solutions at \(\pazocal{S}^\star\). 

		\item \label{ass:PPPA:2} \(\M \in \sym{n}\) is a symmetric positive semidefinite matrix
		such that
	    \begin{equation}\label{eq:PPPA:eigcond}
            1 + \etamin > 0
            \quad\text{with}\quad
			\etamin
				{}\dfn{}
			\lambda_{\mathrm{min}}(\tp \basisP \DRSRho \M \basisP),
                \ifjota\else
                \footnote{As $\tp \basisP \DRSRho \M \basisP$ is similar to a symmetric matric, its eigenvalues are real \cite[Rem. 2.2]{evens2023convergence}.}
                \fi
	    \end{equation}
		where \(\basisP\) is any orthonormal basis for the range of \(\M\).
	\end{enumeratass}
\end{assumption}
As observed in \cite[Rem. 2.2]{evens2023convergence}, \cref{ass:PPPA:2} holds if and only if
\(
    (U^\top \M U)^{-1} + U^\top \DRSRho U \succ 0
\),
which is for instance trivially satisfied when $\DRSRho$ is positive semidefinite.
In contrast to the convergence analysis techniques relying on firm nonexpansiveness of the resolvent mapping, the analysis of \cite{evens2023convergence} relies on a projective interpretation of the preconditioned proximal point algorithm, which dates back to \cite{solodov1996Modified,solodov1999hybrid,konnov1997Class}.
Most notably, it was demonstrated in \cite[Lem. 2.3]{evens2023convergence} that the update rule for the (shadow) sequence generated by \ref{eq:PPPA-intro} can be interpreted as a relaxed projection onto a certain halfspace, and that if any iterate belongs to this halfspace, which contains the set of projected oblique weak Minty solutions \(\projQ\pazocal{S}^\star\), this implies its optimality.
Based on this insight, the following convergence result for \ref{eq:PPPA-intro} was established.

\begin{theorem}[\protect{\cite[Thm. 2.4]{evens2023convergence}} convergence of \ref{eq:PPPA-intro}]\label{thm:pppa}
    Suppose that \cref{ass:PPPA} holds, and consider a sequence $\seq{z^k, \bar z^k}$ generated by \ref{eq:PPPA-intro} starting from $z^0 \in \R^n$ with relaxation parameters $\lambda_k \in \bigl(0,2(1 + \etamin)\bigr)$ such that $\liminf_{k\to\infty}\lambda_k\bigl(2(1+\etamin)-\lambda_k\bigr)>0$, where \(\etamin\) is defined as in \eqref{eq:PPPA:eigcond}.
    Then, either a point $\bar z^k \in \zer T$ is reached in a finite number of iterations or 
    the following hold for the sequence $\seq{z^k, \bar z^k}$.
    \begin{enumerate}
        \item\label{it:pppa:v} $\bar{v}^k:=\M(z^k-\bar{z}^k)\in T\bar{z}^k$ for all \(k\) and $(\bar{v}^k)_{k\in N}$ converges to zero.
        \item\label{it:pppa:subseq} Every limit point (if any) of $\seq{\bar z^k}$ belongs to $\zer T$.
        \item\label{it:pppa:bounded} The sequences $(\projQ z^k)_{k\in\N}$, $(\projQ\bar{z}^k)_{k\in\N}$ are bounded and their limit points belong to $\projQ\zer T$.
       \item\label{it:pppa:full} 
        If in \cref{ass:PPPA:1} the set
        $
            \projQ 
            \pazocal{S}^\star
        $
        is equal to
        $
            \projQ 
            \zer T
        $, then the sequences $\seq{\projQ z^k}$, $\seq{\projQ \bar z^k}$
        converge to some element of $\projQ\zer T$.
        If additionally $(\M + T)^{-1}\circ\M$ is (single-valued) continuous,
        then $\seq{\bar{z}^k}$ converges to some 
        $
            z^\star \in \zer T
        $.
        Finally, if $\lambda_k$ is additionally uniformly bounded in the interval $(0,2)$, then $\seq{z^k}$ also converges to 
        $
            z^\star \in \zer T
        $.
    \end{enumerate}
\end{theorem}
The statement of \Cref{it:pppa:full} has been slightly extended compared to \cite[Thm. 2.4]{evens2023convergence}, as the convergence of $\seq{\projQ z^k}$ to some element of $\projQ\zer T$ also implies the convergence of $\seq{\projQ \bar z^k}$ to the same point, since \(\seq{\bar{v}^k}\) converges to zero.
We also remark that by \cref{it:pppa:bounded} the boundedness of \(\seq{\z^k}\) is already guaranteed if the preconditioned resolvent is locally bounded in the sense of \cite[Def.~5.14]{rockafellar2009Variational}.

As discussed in \cite{evens2023convergence}, \Cref{ass:PPPA:1} can be further relaxed by only requiring $T$ to have $\DRSRho$\hyp{}oblique weak Minty solutions at \(\pazocal{S}^\star\) on \(\pazocal{U} = (\range{(\M+T)^{-1}\M} \times \range{\M})\) instead.
Under this relaxed assumption, all results from \Cref{thm:pppa} remain valid, as the proof of \Cref{thm:pppa} only involves invoking \eqref{def:WMVI} at points in this restricted set.
This relaxation will prove to be relevant in \Cref{ex:CP:example5}.

\section{Chambolle--Pock under oblique weak Minty}\label{sec:CP}

In the monotone setting, it is well-known that \ref{eq:CP} can be interpreted as applying \ref{eq:PPPA-intro} to the primal-dual operator $\Tpd$ \cite{esser2010general,he2012Convergence,condat2013primal,latafat2017Asymmetric,bredies2022degenerate}. 
Relying upon the abstract duality framework from \cite{attouch1996general}, \cite[Sec. 6.9]{auslender2006asymptotic}, this equivalence can be extended to the nonmonotone setting.
Within this framework, inclusion problems \eqref{prob:composite} and \eqref{eq:primaldual} are labelled as the \emph{primal} and the \emph{primal-dual inclusion}, respectively. Related to these two inclusions is the \emph{dual inclusion}, given by
\begin{equation}\label{eq:dual}\tag{D-I}
    \text{find} \quad y \in \R^n \quad \text{such that} \quad 
    0\in \Td y 
        {}\coloneqq{}
    (-L)A^{-1}(-L^\top)(y)+B^{-1}(y),
\end{equation}
A fundamental equivalence property for these inclusions is summarized below.

\begin{proposition}[\protect{\cite[Prop. 6.9.2]{auslender2006asymptotic}}]\label{prop:auslender}
    Let $(x, y) \in \R^n \times \R^m$. The following statements are 
    \ifjota
\else
        equivalent.
        \begin{enumerate}
            \item $(x,y) \in \zer \Tpd$.
            \item $x \in \zer \Tp$ and $y \in \zer \Td$.
            \item $(x,- L^\top y)\in\graph A$ and $(L x,y)\in\graph B$.
        \end{enumerate}
    \fi
    Furthermore, it holds that \(\zer\Tp = \{x \mid \exists y : (x,y) \in\zer \Tpd\}\) and \(\zer\Td = \{y \mid \exists x : (x,y) \in\zer \Tpd\}\).
\end{proposition}
A solution of the primal inclusion \eqref{prob:composite} (and of the dual inclusion \eqref{eq:dual}) can thus be obtained by finding a solution of the associated primal-dual inclusion. 
Now, consider applying \ref{eq:PPPA-intro} to the primal-dual inclusion \eqref{eq:primaldual}, with the preconditioner $\M$ given by \eqref{eq:CP:P}. Then, each iteration corresponds to first finding $\bar x^k$ and $\bar y^k$ satisfying the inclusions
\ifjota
\else
    \begin{align*}
        \tfrac1\gamma x^k - L^\top y^k 
            {}\in{}
        (\tfrac1\gamma \id + A)\x^k
        \quad\text{and}\quad
        \tfrac1\tau y^k - L x^k 
            {}\in{}
        - 2L \x^k + (\tfrac1\tau \id + B^{-1})\y^k
    \end{align*}
\fi
and then performing a relaxation step \(z^{k+1} = z^k + \lambda_k(\z^k - z^k)\). Multiplying the two relations by \(\gamma\) and \(\tau\), respectively, and reordering the terms, the update rule for \ref{eq:CP} is retrieved. This result is summarized in the following lemma.

\begin{lemma}[equivalence of \ref{eq:CP} and \ref{eq:PPPA-intro}]\label{lem:equivalence:PPPA-PDHG}
    Let $z^0 = (x^0, y^0) \in \R^{n + m}$ be the initial guess for \ref{eq:CP} and for \ref{eq:PPPA-intro} applied to the primal-dual inclusion \eqref{eq:primaldual}, with the preconditioner $\M$ given by \eqref{eq:CP:P}.
    Then, the sequences \(\seq{z^k} = \seq{x^k, y^k}\), \(\seq{\bar z^k} = \seq{\bar x^k, \bar y^k}\) generated by \ref{eq:CP} 
    satisfy update rule \ref{eq:PPPA-intro}.
\end{lemma}

Leveraging this connection, we will establish the convergence of \ref{eq:CP} based on \cref{thm:pppa} for \ref{eq:PPPA-intro}.
In contrast to the classical stepsize condition $\gamma\tau\nrm{L}^2 \leq 1$ in the monotone setting, our upcoming analysis will demonstrate that the stepsize condition on $\gamma$ and $\tau$ for \ref{eq:CP} in general does not only depend on $\nrm{L}$, i.e., the largest singular value of $L$, but also its other singular values. 
Therefore, let \(r\) denote the rank of \(L\), and without loss of generality, let \(\sigma_1,\hdots,\sigma_d\) denote its distinct strictly positive singular values in descending order with respective multiplicities \(m_1,\hdots,m_d\). Then, it holds that $r = \sum_{i=1}^d m_i$. Define $\SL = \sigma_1 \I_{m_1} \oplus \cdots \oplus \sigma_d \I_{m_d} \in \R^{r \times r}$ and consider the singular value decomposition 
\begin{equation}\label{eq:L:svd}
    L = 
    \ifjota
\else
        \begin{bmatrix}
            \UL & \UL^\prime
        \end{bmatrix}
        \begin{bmatrix}
            \SL & \\
            & 0
        \end{bmatrix}
        \begin{bmatrix}
            \tp \VL \\
            {\VL^\prime}^\top
        \end{bmatrix},
    \fi
    \;\;
    \UL_
        {}={}
    \begin{bmatrix}
        \UL_{1} & \cdots & \UL_{d}
    \end{bmatrix},
    \;\;
    \VL
        {}={}
    \begin{bmatrix}
        \VL_{1} & \cdots & \VL_{d}
    \end{bmatrix},
\end{equation}
where the zero matrix is in $\R^{(m-r)\times(n-r)}$,
\(\UL_i \in \R^{m \times m_i}\) and $\VL_i \in \R^{n \times m_i}$, $i \in [d]$, have orthonormal columns that span the eigenspace corresponding to eigenvalue $\sigma_i^2$ of $L\tp L$ and $\tp L L$, respectively,
and $\UL^\prime \in \R^{m \times (m-r)}$ and $\VL^\prime \in \R^{n \times (n-r)}$ have orthonormal columns which span the null space of $L^\top$ and $L$, respectively.
The projection onto the range and the kernel of $L$ and $L^\top$ can be expressed as \cite[Sec. 2.5.2]{golub2013matrix}
\ifjota
\else
    \begin{equation}\label{eq:projXY}
        \proj_{\range{L}} 
            = 
        \UL \UL^\top,
        \quad
        \proj_{\range{\tp L}} 
            =
        \VL \VL^\top,
        \quad
        \proj_{\ker{L}}
            =
        \I_n - \VL \VL^\top = \VL' \tp {\VL'}
        \quad\text{and}\quad
        \proj_{\ker{\tp L}}
            =
        \I_m - \UL \UL^\top = \UL' \tp {\UL'}.
    \end{equation}
\fi
These projections will play a central role in our upcoming analysis. We will work under the following assumptions on the individual operators $\A$ and $\B$ and the (nonzero) matrix $L$.

\begin{assumption} \label{ass:CP}
In problem \eqref{prob:composite}, the following hold.
    \begin{enumeratass}
        \item \label{ass:CP:1} Operators \(A\) and \(B\) are outer semicontinuous.
        \item \label{ass:CP:2} For the selected positive stepsizes the corresponding resolvents have full domain, i.e., 
        $\dom J_{\gamma A} = \R^n$ and $\dom J_{\tau B^{-1}} = \R^m$.
        \item \label{ass:SWMVICP} 
        The set \(\zer \Tpd\) is nonempty and there exist parameters $\DRSrhocom, \DRSrhocom^\prime, \DRSrhomon, \DRSrhomon^\prime \in \R$ and a nonempty set \(\pazocal{S}^\star \subseteq \zer \Tpd\) such that 
        the primal-dual operator $\Tpd$ has $\DRSRho$\hyp{}oblique weak Minty solutions at \(\pazocal{S}^\star\), where
        \begin{equation}\label{eq:WMSDRS}
            \DRSRho
                {}\dfn{}
            \DRSRho_{\rm P} \oplus \DRSRho_{\rm D}
                {}={}
            \Bigl(
                \DRSrhocom \proj_{\range{\tp L}} + \DRSrhocom^\prime \proj_{\ker{L}}
            \Bigr)
            \oplus
            \Bigl(
                \DRSrhomon \proj_{\range{L}} + \DRSrhomon^\prime \proj_{\ker{\tp L}} 
            \Bigr) \in \mathbb{S}^{n+m}.
        \end{equation}
        Moreover, defining
        \begin{align*}
            \delta 
                {}\dfn{}&
            1 
                -
            [\DRSrhocom\DRSrhomon]_- \bigl(\nrm{L}^2 - \sigma_{d}^2\bigr)
                {}={}
            \begin{cases}
                1, &\;\ \text{if } \DRSrhocom\DRSrhomon \geq 0,\\
                1 + \DRSrhocom\DRSrhomon \bigl(\nrm{L}^2 - \sigma_{d}^2\bigr) &\;\ \text{if } \DRSrhocom\DRSrhomon < 0,
            \end{cases}
            \numberthis\label{eq:def:delta}
        \end{align*}
        the following conditions on $\DRSrhocom, \DRSrhocom^\prime, \DRSrhomon$ and $\DRSrhomon^\prime$ hold.
        \ifjota
\else
            \begin{align*}
                [\DRSrhocom]_-[\DRSrhomon]_- < \tfrac{1}{4\nrm{L}^2},
                    \quad
                \DRSrhocom^\prime
                    {}\geq{}
                \begin{cases}
                    0, &\quad \text{if } \DRSrhocom \geq 0,\\
                    {\tfrac{\DRSrhocom}{\delta - \DRSrhocom\DRSrhomon\nrm{L}^2}}, &\quad \text{if } \DRSrhocom < 0,
                \end{cases}
                \quad\text{and}\quad
                \DRSrhomon^\prime
                    {}\geq{}
                \begin{cases}
                    0, &\quad \text{if } \DRSrhomon \geq 0,\\
                    {\tfrac{\DRSrhomon}{\delta - \DRSrhocom\DRSrhomon\nrm{L}^2}}, &\quad \text{if } \DRSrhomon < 0.
                \end{cases}
                \numberthis\label{eq:existence:condition}
            \end{align*}
        \fi
    \end{enumeratass}
\end{assumption}
In \Cref{ass:SWMVICP}, the matrix $\DRSRho$ is defined as in \eqref{eq:WMSDRS} with the understanding that in our convergence analysis we shall set
$
    \DRSrhocom^\prime = \infty
$
when $L$ is full column rank and
$
    \DRSrhomon^\prime = \infty
$
when $L$ is full row rank, as in these cases $\DRSRho$ is not affected by their choice.
Note also that the conditions on $\DRSrhocom^\prime$ and $\DRSrhomon^\prime$ specified in \eqref{eq:existence:condition} can be slightly relaxed (see  
\ifjota
\else
    \Cref{thm:CP:full:strong}).
\fi 
This relaxation, however, results in more complicated stepsize rules for $\gamma$ and $\tau$, and for this reason is not pursued further. 

The imposed block structure on $\DRSRho$ in \eqref{eq:WMSDRS} is not simply an arbitrary choice, but it aligns perfectly with the inherent structure present within the primal-dual operator itself. 
To illustrate this, consider the following lemma, which shows that the blocks $\DRSRho_{\rm P}$ and $\DRSRho_{\rm D}$ from \eqref{eq:WMSDRS} can be interpreted as the primal and the dual blocks of $\DRSRho$, respectively. 
This lemma extends \cite[Lem. 3.2]{evens2023convergence}, which considers the case $L = \I$.
Notably, the converse of this lemma is not true, as demonstrated in the discussion below \cite[Lem. 3.2]{evens2023convergence} through a simple linear example.
The proof is deferred to \Cref{proof:lem:SWMI:primal-dual}.
\begin{lemma}[oblique weak Minty for primal and dual operator]\label{lem:SWMI:primal-dual}
    Let $\DRSRho_{\rm P} \in \sym{n}$ and $\DRSRho_{\rm D} \in \sym{m}$.
    Suppose that
    there exists a nonempty set \(\pazocal{S}^\star \subseteq \zer \Tpd\) such that 
    the primal-dual operator $\Tpd$ has $(\DRSRho_{\rm P} \oplus \DRSRho_{\rm D})$\hyp{}oblique weak Minty solutions at \(\pazocal{S}^\star\)
    and let
    \ifjota
\else
        \[
            \pazocal{S}^\star_P \dfn \{x^\star \mid \exists y^\star : (x^\star,y^\star) \in \pazocal{S}^\star\} \subseteq \zer \Tp
            \;\; \text{and} \;\; 
            \pazocal{S}^\star_D \dfn \{y^\star \mid \exists x^\star : (x^\star,y^\star) \in \pazocal{S}^\star\}\subseteq \zer \Td.
        \]
    \fi
    Then, the primal operator $\Tp$ has $\DRSRho_{\rm P}$\hyp{}oblique weak Minty solutions at $\pazocal{S}^\star_P$ and the dual operator $\Td$ has $\DRSRho_{\rm D}$\hyp{}oblique weak Minty solutions at $\pazocal{S}^\star_D$.
\end{lemma}
\ifjota
\else
    As shown in the proof of \Cref{lem:SWMI:primal-dual}, the quadratic terms  $\qindef{y_\A + L^\top y_\B}{\DRSRho_{\rm P}}$ and $\qindef{x_\B - L x_\A}{\DRSRho_{\rm D}}$ emerging in oblique weak Minty inequality correspond to the primal and the dual problems, respectively. By selecting $\DRSRho_{\rm P}$ and $\DRSRho_{\rm D}$ as in \eqref{eq:WMSDRS}, these terms can be written as 
    \begin{align*}
        \qindef{y_\A + L^\top y_\B}{\DRSRho_{\rm P}}
            {}={}&
        \DRSrhocom\nrm{\proj_{\range{\tp L}}y_\A + L^\top y_\B}^2 + \DRSrhocom^\prime\nrm{\proj_{\ker{L}}y_\A}^2,\\
        \qindef{x_\B - L x_\A}{\DRSRho_{\rm D}}
            {}={}&
        \DRSrhomon\nrm{\proj_{\range{L}} x_\B - L x_\A}^2 + \DRSrhomon^\prime\nrm{\proj_{\ker{\tp L}} x_\B}^2,
    \end{align*}
    reducing to the norm of the scaled sum of a vector belonging to the range of \(L^\top\) and another to its nullspace (resp., range of \(L\) and nullspace of \(L^\top\)). 
    This decomposition proves essential in the proof of \Cref{thm:CP:full}, as it enables to split condition
    \eqref{eq:PPPA:eigcond} into two terms, one depending only on $\DRSrhocom$ and $\DRSrhomon$ and the other only depending on $\DRSrhocom^\prime$ and $\DRSrhomon^\prime$ (see \eqref{eq:proof:CP:condition:blocked}).
\fi

One of the main aspects of the upcoming convergence proof for \ref{eq:CP} is showing that \Cref{ass:PPPA} holds for the operator $\Tpd$ and preconditioner $\M$ from \eqref{eq:CP:P}.
\ifjota
\else
    To this end, the stepsizes $\gamma$ and $\tau$ and the relaxation parameters $\lambda_k$ need to be selected as follows.
\fi

\begin{stepsize}\label{ass:CP:stepsize:rule}
    The stepsizes $\gamma$ and $\tau$ satisfy the bounds provided in \Cref{tab:CP:stepsize:rule}, where $\delta$ is defined as in \eqref{eq:def:delta} and
    \ifjota
\else
        \begin{equation}
            \begin{aligned}
                \gamma_{\rm min} 
                    {}\dfn{}
                \frac{-2\DRSrhocom}{\delta + \sqrt{\delta^2 - 4\DRSrhocom\DRSrhomon\nrm{L}^2}},
                \;\;
                \gamma_{\rm max}
                    {}\dfn{}
                \frac{\delta + \sqrt{\delta^2 - 4\DRSrhocom\DRSrhomon\nrm{L}^2}}{-2\DRSrhomon\nrm{L}^2}
                \;\;\text{and}\;\;
                \tau_{\rm min}(\gamma)
                    {}\dfn{}
                \frac{[-\DRSrhomon(\gamma + \DRSrhocom)]_+}{\gamma(\delta-\DRSrhocom\DRSrhomon\nrm{L}^2) + \DRSrhocom}.
            \end{aligned}
            \numberthis\label{eq:stepsizes:CP:bounds:gam}
        \end{equation}
    \fi
    \ifjota
\else
        \vspace{-0.4cm}
    \fi
    \begin{figure}[H]
        \centering
        \includetikz{Tables/stepsizes-CP}
        \caption{Range of the stepsizes $\gamma$ and $\tau$ for \ref{eq:CP}.}
        \label{tab:CP:stepsize:rule}
    \end{figure}
\end{stepsize}

\newlength{\maxmin}
\setlength{\maxmin}{\widthof{$(\|L\|)$}-\widthof{$(\sigma_{2})$}}

\newlength{\maxmincond}
\setlength{\maxmincond}{\widthof{$\text{if } \min\{\DRSrhocom, \DRSrhomon\} > 0$}-\widthof{$\text{if } \DRSrhocom\DRSrhomon \geq 0$}}

\begin{relaxation}\label{ass:CP:relaxation:rule}
    Define
    \begin{align*}
        \theta_{\gamma \tau} (\sigma) 
            {}\coloneqq{}
        \sqrt{\left(\tfrac{1}{2\gamma}\DRSrhocom - \tfrac{1}{2\tau}\DRSrhomon\right)^2 + \DRSrhocom\DRSrhomon\sigma^2}.\numberthis\label{eq:theta}
    \end{align*}
    Let $\etamin \dfn \min\{\eta, \eta^\prime\}$, where $\eta^\prime$ is defined as in \Cref{tab:CP:degen} and
    \begin{align*}
        \eta
            {}\dfn{}
        \begin{cases}
            \left\{
            \begin{array}{ll}
                \tfrac{1}{2\gamma}\DRSrhocom + \tfrac{1}{2\tau}\DRSrhomon - \theta_{\gamma \tau}(\sigma_{d}),
                &\;\;\, \text{if } \DRSrhocom\DRSrhomon < 0\\
                \tfrac{1}{2\gamma}\DRSrhocom + \tfrac{1}{2\tau}\DRSrhomon - \theta_{\gamma \tau}(\|L\|),
                &\;\;\, \text{if } \DRSrhocom\DRSrhomon \geq 0\hspace*{\maxmincond}\\
            \end{array}
            \right\}
            & \;\;\, \text{if } \gamma\tau < \tfrac{1}{\nrm{L}^2},\\
            \tfrac{1}{\gamma}\DRSrhocom + \tfrac{1}{\tau}\DRSrhomon, & \;\;\, \text{if } \gamma\tau = \tfrac{1}{\nrm{L}^2} \text{ and } d=1,\\
            \left\{
            \begin{array}{ll}
                \tfrac{1}{2\gamma}\DRSrhocom + \tfrac{1}{2\tau}\DRSrhomon - \theta_{\gamma \tau}(\sigma_{d})
                ,
                &\;\;\,\hspace{\maxmin} \text{if } \DRSrhocom\DRSrhomon < 0\\
                \tfrac{1}{2\gamma}\DRSrhocom + \tfrac{1}{2\tau}\DRSrhomon - \theta_{\gamma \tau}(\sigma_{2})
                ,
                &\;\;\,\hspace{\maxmin} \text{if } \min\{\DRSrhocom, \DRSrhomon\} \geq 0\\
                \tfrac{1}{\gamma}\DRSrhocom + \tfrac{1}{\tau}\DRSrhomon,
                &\;\;\,\hspace{\maxmin} \text{otherwise}
            \end{array}
            \right\}
            & \;\;\, \text{if } \gamma\tau = \tfrac{1}{\nrm{L}^2} \text{ and } d>1.
        \end{cases}\numberthis\label{eq:thm:CP:etamin}
    \end{align*} 
    The sequence $\seq{\lambda_k}$ satisfies $\lambda_k\in\bigl(0,2(1 + \etamin)\bigr)$ and $\liminf_{k\to\infty}\lambda_k\bigl(2(1 + \etamin)-\lambda_k\bigr)>0$.
    \begin{figure}[H]
        \centering
        \includetikz{Tables/relaxation-CP}
        \caption{Definition of $\eta^\prime$ in \cref{ass:CP:relaxation:rule}.}
        \label{tab:CP:degen}
    \end{figure}
\end{relaxation}

When $L = \I$ and $\tau = \nicefrac1\gamma$, \Cref{ass:CP:stepsize:rule} matches the stepsize range for Douglas--Rachford splitting from \cite[Thm. 3.3]{evens2023convergence}.
In \Cref{ass:CP:relaxation:rule} there is an interplay between the stepsizes $\gamma$ and $\tau$ and the range of admissible relaxation parameters $\lambda_k$.
In the monotone setting where $\DRSrhocom = \DRSrhocom^\prime = \DRSrhomon = \DRSrhomon^\prime = 0$, this interplay vanishes as the relaxation rule reduces to the classical condition $\lambda_k \in (0, 2)$.
In strongly monotone settings, this interplay allows us to select relaxation parameters beyond the classical upper bound of two.
For instance, when $\DRSrhocom > 0$ and $\DRSrhocom^\prime > 0$, then for small enough stepsizes $\gamma$ the upper bound on $\lambda_k$ will be larger than two (see \cite[Ex. 6.2]{evens2023convergence} for an example in the DRS setting where $L = \I$ and $\tau = \nicefrac1\gamma$).
Conversely, when $\DRSrhomon > 0$ and $\DRSrhomon^\prime > 0$, this upper bound will be larger than two for small enough $\tau$.
Finally, when all $\beta$ parameters are strictly positive this holds for all valid stepsizes $\gamma$ and $\tau$ (see e.g. \Cref{ex:CP:second}).

Having discussed our underlying assumptions and stepsize/relaxation parameter rules, we will now present our main convergence theorem for \ref{eq:CP}. 
The proof relies on carefully decomposing both the preconditioner $\M$ and the oblique weak Minty matrix $\DRSRho$ into two separate, orthogonal matrices. Exploiting the inherent structure present in these orthogonal matrices, the conditions from \Cref{ass:PPPA} are reduced to a set of eigenvalue problems of two\hyp{}by\hyp{}two matrices (see equation \eqref{eq:proof:CP:eigcond:reformulated:blocks}).
{
\def\Pr{{S}}
\def\Pn{{\M^\prime}}
\newcommand\Prs[1]{{\hat{\M}_{#1}}}
\def\PrS{{\hat{\M}}}

\def\Vr{{\DRSRho^\star}}
\def\Vn{{\DRSRho^\prime}}
\newcommand\Vrs[1]{{\hat{\DRSRho}_{#1}}}
\def\VrS{{\hat{\DRSRho}}}

\def\Ur{{\basisP^\star}}
\def\Un{{Z^\prime}}
\newcommand\Urs[1]{{\hat{\basisP}_{#1}}}
\def\UrS{{\hat{\basisP}}}

\begin{theorem}\label{thm:CP:full}
    Suppose that \cref{ass:CP} holds,
    that $\gamma$ and $\tau$ are selected according to \Cref{ass:CP:stepsize:rule}
    and that the relaxation sequence $\seq{\lambda_k}$ is selected according to \Cref{ass:CP:relaxation:rule}.
    Consider the sequences $\seq{z^k} = \seq{x^k,y^k}$ and $\seq{\z^k} = \seq{\x^k,\y^k}$ generated by \ref{eq:CP} starting from $z^0 \in \R^{n+m}$. 
    Then, either a point $\bar z^k \in \zer \Tpd$ is reached in a finite number of iterations or the following hold for the sequences $\seq{z^k}$ and $\seq{\z^k}$.
    \begin{enumerate}
        \item\label{it:CP:v} 
        $\bar{v}^k:=\M(z^k-\z^k)\in \Tpd\z^k$ for all \(k\) and $(\bar{v}^k)_{k\in N}$ converges to zero.
        \item\label{it:CP:subseq} Every limit point (if any) of $(\z^k)_{k\in\N}$ belongs to $\zer \Tpd$.
        \item\label{it:CP:bounded} The sequences $(\projQ z^k)_{k\in\N}$, $(\projQ\z^k)_{k\in\N}$ are bounded and their limit points are in $\projQ \zer \Tpd$.
        \item\label{it:CP:full} If in \cref{ass:SWMVICP} the set 
        $
            \projQ 
            \pazocal{S}^\star
        $
        is equal to
        $
            \projQ 
            \zer T
        $, then the sequences
        $\seq{\projQ z^k}$, $\seq{\projQ \bar{z}^k}$ converge
        to some element of $\projQ\zer \Tpd$.
        If additionally $J_{\gamma A}$ and $J_{\tau B^{-1}}$ are (single-valued) continuous,
        then $\seq{\bar{z}^k}$ converges to some 
        $
            z^\star \in \zer \Tpd
        $.
        Finally, if $\lambda_k$ is additionally uniformly bounded in the interval $(0,2)$, then $\seq{z^k}$ also converges to 
        $
            z^\star \in \zer \Tpd
        $.
    \end{enumerate}
    \begin{proof}
\Cref{ass:PPPA:1} is immediate.
Outer semicontinuity of $\Tpd$ follows from that of $A$ and $B$ \cite[Theorem 5.7(a)]{rockafellar2009Variational}, showing that \cref{ass:PPPA:0} holds.
Since $J_{\gamma A}$ and $J_{\tau B^{-1}}$ have full domain
the preconditioned resolvent $(\M + \Tpd)^{-1}\M$ has full domain owing to \cite[Lemma 12.14]{rockafellar2009Variational}, establishing \Cref{ass:PPPA:0.5}.
It only remains to show that \Cref{ass:PPPA:2} holds.
Before that, we first show that \(\etamin\) as defined in \eqref{eq:PPPA:eigcond} is characterized by the expression given in \cref{ass:CP:relaxation:rule}. 

Let \(Z_i \dfn \VL_i \oplus \UL_i \in \R^{(m+n)\times(2m_i)}\) and \(\Un \dfn \VL^\prime \oplus \UL^\prime \in \R^{(m+n)\times (m+n-2r)}\) where $\UL_i$, $\VL_i$, $\UL^\prime$ and $\VL^\prime$ are given by \eqref{eq:L:svd}. Let
\ifjota
\else
    \[
        Z
            {}\dfn{}
        \begin{bmatrix}
            Z_1 & \cdots & Z_d
        \end{bmatrix}
        \in \R^{(m+n) \times (2r)},
    \]
\fi
which by construction has orthonormal columns.
The preconditioner \(\M\) in \eqref{eq:CP:P} can be decomposed as
\ifjota
\else
    \begin{align*}
        \M
            {}={}
        \begin{bmatrix}
            \tfrac{1}{\gamma}\I_n &-L^\top\\
            -L & \tfrac{1}{\tau}\I_m
        \end{bmatrix}
            {}\stackrel{\dueto{\eqref{eq:L:svd}}}{=}{}
        \begin{bmatrix}
            \tfrac{1}{\gamma}\I_n & -\VL \SL \tp \UL\\
            -\UL \SL \tp \VL & \tfrac{1}{\tau} \I_m
        \end{bmatrix}
            {}\stackrel{\dueto{\eqref{eq:projXY}}}{=}{}
        \underbrace*{
        \begin{bmatrix}
            \tfrac{1}{\gamma}\VL \tp \VL & -\VL \SL \tp \UL\\
            -\UL \SL \tp \VL & \tfrac{1}{\tau} \UL \tp \UL
        \end{bmatrix}
        }_{\eqqcolon \,Z \PrS \tp Z}
            +
        \underbrace*{
        \Bigl(\tfrac{1}{\gamma}\VL' \tp {\VL'}\Bigr) \oplus \Bigl(\tfrac{1}{\tau} \UL' \tp {\UL'}\Bigr)
        }_{\eqqcolon \,Z' \Pn \tp {Z'}}
        ,
        \numberthis\label{eq:proof:CP:preconditioner:similarity}
    \end{align*}
\fi
where
\begin{align*}
    \Prs{i}
        {}={}
    \begin{bmatrix}
        \tfrac{1}{\gamma} & -\sigma_i \\
        -\sigma_i & \tfrac{1}{\tau}
    \end{bmatrix},
    \quad 
    \PrS
        {}={}
    \left( (\PrS_1 \otimes \I_{m_1}) \oplus \cdots \oplus (\PrS_d \otimes \I_{m_d})\right)
    \quad\text{and}\quad
    \Pn  
        {}={}
    \tfrac1\gamma \I_{n-r} \oplus \tfrac1\tau \I_{m-r}.
\end{align*}
Using this decomposition, an orthonormal basis for $\range{\M}$ can be constructed.
The matrices \(\Prs{i} \succ 0\) for $i = 2, \hdots, d$  
since 
\(\gamma\tau\sigma_i^2 < 1\), and \(\Prs{1} \succeq 0\) since by 
\Cref{ass:CP:stepsize:rule} $\gamma\tau\|L\|^2 \leq 1$.
Consequently, an orthonormal basis for \(\range{\Prs{i}}\)
is given by $\Urs{i} = \I_2$ for $i = 2, \hdots, d$.
For $i = 1$, an orthonormal basis is given by $\Urs{1} = \I_2$ if \(\gamma\tau\|L\|^2 < 1\) and if \(\gamma\tau\|L\|^2 = 1\) by 
\begin{align*}
    \UrS_1
        = 
    \sqrt{\nicefrac{\tau}{(\gamma + \tau)}}
    \begin{bsmallmatrix}
        1 \\
        -\sqrt{\nicefrac\gamma\tau}
    \end{bsmallmatrix}.
    \numberthis\label{eq:proof:CP:indef:basisP1}
\end{align*}
Owing to the block-diagonal structure of \(\PrS\), an orthonormal basis for \(\range{\PrS}\) is therefore given by
\(
    \UrS
        =
    (\Urs{1} \otimes \I_{m_1}) 
        \oplus 
    \I_{(2r-2m_1)}
\).
For any $A \in \R^{k \times l}$ and any full row rank $B \in \R^{l \times p}$, it holds that $\rge{AB} = \rge{A}$ \cite[Prop. 2.6.3]{bernstein2009matrix}.
Since $Z$ and $Z^\prime$ are full column rank and
\(\Pn\) is full rank,
this implies that
\(
    \range{Z \PrS \tp Z} = \rge{Z\UrS}
\)
and
\(
    \rge{Z' \Pn {Z'}^\top} = \rge{Z'}
\).
Finally, since $Z\UrS$ has orthonormal columns (as $(Z\UrS)^\top Z\UrS = \I$), an orthonormal basis for
$
    \range{\M}
        =
    \range{Z \hat{P} Z^\top} \oplus \range{Z' \hat{P}' Z'^\top}
$
is given by
\begin{align*}
    \basisP 
        = 
    \begin{bmatrix}
        Z\UrS & \Un
    \end{bmatrix}
        {}={}
    \begin{bmatrix}
        Z_1
        (\Urs{1} \otimes \I_{m_1}) &
        Z_2 & \cdots & Z_d & \Un
    \end{bmatrix}.\numberthis\label{eq:proof:CP:basisP}
\end{align*}
Analogously, the \(\DRSRho\)-oblique weak Minty matrix as defined in \eqref{eq:WMSDRS} can be decomposed as
\ifjota
\else
    \begin{align*}
        \DRSRho
            {}\stackrel{\eqref{eq:projXY}}{=}{}
        \smash{
        \underbrace{
        \left(
            \DRSrhocom \VL {\VL}^\top
        \right)
        {}\oplus{}
        \left(
            \DRSrhomon \UL {\UL}^\top
        \right)
        }_{\eqqcolon\,Z\VrS\tp Z}}
            {}+{}
        \smash{
        \underbrace{
        \Bigl(\DRSrhocom^\prime \VL' \tp {\VL'}\Bigr) \oplus \Bigl(\DRSrhomon^\prime \UL' \tp {\UL'}\Bigr)
        }_{\eqqcolon\,Z' \Vn \tp {Z'}}},
        \numberthis\label{eq:proof:CP:V:similarity}
    \end{align*}
\fi
where
\ifjota
\else
    \begin{align*}
        \Vrs{i}
            {}={}&
        \diag(\DRSrhocom, \DRSrhomon), 
        \quad \VrS
            {}={}
        \left((\Vrs{1} \otimes \I_{m_1}) \oplus  \cdots \oplus (\Vrs{d} \otimes \I_{m_d})\right)
        \quad\text{and}\quad
        \Vn  
            {}={}
        \DRSrhocom^\prime \I_{n-r} \oplus \DRSrhomon^\prime \I_{m-r}.
    \end{align*}
\fi
Since \(Z\) and \(Z'\) both have orthonormal columns, i.e., \(\tp Z Z = \I_{2r}\) and \(\tp {Z'} Z' = \I_{m+n-2r}\), and since $\range{Z}$ and $\range{Z'}$ are orthogonal, 
it follows from \eqref{eq:proof:CP:preconditioner:similarity}-\eqref{eq:proof:CP:V:similarity} that
\ifjota
\else
    \begin{align*}
        \tp \basisP \DRSRho \M \basisP
            {}={}
        \begin{bmatrix} 
            \tp \UrS \tp Z\\
            \tp {Z'}
        \end{bmatrix}  
        (Z \VrS \tp Z + Z' \Vn \tp {Z'})(Z \PrS \tp Z + Z' \Pn \tp {Z'})
        \begin{bmatrix}
            Z\UrS & \Un
        \end{bmatrix}
        {}={}&
        \tp \UrS \VrS \PrS \UrS \oplus \tp {Z'} \Vn \Pn Z'.
    \end{align*}
\fi
Therefore, condition \eqref{eq:PPPA:eigcond} becomes
\begin{align*}
    1 + \etamin
        {}={}
    1 + \lambda_{\mathrm min}(\tp \basisP \DRSRho \M \basisP)
        {}={}&
    1 + \min\Bigl(
    \overbrace*{
        \lambda_{\mathrm min}
        (
            \UrS^\top
            \VrS\PrS
            \UrS
        )
    }^{\eqqcolon\,\eta}
    ,
    \overbrace*{
        \lambda_{\mathrm min}
        (   
            {\Un}^\top
            \Vn\Pn 
            \Un
    )
    }^{\eqqcolon\,\eta^\prime}
    \Bigr) > 0.
    \numberthis\label{eq:proof:CP:condition:blocked}
\end{align*}
Due to the block diagonal structure of $\UrS^\top\VrS\PrS\UrS$, it follows that
\ifjota
\else
    \begin{align*}
        \eta
            {}={}&
        \min
            \set{\lambda_{\mathrm{min}}(
                \Urs{i}^\top
                \Vrs{i}\Prs{i}
                \Urs{i})}_{i=1}^d
            {}={}
        \min
        \left\{
            \lambda_{\mathrm{min}}(\tp {\Urs{1}} \Vrs{1}\Prs{1}\Urs{1}),
            \set{\lambda_{\mathrm{min}}(\Vrs{i}\Prs{i})}_{i=2}^d
        \right\}
        \numberthis\label{eq:proof:CP:eigcond:reformulated:blocks}
        \\
            {}={}&
        \min
            \set{\lambda_{\mathrm{min}}\Bigl( 
                \Urs{1}^\top
                \begin{bmatrix}
                    \tfrac{1}{\gamma}\DRSrhocom & - \DRSrhocom \sigma_1\\
                    - \DRSrhomon \sigma_1 & \tfrac{1}{\tau}\DRSrhomon
                \end{bmatrix}
                \Urs{1}\Bigr), 
                \set{
                    \tfrac{1}{2\gamma}\DRSrhocom + \tfrac{1}{2\tau}\DRSrhomon - \theta_{\gamma\tau}(\sigma_{i})
                    }_{i=2}^d},
    \end{align*}
\fi
where $\theta_{\gamma\tau}(\cdot)$ is defined as in \eqref{eq:theta}. 
Note that \(\theta_{\gamma \tau}(\sigma)\) can also be expressed as
\begin{align*}
    \theta_{\gamma \tau}(\sigma) 
        {}\coloneqq{}
    \sqrt{\left(\tfrac{1}{2\gamma}\DRSrhocom - \tfrac{1}{2\tau}\DRSrhomon\right)^2 + \DRSrhocom\DRSrhomon\sigma^2}
        {}={}
    \sqrt{\left(\tfrac{1}{2\gamma}\DRSrhocom + \tfrac{1}{2\tau}\DRSrhomon\right)^2 - \tfrac{\DRSrhocom\DRSrhomon}{\gamma\tau}(1 - \gamma\tau\sigma^2)
    }.
    \numberthis\label{eq:proof:CP:theta}
\end{align*}
If \(\gamma\tau \nrm{L}^2 < 1\), then since $\Urs{1} = \I_{2}$ it follows that $\eta$ is as in \eqref{eq:thm:CP:etamin} since it reduces to
\begin{align*}
    \eta
        {}={}&
    \min
    \set{
        \tfrac{1}{2\gamma}\DRSrhocom + \tfrac{1}{2\tau}\DRSrhomon - \theta_{\gamma\tau}(\sigma_{i})
    }_{i=1}^d
        {}={}
    \begin{cases}
        \tfrac{1}{2\gamma}\DRSrhocom + \tfrac{1}{2\tau}\DRSrhomon - \theta_{\gamma \tau}(\sigma_{d}),
        &\;\;\, \text{if } \DRSrhocom\DRSrhomon < 0,\\
        \tfrac{1}{2\gamma}\DRSrhocom + \tfrac{1}{2\tau}\DRSrhomon - \theta_{\gamma \tau}(\|L\|),
        &\;\;\, \text{if } \DRSrhocom\DRSrhomon \geq 0,
    \end{cases}
\end{align*}
where equivalence \eqref{eq:proof:CP:theta} was used in the second equality to show for which singular value the minimum is attained.
Otherwise, \(\gamma\tau \nrm{L}^2 = 1\) and we have that    
\(
    \UrS_1
        = 
    \sqrt{\nicefrac{\tau}{(\gamma + \tau)}}
    \begin{bsmallmatrix}
        1 \\
        -\sqrt{\nicefrac\gamma\tau}
    \end{bsmallmatrix},
\)
which also matches \eqref{eq:thm:CP:etamin} since 
\begin{align*}
    \eta
        {}={}&
    \min
    \set{\lambda_{\mathrm{min}}\Bigl( 
        \tfrac{\tau}{\gamma + \tau}\begin{bsmallmatrix}1 & -\sqrt{\nicefrac\gamma\tau}\end{bsmallmatrix}
        \begin{bsmallmatrix}
            \tfrac{1}{\gamma}\DRSrhocom & - \DRSrhocom \sigma_1\\
            - \DRSrhomon \sigma_1 & \tfrac{1}{\tau}\DRSrhomon
        \end{bsmallmatrix}
        \begin{bsmallmatrix}1 \\ -\sqrt{\nicefrac\gamma\tau}\end{bsmallmatrix}
        \Bigr), 
        \set{
            \tfrac{1}{2\gamma}\DRSrhocom + \tfrac{1}{2\tau}\DRSrhomon - \theta_{\gamma\tau}(\sigma_{i})
            }_{i=2}^d}\\
        {}={}&
    \min
    \left\{
    \tfrac{1}{\gamma}\DRSrhocom + \tfrac{1}{\tau}\DRSrhomon,
    \set{
        \tfrac{1}{2\gamma}\DRSrhocom + \tfrac{1}{2\tau}\DRSrhomon - \theta_{\gamma\tau}(\sigma_{i})
        }_{i=2}^d
    \right\}\\
        {}={}&
    \begin{cases}
        \tfrac{1}{\gamma}\DRSrhocom + \tfrac{1}{\tau}\DRSrhomon, & \;\;\, \text{if } d=1,\\
        \left\{
        \begin{array}{ll}
            \tfrac{1}{2\gamma}\DRSrhocom + \tfrac{1}{2\tau}\DRSrhomon - \theta_{\gamma\tau}(\sigma_d)
            ,
            &\;\;\,\hspace{\maxmin} \text{if } \DRSrhocom\DRSrhomon < 0\\
            \tfrac{1}{2\gamma}\DRSrhocom + \tfrac{1}{2\tau}\DRSrhomon - \theta_{\gamma\tau}(\sigma_2)
            ,
            &\;\;\,\hspace{\maxmin} \text{if } \min\{\DRSrhocom, \DRSrhomon\} \geq 0\\
            \tfrac{1}{\gamma}\DRSrhocom + \tfrac{1}{\tau}\DRSrhomon,
            &\;\;\,\hspace{\maxmin} \text{otherwise}
        \end{array}
        \right\},
        & \;\;\, \text{if } d>1.
    \end{cases}
\end{align*}
where we used that $\sqrt{\gamma\tau}\sigma_1 = 1$ in the first equality and equivalence \eqref{eq:proof:CP:theta} to show the second equality. 

By definition of \(\Un\), \(\Vn\) and \(\Pn\), it holds that
\begin{equation}\label{eq:etaprime}
        \eta^\prime
        {}={}
    \lambda_{\mathrm min}
        \left(   
            \bigl(
            \tfrac1\gamma\DRSrhocom^\prime \tp {\VL'} \VL'\bigr)
            \oplus 
            \bigl(
            \tfrac1\tau\DRSrhomon^\prime \tp {\UL'} \UL'\bigr)
    \right)
        {}={}
    \lambda_{\mathrm min}
        \Bigl(
            \bigl(\tfrac1\gamma \DRSrhocom^\prime\I_{n-r}\bigr)
                \oplus
            \bigl(\tfrac1\tau \DRSrhomon^\prime \I_{m-r}\bigr) 
        \Bigr),
\end{equation}
which matches the definition of $\eta^\prime$ provided in \Cref{tab:CP:degen}. 
Having shown that the values of \(\eta\) and \(\eta^\prime\) match the ones provided in \cref{ass:CP:relaxation:rule}, it remains to show that \cref{ass:CP:stepsize:rule} ensures that
$1 + \etamin \dfn 1 + \min\{\eta, \eta^\prime\} > 0$ as required by \cref{ass:PPPA:2}. We proceed with two intermediate claims.  

\begin{claims}
    \item \label{claim:proof:CP:1}
    \textit{
    The pair $(\gamma, \tau) \in \R^2_{++}$ satisfies $ \gamma\tau \in (0, \nicefrac{1}{\nrm{L}^2}]$ and $1 + \eta > 0$ if and only if $\gamma$ and $\tau$ comply with \Cref{ass:CP:stepsize:rule}.
    }\newline
    Consider the six different cases for $\eta$ from \eqref{eq:thm:CP:etamin}.
    For the first and fourth case this follows from \cref{it:cor:quadratic:CP:existence:mixed}, for the second case from \cref{it:cor:quadratic:CP:existence:pos,it:cor:quadratic:CP:existence:neg} since $[\DRSrhocom]_-[\DRSrhomon]_- < \tfrac{1}{4\nrm{L}^2}$ by \Cref{eq:existence:condition}
    and for the fifth case from \cref{it:cor:quadratic:CP:existence:pos}. Finally, for the third and the sixth case use \Cref{lem:quadratic}, by plugging in $\tau = \tfrac{1}{\gamma\nrm{L}^2}$ into $1 + \tfrac{1}{\gamma}\DRSrhocom + \tfrac{1}{\tau}\DRSrhomon > 0$ and observing that $\delta = 1$.
    \item \label{claim:proof:CP:2}
    \textit{
        If $\gamma$ and $\tau$ comply with \Cref{ass:CP:stepsize:rule} then $1 + \eta^\prime > 0$.
    }\newline
    Recall that as noted in the discussion after \cref{ass:CP}, the value of \(\DRSrhocom^\prime\) (resp. \(\DRSrhomon^\prime\)) is irrelevant when $L$ is full column rank (resp. full row rank), and was considered equal to $\infty$ in this case.
    In view of \eqref{eq:etaprime}, $\gamma$ and $\tau$
    satisfy $1 + \eta^\prime > 0$ if and only if
    \begin{align*}
        \gamma > -\DRSrhocom^\prime
        \quad\text{and}\quad
        \tau > -\DRSrhomon^\prime.
        \numberthis\label{eq:proof:CP:etaprime:stepsizes}
    \end{align*}
    We distinguish between two cases.
    \begin{proofitemize}
        \item 
        $\min\{\DRSrhocom,\DRSrhomon\} \geq 0$: Then, owing to \Cref{eq:existence:condition} it holds that $\DRSrhocom^\prime \geq 0$ and $\DRSrhomon^\prime \geq 0$, so that \eqref{eq:proof:CP:etaprime:stepsizes} is implied trivially by \Cref{ass:CP:stepsize:rule}.
        \item 
        $\min\{\DRSrhocom,\DRSrhomon\} < 0$:
        We begin by deriving two auxiliary inequalities that are essential for showing that \eqref{eq:proof:CP:etaprime:stepsizes} holds.
        By definition of $\delta$ from \eqref{eq:def:delta} and using 
        $
            \DRSrhocom\DRSrhomon
                =
            [\DRSrhocom\DRSrhomon]_+
                -
            [\DRSrhocom\DRSrhomon]_-
        $
        it holds that for any parameter \(\kappa \geq 1\) that
        \begin{equation}\label{eq:proof:stepsize:delta:1}
            \delta - \kappa \DRSrhocom\DRSrhomon\nrm{L}^2 
                {}={}
            1
                +
            [\DRSrhocom\DRSrhomon]_-\big(\sigma_{d}^2 + (\kappa -1)\|L\|^2\big) - \kappa [\DRSrhocom\DRSrhomon]_+ \nrm{L}^2
                {}\geq{}
            1 - \tfrac{\kappa}4,
        \end{equation}
        where the inequality follows from the fact that either both \(\DRSrhocom\) and \(\DRSrhomon\) are negative, in which case inequality
        \eqref{eq:existence:condition} is used, or only one is negative, in which case the quantity is greater than or equal to 1. 
        Moreover, 
        \begin{align*}
            0 \leq \delta^2 - 4\DRSrhocom\DRSrhomon\nrm{L}^2
                {}={}&
            \bigl(\delta - 2\DRSrhocom\DRSrhomon\nrm{L}^2\bigr)^2 
            +
            4\DRSrhocom\DRSrhomon\nrm{L}^2
            \big(\delta - 1 - \DRSrhocom\DRSrhomon\nrm{L}^2\big)
            \\
                {}={}&
            \bigl(\delta - 2\DRSrhocom\DRSrhomon\nrm{L}^2\bigr)^2 
            -
            4\DRSrhocom\DRSrhomon\nrm{L}^2
            \big( [\DRSrhocom\DRSrhomon]_+\nrm{L}^2
                -
            [\DRSrhocom\DRSrhomon]_-\sigma_d^2\big)
            \\
                {}\leq{}&
            \bigl(\delta - 2\DRSrhocom\DRSrhomon\nrm{L}^2\bigr)^2,
            \numberthis\label{eq:proof:CP:auxiliary}
        \end{align*}
        where the last inequality follows from the trivial observation that for any \(a\in \R\), 
        \(a[a]_+\geq 0\) and \(a[a]_-\leq 0\).
        \begin{proofitemize}
            \item
            $\gamma > -\DRSrhocom^\prime$:
            Using \eqref{eq:proof:CP:auxiliary}, we obtain that
            \(
                    \delta + \sqrt{\delta^2 - 4\DRSrhocom\DRSrhomon\nrm{L}^2}
                        {}\leq{}
                    2(\delta - \DRSrhocom\DRSrhomon\nrm{L}^2)
            \).
            Consequently, by the definition of $\gamma_{\rm min}$ from \eqref{eq:stepsizes:CP:bounds:gam}, it follows that
            \begin{align*}
                &
                \DRSrhocom < 0 \Longrightarrow
                -\DRSrhocom^\prime
                    {}\stackrel{\eqref{eq:existence:condition}}{\leq}{}
                \tfrac{-\DRSrhocom}{\delta - \DRSrhocom\DRSrhomon\nrm{L}^2}
                    {}\leq{}
                \gamma_{\rm min}.
                \numberthis\label{eq:proof:CP:etaprime:implications-gamma}
            \end{align*}
            Since \Cref{ass:CP:stepsize:rule} ensures that
            $\gamma > 0 = \DRSrhocom^\prime$ if $\DRSrhocom \geq 0$ and $\gamma > \gamma_{\rm min}$ if $\DRSrhocom < 0$,
            it thus follows from \eqref{eq:proof:CP:etaprime:implications-gamma}
            that indeed
            $\gamma > -\DRSrhocom^\prime$.

            \item
            $\tau > -\DRSrhomon^\prime$:
            Observe that regardless of the sign of $\DRSrhocom$ it holds that $\gamma(\delta-\DRSrhocom\DRSrhomon\nrm{L}^2) + \DRSrhocom > 0$;
            if $\DRSrhocom \geq 0$, then this holds owing to \eqref{eq:proof:stepsize:delta:1} since $\gamma > 0$, if $\DRSrhocom < 0$, then this inequality remains true since in that case $\gamma > \gamma_{\rm min}$ by \Cref{ass:CP:stepsize:rule} and using \eqref{eq:proof:CP:etaprime:implications-gamma}.
            We proceed by first showing that
            \begin{align*}
                \tau_{\rm min}(\gamma)
                    {}={}
                \tfrac{[-\DRSrhomon(\gamma + \DRSrhocom)]_+}{\gamma(\delta-\DRSrhocom\DRSrhomon\nrm{L}^2) + \DRSrhocom}
                    {}\geq{}
                \tfrac{[-\DRSrhomon]_+}{\delta - \DRSrhocom\DRSrhomon\nrm{L}^2}.
                \numberthis\label{eq:proof:CP:etaprime:implications-tau}
            \end{align*}
            If $\DRSrhomon \geq 0$, then \eqref{eq:proof:CP:etaprime:implications-tau} holds trivially.
            Let $\DRSrhomon < 0$.
            If $\DRSrhocom \geq 0$, then $\gamma + \DRSrhocom > 0$ since $\gamma > 0$.
            If $\DRSrhocom < 0$, then it also holds that $\gamma + \DRSrhocom > 0$ since in that case $\gamma > \gamma_{\rm min} \geq \tfrac{-\DRSrhocom}{1-\DRSrhomon\DRSrhocom\nrm{L}^2} > -\DRSrhocom$ by \Cref{ass:CP:stepsize:rule} and using \eqref{eq:proof:CP:etaprime:implications-gamma}.
            Therefore, it holds that 
            \(
                -\DRSrhomon(\gamma + \DRSrhocom) > 0
            \).
            Having previously established that the denominators in \eqref{eq:proof:CP:etaprime:implications-tau} are strictly positive this inequality therefore is equivalent to
            \begin{align*}
                -\DRSrhomon(\gamma + \DRSrhocom)(\delta - \DRSrhocom\DRSrhomon\nrm{L}^2) \geq -\DRSrhomon(\gamma(\delta-\DRSrhocom\DRSrhomon\nrm{L}^2) + \DRSrhocom).
            \end{align*}
            Substituting $\delta$ from \eqref{eq:def:delta} and using
            $
                \DRSrhocom\DRSrhomon =
                [\DRSrhocom\DRSrhomon]_+
                    -
                [\DRSrhocom\DRSrhomon]_-
            $, rearrange it to obtain
            \(
                [\DRSrhocom\DRSrhomon]_-^2\sigma_{d}^2 + [\DRSrhocom\DRSrhomon]_+^2\nrm{L}^2
                    {}\geq{}
                0
            \), which holds trivially.
            Finally, having established that \Cref{eq:proof:CP:etaprime:implications-tau} holds
            and noting that 
            $\tau > \tau_{\rm min}(\gamma)$
            owing to \Cref{ass:CP:stepsize:rule}, it follows immediately that
            $\tau > -\DRSrhomon^\prime$
            by definition of $\DRSrhomon^\prime$ from \eqref{eq:existence:condition}.
        \end{proofitemize}
    \end{proofitemize}
\end{claims}
Combining \cref{claim:proof:CP:1,claim:proof:CP:2}, it follows that \Cref{ass:CP:stepsize:rule} is equivalent to the preconditioner $\M$ being positive semidefinite
and $1 + \etamin = 1 + \min\{\eta, \eta^\prime\} > 0$, establishing that \Cref{ass:PPPA:2} holds.
Finally, \((\M + \Tpd)^{-1} \circ \M\) is (single-valued) continuous if $J_{\gamma A}$ and $J_{\tau B^{-1}}$ are (single-valued) continuous, since $\bar z^k \in (\M + \Tpd)^{-1}\M z^k$ if and only if
\(
    \bar{x}^k
        {}\in{}
    J_{\gamma \A}\big(x^k - \gamma L^\top y^k\big)
\)
and
\(
    \bar y^k 
        {}\in{}
    J_{\tau \inv\B}\big(y^k + \tau L(2\bar x^k - x^k)\big)
\).
All claims for \ref{eq:CP} follow directly from 
\Cref{thm:pppa}, using the equivalence from \Cref{lem:equivalence:PPPA-PDHG}.
\qedhere
    \end{proof}
\end{theorem}
Note that when $\lambda_k$ is uniformly bounded in the interval $\bigl(0, 2(1 + \etamin)\bigr)$, a rate of $O(\nicefrac1N)$ can be obtained for $\min_{k=0,1,\ldots,N} \|\bar v^k\|^2$ by telescoping (see \cite[thm. 2.4(iv)]{evens2023convergence}).

Observe that \Cref{thm:CP:full} discusses not only the convergence of $\seq{z^k}$, but also of its projection onto the range of the preconditioner $\seq{\proj_{\range{\M}} z^k}$.
Notably, convergence of $\seq{\proj_{\range{\M}} z^k}$ is established under weaker assumptions than for $\seq{z^k}$. 
When $\M$ is positive definite, meaning 
$\gamma \tau < \nicefrac1{\nrm{L}^2}$, this is irrelevant because in this case the range of \(\M\) is full.
However, in the positive semidefinite case, when $\gamma \tau = \nicefrac1{\nrm{L}^2}$, these sequences are no longer identitical.
This observation is not surprising, as it is a natural extension of the convergence results for DRS, i.e., when $L = \I$ and $\gamma = \nicefrac1\tau$.
In particular, in the DRS setting it was shown that the convergence of $\seq{\proj_{\range{\M}} z^k}$ to $\proj_{\range{\M}} \zer \Tpd$ is equivalent to the convergence of the sequence $\seq{s^k} \dfn \seq{x^k - \gamma y^k}$ to a point
$
    s^\star
        {}\in{}
    \left\{
    x^\star - \gamma y^\star\,\middle|\, (x^\star, y^\star) \in \zer \Tpd
    \right\}
$,
and a solution of the primal inclusion is obtained through $J_{\gamma A}(s^\star)$ \cite{lions1979Splitting,svaiter2011Weak,evens2023convergence}.
In the following proposition, this interpretation is generalized to arbitrary $L$ matrices. 

\begin{proposition}[convergent sequences]\label{rem:PDHG:sequences} 
    Using the SVD of $L$ from \eqref{eq:L:svd}, define 
    \ifjota
\else
        \begin{align*}
            \basisP 
                {}={}
            \begin{bmatrix}
                \sqrt{\tfrac{\tau}{\gamma + \tau}}
                Z_1
                \begin{bsmallmatrix}\I_{m_1} \\ -\sqrt{\nicefrac\gamma\tau}\I_{m_1}\end{bsmallmatrix} &
                Z_2 & \cdots & Z_d & \Un
            \end{bmatrix}
                {}={}
            \begin{bmatrix}
                \sqrt{\tfrac{\tau}{\gamma + \tau}}\VL_1 & \VL_2 & 0 & \cdots & \VL_d & 0 & \VL' & 0\\
                -\sqrt{\tfrac{\gamma}{\gamma + \tau}}\UL_1 & 0 & \UL_2 & \cdots & 0 & \UL_d & 0 & \UL'
            \end{bmatrix}.
            \numberthis\label{eq:interpretation:basisP}
        \end{align*}
    \fi
    Consider a sequence $\seq{z^k} = \seq{x^k,y^k}$ generated by \ref{eq:CP} starting from $z^0 \in \R^{n+m}$, where $\gamma\tau=\tfrac{1}{\nrm{L}^2}$ and define
    \begin{align*}
        s^k 
            {}\dfn{}
        \basisP^\top
        \begin{bmatrix}
            x^k\\
            y^k
        \end{bmatrix}
        \quad\text{and}\quad
        \pazocal{T}
            {}\dfn{}
        \left\{
            \basisP^\top
            \begin{bmatrix}
                x^\star\\
                y^\star
            \end{bmatrix}
            \,\middle|\, (x^\star, y^\star) \in \zer \Tpd\right\}.
            \numberthis\label{eq:CP:interpretation:seq}
    \end{align*}
    Then, 
    \ifjota
\else
        the following statements hold.
        \begin{enumerate}
            \item The limit points of $(\projQ z^k)_{k\in\N}$ are in $\projQ \zer \Tpd$ if and only if the limit points of $\seq{s^k}$ are in $\pazocal{T}$.
            \item The sequence $(\projQ z^k)_{k\in\N}$ converges to 
            $\projQ \zer \Tpd$ if and only if $\seq{s^k}$ converges to 
            $\pazocal{T}$.
        \end{enumerate} 
    \fi
    \begin{proof}
        Notice that $\basisP$ is an orthonormal basis for
        $
            \range{\M}
        $,
        obtained by plugging in \eqref{eq:proof:CP:indef:basisP1} into \eqref{eq:proof:CP:basisP}.
        As a result, the claims follow directly from \eqref{eq:CP:interpretation:seq}, using the fact that \(\proj_{\range{\M}} = \basisP \basisP^\top\) and that \(\basisP\) has orthonormal columns.
    \end{proof}
\end{proposition}
Notably, \Cref{rem:PDHG:sequences} along with \Cref{thm:CP:full} establishes the convergence of an $(m+n-m_1)$\hyp{}dimensional sequence $\seq{s^k}$ of \ref{eq:CP} when $\gamma\tau = \nicefrac{1}{\nrm{L}^2}$.
Since
$
    s^k
        =
    \nicefrac{1}{\sqrt{1+\gamma^2}}(x^k - \gamma y^k)
$
when $L = \I$ and $\tau=\nicefrac1\gamma$, it follows immediately that
\Cref{thm:CP:full} matches
the convergence results for DRS obtained in \cite[Thm. 3.3]{evens2023convergence}.
A simple example where $\seq{s^k}$ converges while $\seq{z^k}$ diverges is provided in 
\Cref{ex:CP:second}.


In the remainder of this section, two theoretical examples are presented to demonstrate some of the main attributes of the convergence results in \Cref{thm:CP:full}. Supplementary Python code verifying these results can be found on GitHub\footnote{\url{https://github.com/brechtevens/Minty-CP-examples}.}
and the proofs are deferred to \Cref{proof:ex:CP:example5}.

In the first example, the tightness of the bounds on $\gamma$, $\tau$ and $\lambda_k$ from \Cref{thm:CP:full} is demonstrated through a simple system of linear equations. In this setting, the iterations of \ref{eq:CP} can be expressed as a linear dynamical system,
so that tight bounds can be obtained by ensuring stability.
In this example, an artificial parameter $c$ is introduced when splitting the problem into the form $A + L^\top B L$.
While this parameter may appear inconsequential at first sight, it controls the value of $\DRSrhomon^\prime$ within the problem
(see \Cref{it:CP:example5:2}).

\begin{example}[saddle point problem]\label{ex:CP:example5}
    Consider the problem of finding a zero of 
    \ifjota\else
        the following structured linear inclusion
    \fi
    \begin{equation}\label{eq:ex5-Tp}
        0 \in \Tp x
            {}={}
        \begin{bmatrix}b\ell^2 & a \\ -a & b\ell^2\end{bmatrix}x
            {}={}
        \overbrace{\begin{bmatrix}0 & a \\ -a & 0\end{bmatrix}}^A x + 
        \overbrace*{\begin{bmatrix}\ell & 0 & 0\\ 0 & \ell & 0\end{bmatrix}}^{\tp L}
        \overbrace*{\begin{bmatrix}b & 0 & 0\\ 0 & b & 0\\ 0 & 0 & c\end{bmatrix}}^B
        \overbrace*{\begin{bmatrix}\ell & 0 \\ 0 & \ell\\ 0 & 0\end{bmatrix}}^L x,
    \end{equation}
        where $a, b, \ell \in \R \setminus \{0\}$ and
        $
            c \geq 0
        $.
        Suppose that
        $
            a^2 \neq b^2 \ell^4
        $
        when $b < 0$.
        Any solution to the inclusion problem $0 \in \Tp x$ is a minimax solution of
        \(
            f(x_1,x_2) \dfn ax_1x_2 + \tfrac{b\ell^2}{2}(x_1^2 - x_2^2) 
        \)
        when \(b > 0\) and a maximin solution 
        when \(b < 0\). 
        Consider the sequence $\seq{z^k} = \seq{x^k,y^k}$ generated by applying \ref{eq:CP} to \eqref{eq:ex5-Tp} with $\tau = \tfrac{1}{\gamma \nrm{L}^2}$ and fixed relaxation parameter $\lambda$.
        Then, the following assertions hold.
    \begin{enumerate}
        \item \label{it:CP:example5:1} 
        The sequence
        $\seq{\projQ z^k}$ converges if and only if $\lambda \in \left(0, \bar{\lambda}\right)$ and 
        $\seq{z^k}$ converges if and only if 
        $\lambda \in \left(0, \min\bigl\{2, \bar{\lambda}\bigr\}\right)$, where
        \begin{align*}
            \begin{aligned}
                \bar{\lambda} 
                    {}\dfn{} 
                \min\left\{
                    2\Bigl(1 + \tfrac{b\ell^2}{\gamma (a^2 + b^2\ell^4)} + \tfrac{b a^2 \ell^2 \gamma}{a^2 + b^2\ell^4}\Bigr),
                    2\left(1+\gamma c \ell^2\right)
                \right\}.
            \end{aligned}\numberthis\label{eq:example5:lambda}
        \end{align*}
        This upper bound is strictly positive 
        if and only if 
        \begin{align*}
            \gamma \in
            \begin{cases}
                (0, +\infty), &\quad \text{if } b > 0,\\
                \left(      
                    \tfrac{-2b\ell^2}{a^2 + b^2\ell^4 + |a^2 - b^2\ell^4|},
                    \tfrac{a^2 + b^2\ell^4 + |a^2 - b^2\ell^4|}{-2ba^2},                         
                \right), &\quad \text{if } b < 0.
            \end{cases}
            \numberthis\label{eq:example5:gamma}
        \end{align*}
        \item \label{it:CP:example5:2}
        The operator $\Tpd$ has a 
        \(
            \DRSRho
                =
            \tfrac{b\ell^2}{a^2 + b^2\ell^4} \proj_{\range{\tp L}}
            \oplus\
            \Bigl(
                \tfrac{b a^2}{a^2 + b^2\ell^4} \proj_{\range{L}} + c \proj_{\ker{\tp L}} 
            \Bigr) 
        \)\hyp{}oblique weak Minty solution at zero on
        \(\pazocal{U} = \R^5 \times \range{\M}\).
        \item \label{it:CP:example5:3}
        Recalling that \Cref{ass:SWMVICP} only needs to hold on $\pazocal{U}$ (see the discussion below \Cref{thm:pppa}),
        all assumptions of \Cref{it:CP:full} hold. Moreover, this result
        is tight in the sense that 
        \Cref{ass:CP:stepsize:rule} corresponds to
        \eqref{eq:example5:gamma}
        and
        the bounds from 
        \Cref{ass:CP:relaxation:rule}
        match those
        from \eqref{eq:example5:lambda}.
        \item \label{it:CP:example5:4}
        The range of parameters \(a\), \(b\), \(c\) and \(\ell\) for which \ref{eq:CP} converges 
        includes cases where neither the primal, nor the dual, nor the primal-dual problem are monotone. For instance, when $a = 10$, 
        $
            b = -\tfrac14
        $,
        $c = 0$
        and $\ell = 2$, by \cref{it:CP:example5:1}, $\seq{z^k}$ then converges if and only if
        \(
            \gamma \in 
            \bigl(
                \tfrac1{100}, 1
            \bigr) 
        \)
        and
        \(
            \lambda 
                {}\in{}
            \bigl(
                0, 2 - \tfrac{2}{101\gamma} - \tfrac{200\gamma}{101}
            \bigr).
        \)\qedhere
    \end{enumerate}
\end{example}
\begin{figure}
    \centering
    \ifjota
\else
        \includetikz{Examples/second/shadow}
    \fi
    \caption{
    \ifjota
\else
        Convergence of the sequence
        $
            \seq{s^k}
                =
            \bigl(
                \VL_1^\top x^k - \UL_1^\top y^k
                , \;
                \VL_{2:}^\top x^k
                , \;
                \UL_{2:}^\top y^k
            \bigr)_{k \in \N}
        $
        from \Cref{ex:CP:second} for $n = 3$, $\ell_2 = \nicefrac12$, $\ell_3 = \nicefrac{1}{5}$ and $\lambda = 2.1$.
        (a)
        Norm of the sequence $\seq{z^k} = \seq{x^k, y^k}$. This sequence does not converge, since $\lambda$ has been selected larger than two (see \Cref{it:CP:full}).
        (b)
        Norm of the sequence $\seq{\VL_1^\top x^k - \UL_1^\top y^k}$, which converges to zero.
        (c)
        Visualization of the primal sequences $\seq{x^k}$ and $\seq{\VL_{2:}^\top x^k}$. It can be seen that although $\seq{x^k}$ does not converge (its first coordinate diverges), its projection onto the 2-dimensional space spanned by the columns of $\VL_{2:}$ does converge to zero (marked by a green dot).
        (d)
        Visualization of the dual sequences $\seq{y^k}$ and $\seq{\UL_{2:}^\top y^k}$. Analogous to the primal setting, $\seq{y^k}$ diverges while $\seq{\UL_{2:}^\top y^k}$ converges to zero.
    \fi
    }
    \label{fig:ex:second:sequences}
\end{figure}
The second example focusses on a particular instance of \Cref{thm:CP:full} where $\DRSrhocom$ and $\DRSrhomon$ are both strictly positive, the number of distinct singular values of $L$ is strictly larger than $1$ and 
$\gamma\tau\nrm{L}^2 = 1$.
Then, the admissible range for $\lambda$ depends on the second largest singular value of $L$ (see \Cref{ass:CP:relaxation:rule}).
The following example demonstrates that this dependence is not merely an artifact of our analysis, but is also observed in practice.
Moreover, \Cref{fig:ex:second:sequences} illustrates that using relaxation parameters greater than two may cause the non-projected sequence $\seq{z_k}$ to diverge, in line with \Cref{it:CP:full}.

\begin{example}[influence of singular values]\label{ex:CP:second}
    Let $n \in \{2, 3\hdots\}$ and $L = \diag\left(1, \ell_2 \hdots \ell_n\right)$, where $|\ell_k| < 1$, $\forall k \in \{2, \hdots, n\}$. Let
    \ifjota
\else
        $$
            A = \diag\Bigl(
                1, 
                1 + \sqrt{
                    1 - \ell_{\smash{2}}^2
                }, 
                \hdots,
                1 + \sqrt{
                    1 - \ell_{\smash{n}}^2
                }
            \Bigr)
            \quad\text{and}\quad
            B = \diag\Bigl(
                1,
                \tfrac{1}{1 + \sqrt{1 - \ell_2^2}},
                \hdots,
                \tfrac{1}{1 + \sqrt{1 - \ell_{\smash{n}}^2}}
                \Bigr).
        $$ 
    \fi
    Consider the sequences $\seq{z^k} = \seq{x^k, y^k}$ and $\seq{s^k}$ generated by applying \ref{eq:CP} to $0 \in Ax + L^\top B Lx$ with $\gamma = \tau = 1$ and constant $\lambda$, where $s^k$ is defined as in \eqref{eq:CP:interpretation:seq}.
    Then, the following assertions hold.
    \newcommand{\setlenwidths}[2]{%
    \newlength{\lenA}%
    \newlength{\lenB}%
    \setlength{\lenA}{#1}%
    \setlength{\lenB}{#2}%
    }
    \ifjota
\else
        \setlenwidths{0.53\textwidth}{0.41\textwidth}
    \fi
    \begin{enumerate}[leftmargin=*]
        \item \label{ex:CP:second:Tpd}
        \begin{minipage}[t]{\lenA}
        The associated primal-dual operator $\Tpd$ has a $\bigl(\tfrac12 \I_n \oplus \tfrac12 \I_n\bigr)$\hyp{}oblique weak Minty solution at $(0, 0) = \zer \Tpd$.
        \item \label{ex:CP:second:theorem}
        By \Cref{it:CP:full} and \Cref{rem:PDHG:sequences},
        \ifjota
\else
            both the sequences
        \fi
        $(\projQ z^k)_{k\in\N}$ and $\seq{s^k}$ converge to zero
        if $\lambda$ is selected according to \Cref{ass:CP:relaxation:rule}, i.e., if $\lambda \in (0, \bar \lambda)$, where
        \begin{align*}
            \bar \lambda 
                {}={}&
            2\left(
                1 + \tfrac{1}{2\gamma}\DRSrhocom + \tfrac{1}{2\tau}\DRSrhomon - \theta_{\gamma \tau}(\max\{|\ell_{2}|, \hdots, |\ell_{n}|\})
            \right)
            \\
                {}={}& 
            3 - \max\{|\ell_{2}|, \hdots, |\ell_{n}|\}.
        \end{align*}
        \item \label{ex:CP:second:spectral}
            Let $n = 3$, $\ell_2 \in (0, 1)$ and $\ell_3 = \nicefrac{1}{5}$. 
            Then, the set of relaxation parameters for which the sequences $\seq{\projQ z^k}$ and $\seq{s^k}$ converge is almost entirely covered by \Cref{it:CP:full} (see \Cref{fig:ex:second:bound}).
        \end{minipage}%
        \begin{minipage}[t]{\lenB}
            \begin{figure}[H]
                \centering
                \ifjota
\else
                    \vspace{-0.69cm}
                    \includetikz{Examples/second/spectral}%
                \fi
                \caption{The upper bounds $\bar \lambda$ and $\bar \lambda_{\rm spectral}$ for \Cref{ex:CP:second}, where $\bar \lambda_{\rm spectral}$ is obtained by examining the spectral radius of the algorithmic operator
                \(
                    \proj_{\range{\M}}\left(\I + \lambda\bigl((\M + \Tpd)^{-1}\M - I\bigr)\right)
                \) 
                (see \eqref{eq::CP:second:spectral:lambda}).
                }%
                \label{fig:ex:second:bound}
            \end{figure}
        \end{minipage}
    \end{enumerate}
\end{example}


\section{Semimonotone operators}\label{sec:semi}

In this section, we provide calculus rules for the class of \((\Mon,\Com)\)-semimonotone operators defined in \Cref{def:semimonotonicity}, generalizing the class of \((\mon,\com)\)-semimonotone operators introduced in \cite[Sec. 4]{evens2023convergence}.
Sufficient conditions for the convergence of \ref{eq:CP} applied to \eqref{prob:composite} for \((\Mon,\Com)\)-semimonotone operators $A$ and $B$ will be provided in \Cref{sec:CP:semi}.
The proofs of the calculus rules in this section are deferred to \Cref{proof:prop:semi:young}.

For some choices of $\Mon$ and $\Com$, it follows from the Fenchel-Young inequality that all operators satisfy the definition of $(\Mon, \Com)$\hyp{}semimonotonicity, as stated below.

\begin{proposition}\label{prop:semi:young}
    Let $\Mon,\Com\in\sym{n}$. 
    If $\Mon \prec 0$, $\Com \prec 0$ and $\Mon\preceq\tfrac{1}{4}\Com^{-1}$, then 
    all operators $A : \R^n \rightrightarrows \R^n$ 
    satisfy the definition of \((\Mon, \Com)\)\hyp{}semimonotonicity.
\end{proposition}

The following proposition demonstrates that the maximality of a semimonotone operator is preserved when the semimonotonicity moduli $\Mon$ and $\Com$ are tightened.
\begin{proposition}[preservation of maximality]\label{prop:semi:maximality}
    Let operator \(T : \R^n \rightrightarrows \R^n\) be \((\Mon, \Com)\)\hyp{}semimonotone. If \(T\) is maximally \((\Mon', \Com')\)\hyp{}semimonotone for some $\Mon' \preceq \Mon$ and $\Com' \preceq \Com$, then \(T\) is maximally \((\Mon, \Com)\)\hyp{}semimonotone.
\end{proposition}
\def\D{D}
\def\F{T}
In what follows, various basic properties of $(\Mon, \Com)$-semimonotone operators will be provided.
For instance, by definition, their inverses belong to the same class of operators, with the roles of \(\Mon\) and \(\Com\) reversed.
Additionally, the following proposition analyzes scaling and shifting of semimonotone operators, as well as the cartesian product of two semimonotone operators.
\begin{proposition}[inverting, shifting, scaling and cartesian product]\label{lem:calculus}
    Let operator $A : \R^n \rightrightarrows \R^n$ be (maximally) $(\MonA,\ComA)$\hyp{}semimonotone \optional{at $(\other{x}_A, \other{y}_A) \in \gph A$} and operator $B : \R^m \rightrightarrows \R^m$ be (maximally) $(\MonB,\ComB)$\hyp{}semimonotone \optional{at $(\other{x}_B, \other{y}_B) \in \gph \B$}. Let $\alpha \in \R_{++}$.
    Then, the following hold.
    \begin{enumerate}
        \item\label{prop:calculus:inverse} The inverse operator $A^{-1}$ is (maximally) $(\Com_A, \Mon_A)$\hyp{}semimonotone [at $(\other{y}_A, \other{x}_A)\in \graph A^{-1}$].
        \item\label{lem:calculus:scaling}
        For all $u,w \in \R^n$, operator $T(x) \dfn w + \alpha \A (x+u)$ is (maximally) $(\alpha\MonA,\alpha^{-1}\ComA)$\hyp{}semimonotone \optional{at $(\other{x}_A-u,w + \alpha\other{y}_A)$}.
        \item\label{lem:calculus:cartesian} 
        Operator 
        $T \dfn A \times B$
        is (maximally) $\bigl(\MonA \oplus \MonB,\ComA \oplus \ComB\bigr)$\hyp{}semimonotone \optional{at $(\other{x}, \other{y}) \in \gph T$ where $\other{x} = (\other{x}_\A, \other{x}_\B)$ and $\other{y} = (\other{y}_\A,\other{y}_\B)$}.
    \end{enumerate}
\end{proposition}
In \Cref{def:semimonotonicity}, there is some freedom in selecting the matrices $\Mon$ and $\Com$, which might lead to a tradeoff between both. One particular class of operators for which this is true is the class of linear operators. This is summarized in the following proposition, which generalizes \cite[Prop. 5.1]{bauschke2021Generalized} for $\mu$\hyp{}monotone and $\rho$\hyp{}comonotone operators and \cite[Prop. 4.5]{evens2023convergence} for $(\mon,\com)$\hyp{}semimonotone operators.

\begin{proposition}[linear operator]\label{prop:semi:lin}
    Let $D \in \R^{n \times n}$ and let $\Mon, \Com\in\sym{n}$. Then, $D$ is $(\Mon, \Com)$\hyp{}semimonotone if and only if $\tfrac{1}{2}(D + D^\top) - \Mon - D^\top \Com D \succeq 0$.
\end{proposition}
Given a certain matrix $D$ and a desired semimonotonicity modulus $\Mon$, it might be difficult to determine whether there exists an $\Com$ satisfying $D^\top \Com D \preceq \tfrac{1}{2}(D + D^\top) - \Mon$, as this corresponds to solving a linear matrix inequality (LMI). 
The study of LMIs in general form has been extensively explored within the control and systems theory communities, leading to well-known results such as the Kalman--Yakubovich--Popov lemma, Finsler’s lemma and the (nonstrict) projection lemma \cite{boyd1994linear,helmersson1995methods,balakrishnan2003semidefinite,meijer2023nonstrict}.
In this work, we rely on a particular result for LMIs of the form $\tp \D X \D \preceq Y$, which is due to \cite{tian2010equalities,tian2013analytical} and relies upon the classical results from 
\cite{penrose1955generalized,khatri1976hermitian,baksalary1989properties} for the linear matrix equality $D^\top X D = Y$.
\begin{proposition}[symmetric solution of \(\tp \D X \D \preceq Y\)]\label{lem:LMI:solution}
    Let \(\D \in \R^{m \times n}\) and \(Y \in \sym{n}\).
    Then, the following hold.
    \begin{enumerate}
        \item \label{it:DXD:existence} The set of solutions \(C \dfn \set{X \in \sym{m}}[\tp \D X \D \preceq Y]\) is nonempty if and only if 
        \begin{align*}
            \proj_{\ker{\D}} Y \proj_{\ker{\D}} \succeq 0
            \quad \text{and} \quad
            \rank(\proj_{\ker{\D}} Y \proj_{\ker{\D}}) = \rank(\proj_{\ker{\D}} Y).\numberthis\label{eq:LMI:condition}
        \end{align*}
        \item \label{it:DXD:optimal} If \eqref{eq:LMI:condition} holds, then
            \(X^\star\in C\), where
            \begin{align*}
                X^\star
                    {}={}
                \begin{bmatrix}
                    0 & \I
                \end{bmatrix}
                \begin{bmatrix}
                    -Y & \tp \D \\ 
                    \D & 0
                \end{bmatrix}^\dagger
                \begin{bmatrix}
                    0 \\ \I
                \end{bmatrix}
                    {}={}
                (\D^\dagger)^\top \left(Y - Y \proj_{\ker{\D}}(\proj_{\ker{\D}} Y \proj_{\ker{\D}})^\dagger \proj_{\ker{\D}} Y\right) \D^\dagger.
                \numberthis\label{eq:LMI:solution:best:matrix}
            \end{align*}
        Moreover,
        \(\tp \D X \D \preceq \tp \D X^\star \D \preceq Y\) for all \(X \in C\). 
        \item \label{it:DXD:consistent} If the matrix equation \(\tp \D X \D = Y\) is consistent, i.e. if \(\range{Y} \subseteq \range{\tp \D}\),
        then \(X^\star = (\D^\dagger)^\top Y \D^\dagger\) is the solution of \(\tp \D X \D = Y\) with minimal trace \(\trace X^2\).
    \end{enumerate}
\end{proposition}
Applying this result to \Cref{prop:semi:lin}, the following corollary for linear operators is obtained.
\begin{corollary}[linear operator]\label{cor:semi:lin}
    Let $D \in \R^{n \times n}$ and $\Mon\in\sym{n}$.
    Then, the following hold.
    \begin{enumerate}
        \item \label{prop:semi:lin:existence} There exists $\Com \in \sym{n}$ such that $D$ is \((\Mon, \Com)\)\hyp{}semimonotone if and only if
        \begin{align*}
            \proj_{\ker{D}} \Mon \proj_{\ker{D}} \preceq 0
            \ifjota
\else
                \quad\text{and}\quad
            \fi
            \rank(\proj_{\ker{D}} \Mon \proj_{\ker{D}}) = \rank\bigl(\proj_{\ker{D}} (\tfrac12 D - \Mon)\bigr).\numberthis\label{eq:LMI:condition:lin}
        \end{align*}
        \item \label{prop:semi:lin:optimal} If \eqref{eq:LMI:condition:lin} holds, then
        $D$ is \((\Mon, \Com^\star)\)\hyp{}semimonotone,
        where
        \begin{align*}
            \Com^\star
                {}={}&
            \begin{bmatrix}
                0 & \I
            \end{bmatrix}
            \begin{bmatrix}
                \Mon - \tfrac{1}{2}(D + D^\top) & \tp D \\ 
                D & 0
            \end{bmatrix}^\dagger
            \begin{bmatrix}
                0 \\ \I
            \end{bmatrix}.
        \end{align*}
        In particular, when $D$ is either symmetric or skew-symmetric, it holds that
        \begin{align*}
            \Com^\star
                {}={}&
            2(D + D^\top)^\dagger - {D^\dagger}^\top \Mon \D^\dagger + (\D^\dagger)^\top \Mon \proj_{\ker{\D}}(\proj_{\ker{\D}} \Mon \proj_{\ker{\D}})^\dagger \proj_{\ker{\D}} \Mon\D^\dagger.
            \ifjota\else
                \numberthis\label{eq:LMI:solution:best:matrix:lin:sym}
            \fi
        \end{align*}
    \end{enumerate}
\end{corollary}
Note that $\Com^\star$ can be seen as an optimal choice for $\Com$, as it solves the LMI from \cref{prop:semi:lin} as tightly as possible.
A second consequence of \Cref{lem:LMI:solution} is the following result on operators of the form \(\D \F \D^\top\).

\begin{corollary}[semimonotonicity of \(\D \F \D^\top\)]\label{lem:calculus:composition:1}
    Let 
    \(
        \D \in \R^{n\times m}
    \),
    \(
        \Mon, Y \in \sym{m}
    \)
    and let operator \(\F : \R^m \rightrightarrows \R^m\) be \((\Mon,Y)\)\hyp{}semimonotone \optional{at \((\D^\top\other{x}, \other{y}) \in \gph \F\)}.
    If \eqref{eq:LMI:condition} holds for $\D$ and $\F$, 
    then \(\D \F \D^\top\) is \((\D \Mon \D^\top,X^\star)\)\hyp{}semimonotone \optional{at \((\other{x}, \D \other{y})\)}, where \(X^\star\) is   
    given by \eqref{eq:LMI:solution:best:matrix}.
\end{corollary}
Leveraging the previous result for the semimonotonicity of $\D \F \D^\top$, the semimonotonicity of the sum and parallel sum of two semimonotone operators is investigated next. 
Recall that the parallel sum of two operators $\A,\B : \R^n\rightrightarrows\R^n$ is defined as 
\(
    \A \Box \B \dfn (\A^{-1} + \B^{-1})^{-1}
\).
We shall also use the parallel sum between symmetric matrices $X, Y \in \sym{n}$ as defined below.

\begin{definition}[parallel sum between symmetric matrices]
    Let $X, Y \in \sym{n}$. We say that $X$ and $Y$ are parallel summable if 
    \(
        \range{X} \subseteq \range{X+Y}
    \)
    or equivalently 
    \(
        \range{Y} \subseteq \range{X+Y}
    \). 
    For parallel summable matrices $X$ and $Y$, their parallel sum is defined as 
    \cite[Cor. 9.2.5]{mitra2010Matrix}
    \begin{align*}
        X \Box Y
            \dfn
        X(X+Y)^\dagger Y
            {}={}
        Y(X+Y)^\dagger X
            {}={}
        X - X(X+Y)^\dagger X
            {}={}
        Y - Y(X+Y)^\dagger Y.
        \numberthis\label{eq:def:parsum:prop}
    \end{align*}
    For scalars $\alpha,\beta \in \R$, we say that $\alpha$ and $\beta$ are parallel summable if either $\alpha = \beta = 0$ or $\alpha + \beta \neq 0$ and their parallel sum is defined as 
    \[
        \alpha \Box \beta 
            {}\coloneqq{}
        \alpha(\alpha+\beta)^\dagger \beta
            {}={}
        \begin{cases}
            0, & \qquad\text{if } \alpha = \beta = 0,\\
            \frac{\alpha\beta}{\alpha+\beta}, & \qquad\text{otherwise}.
        \end{cases}
    \]
\end{definition}

In addition. consider the following set, which will be referred to as the effective domain of the parallel sum between two symmetric matrices.

\begin{definition}[effective domain of parallel sum]\label{def:parsum:domain}
    The set 
    \begin{align}\label{eq:parsum:dom}
        \parsumset \dfn \set{(X,Y) \in \mathbb{S}^n\times \mathbb{S}^n}[X+Y\succeq 0,\; X\text{ and }Y\text{ are parallel summable}]
    \end{align}
    is the effective domain of the parallel sum between two symmetric (possibly indefinite) matrices.
    Let \(X = \alpha \I_n\) and \(Y = \beta \I_n\) where $\alpha, \beta \in \R$. Then, $(X,Y) \in \parsumset$ reduces to
    \ifjota
\else
        $$
            (\alpha,\beta) \in \parsumset = \set{(\alpha, \beta)}[\alpha+\beta > 0 \text{ or } \alpha = \beta = 0].
        $$
    \fi
\end{definition}
In the upcoming proposition,
it is shown that the sum and parallel sum of two semimonotone operators are also semimonotone operators, if the involved semimonotonicity matrices belong to the effective domain of the parallel sum. 
This result generalizes \cite[Prop. 4.7]{evens2023convergence} for the sum of two $(\mon,\com)$-semimonotone operators.

\begin{proposition}[sum and parallel sum]\label{prop:calculus:sum+parsum}
    Let operator $A : \R^n \rightrightarrows \R^n$ be $(\MonA,\ComA)$\hyp{}semimonotone \optional{at $(\other{x}_A, \other{y}_A) \in \gph A$} and operator $B : \R^n \rightrightarrows \R^n$ be $(\MonB,\ComB)$\hyp{}semimonotone \optional{at $(\other{x}_B, \other{y}_B) \in \gph \B$}.
    \begin{enumerate}
        \item\label{prop:calculus:sum} If $(\ComA,\ComB) \in \parsumset$ \optional{and $\other{x}_A = \other{x}_B \eqqcolon \other{x}$}, then $A+B$ is $(\MonA+\MonB,\ComA \Box \ComB)$\hyp{}semimonotone \optional{at $(\other{x}, \other{y}_A+\other{y}_B)$}.
        \item\label{prop:calculus:parsum} If $(\MonA,\MonB) \in \parsumset$ \optional{and $\other{y}_A = \other{y}_B \eqqcolon \other{y}$}, then $A \Box B$ is $(\MonA \Box \MonB, \ComA+\ComB)$\hyp{}semimonotone \optional{at $(\other{x}_A+\other{x}_B, \other{y})$}.
    \end{enumerate}
\end{proposition}
When one of the two involved operators is linear, more precise statements for the resulting semimonotonicity matrices can be derived. For instance, consider the following lemma for the sum of a semimonotone operator and a skew-symmetric matrix.
This result will be used later in \Cref{thm:calculus:primaldual} for analyzing the primal-dual operator $\Tpd$.

\begin{lemma}[sum with skew-symmetric matrix]\label{lem:calculus:sum:linear}
    Let $\D \in \R^{n\times n}$ be a skew-symmetric matrix and operator $\F : \R^n \rightrightarrows \R^n$ be $(\tp \D\Mon \D,\Com + \Com^\prime)$\hyp{}semimonotone \optional{at $(\other{x}, \other{y}) \in \gph \F$}, where
    \(
        \range{\Com'} \subseteq \ker{\D}
    \)
    and \((\Mon, \Com) \in \parsumset\). Then, $\F + \D$ is $(0,\Com^\prime + \Mon \Box \Com)$\hyp{}semimonotone \optional{at $(\other{x}, \other{y} + D \other{x})$}.
\end{lemma}

The calculus rules presented in this section can be used both to generate and to identify semimonotone operators. 
We remark that it is also possible to confirm semimonotonicity of operators either directly by definition or numerically through their \emph{scaled relative graph}, as detailed in \cite{quan2024scaled}.
    
\section{Chambolle--Pock for semimonotone operators}\label{sec:CP:semi}

In \Cref{sec:CP}, convergence of \ref{eq:CP} was established under an oblique weak Minty assumption on the underlying primal-dual operator.
 This section provides a set of sufficient conditions for the convergence of \ref{eq:CP} for composite inclusion problems involving semimonotone operators. 

\subsection{Existence of oblique weak Minty solutions}
The main tool for establishing simplified conditions for \ref{eq:CP} for semimonotone operators is the following calculus rule, which connects the semimonotonicity of the individual operators \(A\) and \(B\) to the existence of $\DRSRho$\hyp{}oblique weak Minty solutions of the primal-dual operator $\Tpd$.

\begin{theorem}[primal-dual operator]\label{thm:calculus:primaldual}
    In the primal-dual inclusion~\eqref{eq:primaldual}, suppose that there exists a 
    nonempty set \(\pazocal{S}^\star \subseteq \zer \Tpd\) 
    and matrices $\ComA, \ComA', \ComB \in \sym{n}$ and $\MonA, \MonB, \MonB' \in \sym{m}$
    such that for every $z^\star = (x^\star, y^\star) \in \pazocal{S}^\star$ the following hold.
    \begin{enumerate}
        \item\label{eq:calculus:primaldual:range}
        $
            \range{\ComA} \subseteq \rge{\tp L} 
        $,
        $
            \range{\ComA'} \subseteq \ker{L} 
        $,
        $
            \range{\MonB} \subseteq \rge{L}
        $
        and
        $
            \range{\MonB'} \subseteq \ker{\tp L}.
        $
        \item\label{eq:calculus:primaldual:parsum}
        $
            (\MonA, \MonB) \in \parsumset
        $
        and
        $
            (\ComA, \ComB) \in \parsumset
        $.
        \item Operator $A$ is $(\tp L\MonA L,\ComA + \ComA')$\hyp{}semimonotone at $(x^\star, -L^\top y^\star) \in \gph A$.
        \item Operator $B$ is $(\MonB + \MonB',L \ComB \tp L)$\hyp{}semimonotone at $(L x^\star, y^\star) \in \gph \B$.
    \end{enumerate}
    Then, 
    $\Tpd$
    has
    $((\ComA \Box \ComB + \ComA') \oplus (\MonA \Box \MonB + \MonB'))$\hyp{}oblique weak Minty solutions at $\pazocal{S}^\star$.
    \begin{proof}
        Let $(x^\star, y^\star) \in \pazocal{S}^\star$ and decompose the primal-dual operator as $\Tpd = T + D$, where 
        $T \dfn A \times B^{-1}$
        and $D(x,y) \dfn (L^\top y, -L x)$.
        By \Cref{prop:calculus:inverse}, the inverse operator
        $B^{-1}$ is $(L \ComB \tp L, \MonB + \MonB')$\hyp{}semimonotone at $(y^\star,Lx^\star) \in \gph B^{-1}$.
        Therefore, it follows from \Cref{lem:calculus:cartesian} for the Cartesian product of two semimonotone operators that 
        $T$ is $\bigl(D^\top (\MonA \oplus \ComB)D,(\ComA + \ComA') \oplus (\MonB + \MonB')\bigr)$\hyp{}semimonotone at $\bigl((x^\star, y^\star), (-L^\top y^\star, Lx^\star)\bigr) \in \gph T$.
        Having established the semimonotonicity of $T$,
        it follows from \Cref{lem:calculus:sum:linear} and skew-symmetry of $D$ that 
        $
            \Tpd = T + D
        $   
        is
        $\bigl(0,(\ComA \Box \ComB + \ComA') \oplus (\MonA \Box \MonB + \MonB')\bigr)$\hyp{}semimonotone at $((x^\star, y^\star), 0)$.
        The claim then follows by \cref{def:SWMVI}.
    \end{proof}
\end{theorem}
Suppose that the underlying assumptions from \Cref{thm:calculus:primaldual} hold, in which case $\rge{\ComA \Box \ComB} \subseteq \rge{\ComA} \subseteq \rge{\tp L}$ and $\rge{\MonA \Box \MonB} \subseteq \rge{\MonB} \subseteq \rge{L}$.
Then, by virtue of \eqref{eq:projXY} and the particular form of 
\(
    \DRSRho = \bigl((\ComA \Box \ComB + \ComA') \oplus (\MonA \Box \MonB + \MonB')\bigr)
\)
from \Cref{thm:calculus:primaldual},
the primal-dual operator $\Tpd$ has $\DRSRho$\hyp{}oblique weak Minty solutions at \(\pazocal{S}^\star\), with \(\DRSRho\) given by \eqref{eq:WMSDRS} and 
\ifjota
\else
    \begin{equation}\label{eq:calculus:primaldual:params}
        \DRSrhocom = \lambda_{\rm min}(\VL^\top(\ComA \Box \ComB)\VL),
        \;\;
        \DRSrhomon = \lambda_{\rm min}(\UL^\top(\MonA \Box \MonB)\UL),
        \;\;
        \DRSrhocom' = \lambda_{\rm min}({\VL'}^\top\ComA'\VL'),
        \;\;
        \DRSrhomon' = \lambda_{\rm min}({\UL'}^\top\MonB'\UL'),
    \end{equation}
\fi
where $\VL, \VL', \UL, \UL'$ are defined as in \eqref{eq:L:svd}.
Hence, \Cref{ass:SWMVICP} holds if the parameters from \eqref{eq:calculus:primaldual:params}
satisfy condition \eqref{eq:existence:condition}.


As an implication of \Cref{thm:calculus:primaldual},
one can easily verify the underlying assumptions of \Cref{thm:CP:full}. 
For instance,
the following proposition applies to constrained, nonconvex quadratic optimization problems, with its proof provided in \Cref{proof:prop:constrainedQP}.

\begin{proposition}[constrained QP]\label{prop:constrainedQP}
    Consider the quadratic program
    \begin{align*}
        \minimize_{x\in\R^n} \; \tfrac12 x^\top Q x + q^\top x\quad
        \stt \quad L x \in C,\numberthis\label{eq:prop:QP}
    \end{align*}
    where 
    $Q \in \sym{n}$, $q \in \R^n$. $L \in \R^{m \times n}$ and $C \dfn \set{x\in \R^m}[l_i \leq x_i \leq u_i, \; i = 1, \dots, m]$, where $l, u \in \R^m$.
    The associated first-order optimality condition is given by
    $
        0 \in Ax + L^\top B L x,
    $
    where $A : x \mapsto Qx + q$ and $B \dfn N_{C}$. 
    Suppose that $L$ is either full column rank, or that
    $
        \proj_{\range{L^\top}} Q \proj_{\ker{L}} = 0
    $
    and
    $
        \proj_{\ker{L}} Q \proj_{\ker{L}} \succeq 0
    $.
    Then, the following hold.
    \begin{enumerate}
        \item \label{it:prop:constrainedQP:A} Operator $A$ is $(L^\top \MonA L, \ComA^\prime)$-semimonotone,
        where 
        $\MonA = {L^\dagger}^\top Q L^\dagger$
        and
        $
            \ComA^\prime = \proj_{\ker{L}} Q^\dagger \proj_{\ker{L}} \succeq 0
        $.
        \item \label{it:prop:constrainedQP:B}
        Let
        $
            (x^\star, y^\star) 
                \in
            \zer \Tpd
        $
        and
        $
            \MonB(y^\star)
                \coloneqq
            \diag\Bigl(\frac{|y^\star_1|}{u_1 - l_1}, \dots, \frac{|y^\star_n|}{u_n - l_n}\Bigr)
        $.
        Then, operator $B$ is $\Bigl(\MonB(y^\star), 0\Bigr)$-semimonotone at $(Lx^\star, y^\star) \in \gph B$.
    \end{enumerate}
    Moreover, if there exists a primal-dual pair
    $
        (x^\star, y^\star)  
            \in
        \zer \Tpd
    $
    satisfying
    \(
        \left(
            \MonA,
            \MonB(y^\star)
        \right) \in \parsumset,
    \)
    the following hold.
    \begin{enumerate}[resume]
        \item \label{it:prop:constrainedQP:Tpd}
        Let $\VL'$ and $\UL$ be as in \eqref{eq:L:svd}.
        The pair
        $
            (x^\star, y^\star) 
        $
        is a $\DRSRho$\hyp{}oblique weak Minty solution of $\Tpd$
        with $\DRSRho$ given by \eqref{eq:WMSDRS}, where 
        \(
            \DRSrhocom = 0
        \),
        \(
            \DRSrhomon = \lambda_{\rm min}\left(\UL^\top\left(
                \MonA
                \Box 
                \MonB(y^\star)
            \right)\UL\right)
        \),
        \[    
            \DRSrhocom^\prime 
                {}={}
            \begin{cases}
                \infty, &\;\; \text{if } \rank L = n,\\
                \lambda_{\rm min}\Bigl({\VL'}^\top
                        Q^\dagger
                    \VL'\Bigr) \geq 0, &\;\; \text{otherwise.}
            \end{cases}
            \quad\text{and}\quad
            \DRSrhomon^\prime 
                {}={}
            \begin{cases}
                \infty, &\;\; \text{if } \rank L = m,\\
                0, &\;\; \text{otherwise,}
            \end{cases}
        \]
        \item \label{it:prop:constrainedQP:convergence}
        Let the stepsizes $\gamma, \tau$ and fixed relaxation parameter $\lambda$ satisfy
        \begin{align*}
            \gamma \in 
            \bigl(0, \tfrac{1}{[\DRSrhomon]_-\nrm{L}^2}\bigr),
            \;
            \tau \in 
            \bigl([\DRSrhomon]_-, \tfrac{1}{\gamma\nrm{L}^2}\bigr),
            \;
            \lambda \in 
            \left(0, 2(1 + \min\{-\tfrac{1}{\tau} [\DRSrhomon]_-, \tfrac{2}{\gamma} \DRSrhocom^\prime\})\right).
            \numberthis\label{eq:QP:convergence:stepsizes}
        \end{align*}
        Then, \cref{ass:CP,ass:CP:stepsize:rule,ass:CP:relaxation:rule} hold, the resolvents $J_{\gamma A}$ and $J_{\tau B^{-1}}$ are (single-valued) continuous and the convergence properties of \ref{eq:CP} follow from \Cref{thm:CP:full}.
    \end{enumerate}
\end{proposition}

Notably, \Cref{prop:constrainedQP} covers examples of nonconvex quadratic programs for which the matrix $L$ is rank-deficient.
For instance, consider the following two examples with an indefinite $Q$ matrix, where $L$ is full row rank in the first example, while it is neither full row nor full column rank in the second.
The proofs are deferred to \Cref{proof:ex:constrainedQP:Semi}.
\begin{example}\label{ex:constrainedQP:Semi}
    Consider the QP from \eqref{eq:prop:QP},
    where 
    \ifjota
\else
        $
            Q
                {}={}
            \diag(1, -1, 2)
        $,
        $
            q
                {}={}
            \begin{bsmallmatrix}
                -1&
                1&
                -1
            \end{bsmallmatrix}^\top
        $,
        $
        L = 
        \begin{bsmallmatrix}
            1 & \nicefrac14 & 0\\
            0 & 1 & 0
        \end{bsmallmatrix}
        $
        and
        $
            C \dfn \set{x\in \R^2}[2 \leq x_i \leq 4, i = 1, 2]
        $.
    \fi
    Let operators $A$ and $B$ be defined as in \Cref{prop:constrainedQP}.
    Then, the global minimizer of \eqref{eq:prop:QP} is given by         
    $
        x^\star 
            = 
        \begin{bsmallmatrix}
            1 &
            4 &
            \nicefrac12
        \end{bsmallmatrix}^\top
    $
    and the following hold.
    \begin{enumerate}
        \item \label{it:QP:A:full:Semi} 
        \ifjota\else
            Operator 
        \fi
        $A$ is 
        $
            \bigl(\diag(1, -1, 0), \diag(0,0,\nicefrac12)\bigr)
        $\hyp{}semimonotone.
        \item \label{it:QP:B:full:Semi} 
        \ifjota\else
            Operator 
        \fi
        $B$ is 
        $
            \bigl(\diag(0, \nicefrac32), 0\bigr)
        $\hyp{}semimonotone at 
        $
            (Lx^\star, -{L^\dagger}^\top A x^\star) 
                = 
            \left(
                \begin{bsmallmatrix}
                    2 \\
                    4
                \end{bsmallmatrix}
                , 
                \begin{bsmallmatrix}
                    0 \\
                    3
                \end{bsmallmatrix}
            \right) 
            \in \gph B
        $.
        \item \label{it:QP:Tpd:full:Semi}
        The
        \ifjota\else
            primal-dual
        \fi
        pair
        \((x^\star, -{L^\dagger}^\top A x^\star) \in \zer \Tpd\)
        is a $\DRSRho$\hyp{}oblique weak Minty solution of $\Tpd$
        with $\DRSRho$ given by \eqref{eq:WMSDRS}, where 
        $
            \DRSrhocom = 0,
            \DRSrhomon = -3,
            \DRSrhocom^\prime = \nicefrac12
        $
        and
        $
            \DRSrhomon^\prime = \infty
        $.
        \item \label{it:QP:CP:full:Semi}
        The sequence $\seq{z^k} = \seq{x^k,y^k}$ generated by \ref{eq:CP} with fixed relaxation parameter $\lambda$ converges for 
        \ifjota
\else
            $$
                \gamma \in \Bigl(0, \tfrac{-1}{\DRSrhomon\nrm{L}^2}\Bigr) \approx (0, 0.26),
                \quad
                \tau \in \Bigl(-\DRSrhomon, \tfrac{1}{\gamma\nrm{L}^2}\Bigr] \approx (3, \tfrac{0.779}{\gamma}],
                \quad
                \lambda \in (0, 2 + \tfrac{2}{\tau} \DRSrhomon) = (0, 2 - \tfrac{6}{\tau}),
            $$
            where we used that $\nrm{L}^2 = \tfrac{33 + \sqrt{65}}{32} \approx 1.28$. \qedhere
        \fi
    \end{enumerate}
\end{example}

\begin{example}\label{ex:constrainedQP:semi}
    Consider the QP from \eqref{eq:prop:QP}, where
    $
        Q
            {}={}
        \diag(-3, -2, 1)
    $,
    $
        q
            {}={}
        \begin{bsmallmatrix}
            0&
            1&
            0
        \end{bsmallmatrix}^\top
    $,
    $
        L = 
        \begin{bsmallmatrix}
            1 & 0 & 0\\
            1 & 1 & 0\\
            1 & -1 & 0
        \end{bsmallmatrix}
    $
    and
    $
        C = \set{x\in \R^2}[\nicefrac12 \leq x_i \leq 1, \; i = 1, 2, 3].
    $
    Let operators $A$ and $B$ be defined as in \Cref{prop:constrainedQP}.
    Then, the global minimizer is given by         
    $
        x^\star 
            = 
        \begin{bsmallmatrix}
            1 &
            0 &
            0
        \end{bsmallmatrix}^\top
    $
    and the following hold.
    \begin{enumerate}
        \item \label{it:QP:A:semi}
        \ifjota\else
            Operator 
        \fi
        $A$ is 
        \ifjota
\else
            $
                \bigl(\monA L^\top L, \diag(0,0,1)\bigr)
            $\hyp{}semimonotone where $\monA = -1$.
        \fi
        \item \label{it:QP:B:semi}
        \ifjota\else
            Operator 
        \fi
        $B$ is 
        \ifjota
\else
            $
                \bigl(\monB, 0\bigr)
            $%
        \fi
        \hyp{}semimonotone at 
        $
            (Lx^\star, -{L^\dagger}^\top A x^\star) 
                = 
            \left(
                \begin{bsmallmatrix}
                    1 &
                    1 &
                    1
                \end{bsmallmatrix}^\top
                , 
                \begin{bsmallmatrix}
                    1 &
                    1 &
                    1
                \end{bsmallmatrix}^\top
            \right) 
            \in \gph B
        $\ifjota
\else
            \ where $\monB = 2$.
        \fi
        \item \label{it:QP:Tpd:full:semi}
        The
        \ifjota\else
            primal-dual
        \fi
        pair
        \((x^\star, -{L^\dagger}^\top A x^\star) \in \zer \Tpd\)
        is a $\DRSRho$\hyp{}oblique weak Minty solution of $\Tpd$
        with $\DRSRho$ given by \eqref{eq:WMSDRS}, where 
        $
            \DRSrhocom = 0,
            \DRSrhomon = -2,
            \DRSrhocom^\prime = 1
        $
        and
        $
            \DRSrhomon^\prime
                = 
            0
        $.
        \item \label{it:QP:CP:semi}
        The sequence $\seq{z^k} = \seq{x^k,y^k}$ generated by \ref{eq:CP} with fixed relaxation parameter $\lambda$ converges for 
        \ifjota
\else
            $$
                \gamma \in \bigl(0, \tfrac{-1}{(\monA \Box \monB)\nrm{L}^2}\bigr) = (0, \tfrac16),
                \quad
                \tau \in \bigl(-(\monA \Box \monB), \tfrac{1}{\gamma\nrm{L}^2}\bigr] = (2, \tfrac{1}{3\gamma}],
                \quad
                \lambda \in \bigl(0, 2 + \tfrac{2}{\tau} (\monA \Box \monB)\bigr) = (0, 2 - \tfrac{4}{\tau}).\qedhere
            $$
        \fi
    \end{enumerate}
\end{example}


\subsection{Sufficient conditions for convergence of CPA}
\Cref{eq:calculus:primaldual:range} imposes range conditions on the semimonotonicity matrices of $A$ and $B$.
In this subsection, it is shown that this can be achieved by imposing a certain structure on these matrices. 
In particular, consider the following set of assumptions.

\begin{assumption}\label{ass:SWMVICPsemimonotone}
    In problem~\eqref{prob:composite}, suppose that \(\zer \Tpd\) is nonempty and that there exists a nonempty set \(\pazocal{S}^\star\subseteq \zer \Tpd\) such that for every $z^\star = (x^\star, y^\star) \in \pazocal{S}^\star$ it holds that $\A$ is $(\monA \tp L L, \comA \I_n)$\hyp{}semimonotone at $(x^\star, -L^\top y^\star) \in \gph A$, $\B$ is
    $(\monB \I_m, \comB L \tp L)$\hyp{}semimonotone at $(L x^\star, y^\star) \in \gph B$
    and the semimonotonicity moduli $(\monA, \monB, \comA, \comB) \in \R^4$ satisfy either one of the following conditions.
    \begin{enumerate}
        \item\label{ass:CP:semi:case:1} (either) \(\monA = \monB = 0\) and $\comA = \comB = 0$.
        \item\label{ass:CP:semi:case:2} (or) $\monA + \monB > 0$ and $\comA = \comB = 0$.
        \item\label{ass:CP:semi:case:3} (or) $\comA + \comB > 0$ and \(\monA = \monB = 0\).
        \item\label{ass:CP:semi:case:nonmon} (or) $\monA + \monB > 0$, $\comA + \comB > 0$ and 
        $
            [\monA \Box \monB]_-[\comA \Box \comB]_- < \tfrac{1}{4\nrm{L}^2}
        $.
    \end{enumerate}
\end{assumption}
\Cref{ass:SWMVICPsemimonotone} provides a set of sufficient conditions on the operators $A$ and $B$ so that the primal-dual operator admits oblique weak Minty solutions.
This key result is stated in the following corollary.

\begin{corollary}\label{prop:semi:WMVI}
    Suppose that \Cref{ass:SWMVICPsemimonotone} holds. Then, 
    the primal-dual operator $\Tpd$ has $\DRSRho$\hyp{}oblique weak Minty solutions at \(\pazocal{S}^\star\), where $\DRSRho$ is given by \eqref{eq:WMSDRS} and
    \ifjota
\else
        \begin{align*}
            \DRSrhocom = \comA \Box \comB,
            \quad
            \DRSrhomon = \monA \Box \monB,
            \quad
            \DRSrhocom^\prime = 
            \begin{cases}%
                \infty,  &\quad \text{if } \rank L = n,\\
                \comA, &\quad \text{if } \rank L < n,
            \end{cases}, 
            \quad 
            \DRSrhomon^\prime = 
            \begin{cases}
                \infty,  &\quad \text{if } \rank L = m,\\
                \monB, &\quad \text{if } \rank L < m.
            \end{cases}
            \numberthis\label{eq:DRSRho:semimonotonicity}
        \end{align*}%
    \fi
    \begin{proof}
        Note that
        $
            (\monA, \monB) \in \parsumset
        $
        and
        $
            (\comA,\comB) \in \parsumset
        $
        by \Cref{ass:SWMVICPsemimonotone} and that
        $
            \comA \I_n = \comA \proj_{\rge{\tp L}} + \DRSrhocom^\prime \proj_{\ker{L}}
        $
        and
        $
            \monB \I_m = \monB \proj_{\rge{L}} + \DRSrhomon^\prime \proj_{\ker{\tp L}}
        $.
        Therefore, the claim follows directly from \Cref{thm:calculus:primaldual}, with $\MonA = \monA \proj_{\rge{L}}$, $\ComA = \comA \proj_{\rge{\tp L}}$, $\ComA^\prime = \DRSrhocom^\prime \proj_{\ker{L}}$, $\MonB = \monB \proj_{\rge{L}}$, $\MonB^\prime = \DRSrhomon^\prime \proj_{\ker{\tp L}}$ and $\ComB = \comB \proj_{\rge{\tp L}}$.
    \end{proof}
\end{corollary}
\setlength{\maxmin}{\widthof{$\min\Bigl\{1 + \Delta_{\gamma,\tau} - \theta_{\gamma \tau}(\sigma_{d}), \eta^\prime\Bigr\}$}-\widthof{$1 + 2\Delta_{\gamma,\tau}$}}
\setlength{\maxmincond}{\widthof{$\text{if } \max\{\monA\monB, \comA\comB\} \leq 0$}-\widthof{$\text{otherwise}$}}

Leveraging this key result, the convergence of \ref{eq:CP} for semimonotone operators is analyzed in \Cref{cor:CP:semi}, for which the stepsize rule corresponds to a simple look-up table.
The proof is deferred to \Cref{proof:cor:CP:semi}.

\begin{stepsize}\label{ass:CP:stepsize:rule:semi}
    The stepsizes $\gamma$ and $\tau$ satisfy the bounds provided in \Cref{tab:CP:stepsize:rule:semi}, where
    $\gamma_{\rm min}$, $\gamma_{\rm max}$ and $\tau_{\rm min}(\gamma)$ are defined as in
    \eqref{eq:stepsizes:CP:bounds:gam} with $\DRSrhomon = \monA \Box \monB$ and $\DRSrhocom = \comA \Box \comB$.
    \ifjota
\else
        \vspace{-0.1cm}
    \fi
    \begin{figure}[H]
        \centering
        \includetikz{Tables/stepsizes-CP-semi}
        \caption{Range of the stepsizes $\gamma$ and $\tau$ for \ref{eq:CP} involving semimonotone operators.}
        \label{tab:CP:stepsize:rule:semi}
    \end{figure}
\end{stepsize}

\begin{relaxation}\label{ass:CP:relaxation:rule:semi}
    Let $\eta^\prime$ be defined as in \Cref{tab:CP:degen:semi} and define
    \[
        \Delta_{\gamma,\tau} \dfn \tfrac{1}{2\gamma}(\comA \Box \comB) + \tfrac{1}{2\tau}(\monA \Box \monB)
        \quad\text{and}\quad
        \theta_{\gamma\tau}(\sigma) 
            {}\dfn{} 
        \sqrt{\Delta_{\gamma,\tau}^2 + (\monA \Box \monB)(\comA \Box \comB)\sigma^2}.
    \]
    The sequence $\seq{\lambda_k}$ satisfies $\lambda_k\in\bigl(0,2(1 + \etamin)\bigr)$ and $\liminf_{k\to\infty}\lambda_k\bigl(2(1 + \etamin)-\lambda_k\bigr)>0$, where
    \begin{align*}
        \etamin
            {}\dfn{}
        \begin{cases}
            \left\{
            \begin{array}{ll}
                \Delta_{\gamma,\tau} - \theta_{\gamma \tau}(\|L\|),
                &\quad \text{if } \monA\monB\comA\comB \geq 0\\
                \min\Bigl\{\Delta_{\gamma,\tau} - 
                \theta_{\gamma \tau}(\sigma_{d}), \eta^\prime\Bigr\},
                &\quad \text{otherwise}\hspace*{\maxmincond}
            \end{array}
            \right\}
            & \quad \text{if } \gamma\tau < \tfrac{1}{\nrm{L}^2},\\
            \left\{
            \begin{array}{ll}
                2\Delta_{\gamma,\tau},\hspace*{\maxmin}
                &\quad \text{if } \max\{\monA\monB, \comA\comB\} \leq 0\\
                \min\Bigl\{2\Delta_{\gamma,\tau}, \eta^\prime\Bigr\},
                &\quad \text{otherwise}\hspace*{\maxmincond}
            \end{array}
            \right\}
            & \quad \text{if } \gamma\tau = \tfrac{1}{\nrm{L}^2} \text{ and } d = 1,\\
            \left\{
                \begin{array}{ll}
                    2\Delta_{\gamma,\tau},\hspace*{\maxmin}
                    &\quad \text{if } \max\{\monA\monB, \comA\comB\} \leq 0\\
                    \Delta_{\gamma,\tau} - \theta_{\gamma \tau}(\sigma_2)
                    ,
                    &\quad \text{if } \min\{\monA\monB, \comA\comB\} > 0\\
                    \min\Bigl\{\Delta_{\gamma,\tau} - \theta_{\gamma \tau}(\sigma_{d}), \eta^\prime\Bigr\}
                    ,
                    &\quad \text{otherwise}
                \end{array}
                \right\}
                & \quad \text{if } \gamma\tau = \tfrac1{\nrm{L}^2} \text{ and } d > 1.
        \end{cases}
    \end{align*} 
    \ifjota
\else
        \vspace{-0.3cm}
    \fi
    \begin{figure}[H]
        \centering
        \includetikz{Tables/relaxation-CP-semi}
        \caption{Definition of $\eta^\prime$ in \cref{ass:CP:relaxation:rule:semi}.}
        \label{tab:CP:degen:semi}
    \end{figure}
\end{relaxation}

\begin{corollary}[convergence of \ref{eq:CP} under semimonotonicity]\label{cor:CP:semi}
    Suppose that \Cref{ass:CP:1}, \Cref{ass:CP:2} and \Cref{ass:SWMVICPsemimonotone} hold, that $\gamma$ and $\tau$ are selected according to \Cref{ass:CP:stepsize:rule:semi} 
    and that the relaxation sequence $\seq{\lambda_k}$ is selected according to \Cref{ass:CP:relaxation:rule:semi}.
    Then, all the claims of \cref{thm:CP:full} hold.
\end{corollary}

\Cref{cor:CP:semi} serves as a universal framework for analyzing the convergence of \ref{eq:CP}, both in monotone and nonmonotone settings. Notably, it encompasses and extends many of the existing results in literature. Several examples are provided below.

\begin{remark}[connection to existing theory]\label{rem:CP:cor:connections-literature}
    Case (i) of \cref{ass:SWMVICPsemimonotone} can be interpreted as a relaxation of the classical monotonicity assumption for \ref{eq:CP} \cite{chambolle2011firstorder}, as it does not impose monotonicity between any two pairs in the graph of $A$ and $B$. In case (ii) of \cref{ass:SWMVICPsemimonotone}, a monotone problem is split in a nonmonotone fashion. In the optimization setting, this setting was studied in \cite{thomas2015primaldual}. To see this, let $g$ be a proper lsc $\mu_g$-convex function with $\mu_g > 0$ and $h$ be a proper lsc $\mu_h$-convex function. Then, $A = \partial g$ is $(\nicefrac{\mu_g}{\nrm{L}^2}L^\top L,0)$-semimonotone and $B = \partial h$ is $(\mu_h, 0)$-semimonotone. Then, \Cref{cor:CP:semi} requires that $\nicefrac{\mu_g}{\nrm{L}^2} + \mu_h > 0$, which matches \cite[Thm. 2.8]{thomas2015primaldual}. 
    Note that case (iii) of \cref{ass:SWMVICPsemimonotone} can be interpreted as the dual counterpart of case (ii), as the assumptions of the latter hold for $A$ and $B$ in the primal inclusion problem if and only if the assumptions of the former hold for $B^{-1}$ and $A^{-1}$ in the dual one (see \eqref{eq:dual}).
    Up to the knowledge of the authors, no particular instances of case (iv) of \cref{ass:SWMVICPsemimonotone} have been covered in literature, even in the minimization setting.
\end{remark}


In what follows, \Cref{cor:CP:semi} will be applied to several examples previously discussed in this paper.
First, note that the convergence result from \Cref{it:QP:CP:semi} can be more conveniently derived as a consequence of \Cref{cor:CP:semi}, since \Cref{ass:SWMVICPsemimonotone} is satisfied for this problem with $\monA = -1$, $\monB = 2$, $\comA = 1$ and $\comB = 0$.
Next, we revisit the two theoretical examples of \Cref{sec:CP}, this time under the lens of semimonotonicity. 
First of all, consider the linear inclusion problem from \Cref{ex:CP:example5}, where the parameters $a, b, c, l$ are selected as in \Cref{it:CP:example5:4}. 

\begin{example}[saddle point problem (revisited)]\label{ex:CP:example5:revisited}
    Consider 
    inclusion problem \eqref{eq:ex5-Tp} with $a = 10$, 
    $
        b = -\tfrac14
    $,
    $c = 0$
    and $\ell = 2$. 
    Using \Cref{prop:semi:lin}, it follows that $A$ is $\Bigl(L^\top L, -\tfrac{1}{25}\I_n\Bigr)$\hyp{}semimonotone and $B$ is $\Bigl(-\tfrac{3}{10}\I_m, \tfrac{1}{5}L^\top L\Bigr)$\hyp{}semimonotone. 
    By \Cref{cor:CP:semi}, the sequence $\seq{z^k}$ generated by applying \ref{eq:CP} to \eqref{eq:ex5-Tp} with $\tau = \tfrac{1}{\gamma \nrm{L}^2}$ and fixed relaxation parameter $\lambda$
    converges for
    \ifjota
\else
        \[
            \gamma \in 
            \left(
                \gamma_{\rm min}, \gamma_{\rm max}
            \right) 
                {}\approx{}
            \left(
                0.055, 0.528
            \right)
            \quad\text{and}\quad
            \lambda 
                {}\in{} 
            \left(0, 2 - \tfrac1{10\gamma} - \tfrac{24\gamma}{7}\right).\qedhere
        \]
    \fi
\end{example}
The obtained range of stepsize parameters is only a subset of the tight range obtained in \Cref{it:CP:example5:4}.
However, this should not come as a surprise, since part of the information about operators $A$ and $B$ is lost by analyzing them under the lens of semimonotonicity. 
This is also observed in the second example.

\begin{example}[influence of singular values (revisited)]\label{ex:CP:second:revisited}
    Consider 
    the composite inclusion problem $0 \in Ax + L^\top B Lx$
    from \Cref{ex:CP:second}.
    It follows from \Cref{prop:semi:lin} that $A$ is $\Bigl(\tfrac12 L^\top L, \tfrac12 \I_n\Bigr)$\hyp{}semimonotone and $B$ is $\Bigl(\tfrac12 \I_n, \tfrac12 LL^\top\Bigr)$\hyp{}semimonotone.
    By \Cref{cor:CP:semi} and \Cref{rem:PDHG:sequences}, the sequences
    $(\projQ z^k)_{k\in\N}$ and $\seq{s^k}$ 
    generated by applying \ref{eq:CP} with $\gamma = \tau = 1$ and fixed relaxation parameter $\lambda$,
    converge to zero
    if $\lambda$ is selected according to \Cref{ass:CP:relaxation:rule:semi}, which reduces to
    \ifjota
\else
        \begin{align*}
            \lambda 
                {}\in{}&
            \left(
            0, 
            2\left(
                1 + \tfrac{1}{2\gamma}\DRSrhocom + \tfrac{1}{2\tau}\DRSrhomon - \theta_{\gamma \tau}(\max\{|\ell_{2}|, \hdots, |\ell_{n}|\})
            \right)
            \right)
                {}={}
            \Bigl(
            0, 
            \nicefrac52 - \nicefrac12\max\{|\ell_{2}|, \hdots, |\ell_{n}|\}
            \Bigr).\qedhere
        \end{align*}
    \fi
\end{example}
In \Cref{ex:CP:second:Tpd}, it was shown that $\Tpd$ has $(\tfrac12\I_n \oplus \tfrac12\I_n)$\hyp{}oblique weak Minty solutions at \(\zer \Tpd\).
On the other hand, here it is shown that $A$ is $(\tfrac12 L^\top L, \tfrac12 \I_n)$\hyp{}semimonotone and $B$ is $(\tfrac12 \I_n, \tfrac12 LL^\top)$\hyp{}semimonotone.
Applying \Cref{prop:semi:WMVI}, this implies that $\Tpd$ only has a $(\tfrac14\I_n \oplus \tfrac14\I_n)$\hyp{}oblique weak Minty solution at $(0_n, 0_n) = \zer \Tpd$.
By analyzing $A$ and $B$ under the lens of semimonotonicity, some additional looseness is inevitably introduced. More specifically, the information that $A$ and $B$ are linear and symmetric, that $A = B^{-1}$ and that $A = \nicefrac12 L^\top L + \nicefrac12 A^\top A$ are lost in this process.

Finally, for a more practical application of our theory, we refer the interested reader to \cite[Prop. V.2]{quan2024scaled}, where \Cref{cor:CP:semi} is applied to a nonmonotone common-emitter amplifier circuit.

  \section{Conclusion}\label{sec:conclusion}

In this work, convergence of the Chambolle--Pock algorithm (CPA) was established for a class of nonmonotone problems, characterized by an oblique weak Minty assumption on the associated primal-dual operator.
To facilitate the verification of this underlying assumption, a generalization of the class of semimonotone operators (see \cite{evens2023convergence}) was introduced, and sufficient conditions for the convergence of CPA were provided for inclusion problems involving operators belonging to this class.
When restricting to minimization problems, our results reveal that for certain problem classes no explicit rank or condition number restriction on the linear mapping is required.

It would be interesting to explore if stronger results can be obtained when the operators are known to be subdifferentials. Other future research directions include extensions to the setting where the preconditioning is indefinite,
allowing to cover the extended Chambolle--Pock stepsize range $\gamma\tau\nrm{L}^2 \leq \nicefrac43$ from \cite{yan2024improved,banert2023chambollepock}, 
as well as analyzing other splitting methods in nonmonotone settings.

\begin{appendix}
\ifjota\else
  \section{Supplementary results}\label{sec:supplementary}

\begin{theorem}\label{thm:CP:full:strong}
    Suppose that \cref{ass:CP} holds, with condition \eqref{eq:existence:condition} relaxed to
    \begin{align*}
        [\DRSrhocom]_-[\DRSrhomon]_- < \tfrac{1}{4\nrm{L}^2},
            \quad
        [\DRSrhocom^\prime]_-[\DRSrhomon^\prime]_- < \tfrac{1}{\nrm{L}^2},
            \quad
        [\DRSrhocom^\prime]_- 
            {}>{}
        \tfrac{\delta + \sqrt{\delta^2 - 4\DRSrhocom\DRSrhomon\nrm{L}^2}}{2[\DRSrhomon]_-\nrm{L}^2}
        \quad\text{and}\quad
        [\DRSrhomon^\prime]_-
            {}>{}
        \tfrac{\delta + \sqrt{\delta^2 - 4\DRSrhocom\DRSrhomon\nrm{L}^2}}{2[\DRSrhocom]_-\nrm{L}^2}.
        \numberthis\label{eq:existence:condition:strong}
    \end{align*}
    Suppose that
    $\gamma$ and $\tau$ are selected within the intersection of the bounds provided in \Cref{ass:CP:stepsize:rule} and the set
    \begin{align*}
        \set{(\gamma, \tau) \in \R^2}[
            \gamma \in \bigl(-\lbrack\DRSrhocom^\prime\rbrack_-, \tfrac{1}{-\lbrack\DRSrhomon^\prime\rbrack_-\nrm{L}^2}\bigr),
            \tau \in \bigl(-\lbrack\DRSrhomon^\prime\rbrack_-, \tfrac{1}{\gamma\nrm{L}^2}\bigr\rbrack].
            \numberthis\label{eq:proof:CP:etaprime:stepsizes:extra}
    \end{align*}
    Moreover, suppose that the relaxation sequence $\seq{\lambda_k}$ is selected according to \Cref{ass:CP:relaxation:rule}.
    Then, all the claims from \Cref{thm:CP:full} hold.
    \begin{proof}
        The proof is identical to the proof of \Cref{thm:CP:full}, except that in the item discussing \textit{$1 + \eta > 0$} the occurence of \eqref{eq:existence:condition} is replaced by \Cref{eq:existence:condition:strong} and that the item discussing \textit{$1 + \etamin > 0$} is replaced by the following.
        \begin{proofitemize}
            \item
            \textit{
                The set of pairs $(\gamma, \tau) \in \R^2_{++}$ satisfying $ \gamma\tau \in (0, \nicefrac{1}{\nrm{L}^2}]$ and $1 + \etamin > 0$ is given by 
                the the intersection of the bounds provided in \Cref{ass:CP:stepsize:rule} and the set \eqref{eq:proof:CP:etaprime:stepsizes:extra}, which is nonempty.
            }\newline
            With the understanding that $\DRSrhocom^\prime = \infty$ when $L$ is full column rank and $\DRSrhomon^\prime = \infty$ when $L$ is full row rank, as discussed below \Cref{ass:CP}, the set of pairs $(\gamma, \tau) \in \R^2_{++}$ satisfying $ \gamma\tau \in (0, \nicefrac{1}{\nrm{L}^2}]$ and $1 + \eta^\prime > 0$ is given by \eqref{eq:proof:CP:etaprime:stepsizes:extra}.
            Hence, it only remains to verify that the intersection of the bounds provided in \Cref{ass:CP:stepsize:rule} and the set \eqref{eq:proof:CP:etaprime:stepsizes:extra} is nonempty, which is ensured by the conditions on $\DRSrhocom^\prime$ and $\DRSrhomon^\prime$ provided in  \Cref{eq:existence:condition:strong}.\qedhere
        \end{proofitemize}
    \end{proof}
\end{theorem}
\fi

\section{Auxiliary results}\label{sec:auxiliary}

\begin{fact}[solution of quadratic inequality]\label{lem:quadratic}
    Let \(\DRSrhocom, \DRSrhomon\in\R\), \(\nrm{L} > 0, \sigma_d \in (0, \nrm{L}]\), let $\delta$ be defined as in \eqref{eq:def:delta} and let $\gamma_{\rm min}$ and $\gamma_{\rm max}$ be defined as in \cref{ass:CP:stepsize:rule}.
    Then, the following hold.
    \begin{enumerate}
        \item \label{it:lem:quadratic:existence} There exists a $\gamma > 0$ satisfying
        \(
            \DRSrhomon\nrm{L}^2\gamma^2 + \delta\gamma + \DRSrhocom > 0
        \)  
        \ifjota
\else
            if and only if
        \fi
        $
            [\DRSrhocom]_-[\DRSrhomon]_- < \tfrac{1}{4\nrm{L}^2}
        $.
        \item \label{it:lem:quadratic:solutions} 
        If
        $
            [\DRSrhocom]_-[\DRSrhomon]_- < \tfrac{1}{4\nrm{L}^2}
        $,
        then $\gamma > 0$ satisfies 
        \(
            \DRSrhomon\nrm{L}^2\gamma^2 + \delta\gamma + \DRSrhocom > 0
        \)
        if and only if \(\gamma\) lies within the (nonempty) stepsize interval provided in \Cref{tab:CP:stepsize:rule}.
    \end{enumerate}
\end{fact}

\begin{lemma}\label{cor:CP:stepsize:2}
    Let \(\DRSrhocom, \DRSrhomon\in\R\), \(\nrm{L} > 0\), \(\sigma_d \in (0, \nrm{L}]\) and for any $\sigma \in [\sigma_d, \nrm{L}]$ define the set
    \begin{align*}
        \Gamma(\sigma)
            {}\dfn{}&
        \left\{(\gamma, \tau) \in \R^2_{++} \; \middle| \; \begin{array}{l}
            \gamma\tau \in \bigl(0, \nicefrac{1}{\nrm{L}^2}\bigr] \text{ and } 1 + \tfrac{1}{2\gamma}\DRSrhocom + \tfrac{1}{2\tau}\DRSrhomon
                -
            \theta_{\gamma\tau}(\sigma)
                    >
                0
        \end{array}
            \right\},
    \end{align*}
    where $\theta_{\gamma\tau}(\cdot)$ is defined as in \eqref{eq:theta}. 
    Then, the following hold.
    \begin{enumerate}
        \item \label{it:cor:quadratic:CP:existence:pos} If $\min\{\DRSrhocom,\DRSrhomon\} \geq 0$, then, for any $\sigma \in [\sigma_d, \nrm{L}]$, stepsizes \((\gamma, \tau)\) belonging to the nonempty set $\Gamma(\sigma)$ are as in
        \Cref{tab:CP:stepsize:rule}.
        \item \label{it:cor:quadratic:CP:existence:mixed} If 
        $\min\{\DRSrhocom,\DRSrhomon\} < 0$ and
        $\max\{\DRSrhocom,\DRSrhomon\} \geq 0$, then, stepsizes \((\gamma, \tau)\) belonging to the nonempty set $\Gamma(\sigma_d)$ are as in
        \Cref{tab:CP:stepsize:rule}.
        \item \label{it:cor:quadratic:CP:existence:neg} If $\max\{\DRSrhocom,\DRSrhomon\} < 0$, then the set $\Gamma(\nrm{L})$
        is nonempty 
        if and only if
        $
            [\DRSrhocom]_-[\DRSrhomon]_- < \tfrac{1}{4\nrm{L}^2}
        $, in which case stepsizes \((\gamma, \tau)\) belonging to the nonempty set $\Gamma(\nrm{L})$ are as in
        \Cref{tab:CP:stepsize:rule}.
    \end{enumerate}
        \begin{proof}%
            Throughout the proof, we will use $\gamma_{\rm min}$, $\gamma_{\rm max}$ and $\tau_{\rm min}(\gamma)$ as defined as in \eqref{eq:stepsizes:CP:bounds:gam}.
            Let $\sigma \in [\sigma_d, \nrm{L}]$. Solving the square root inequality 
            $
                1 + \tfrac{1}{2\gamma}\DRSrhocom + \tfrac{1}{2\tau}\DRSrhomon > \theta_{\gamma\tau}(\sigma)
            $,
            it follows that
            \begin{align*}
                \Gamma(\sigma)
                    {}={}
                    \left\{(\gamma, \tau) \in \R^2_{++} \; \middle| \; 
                        \gamma\tau \in \bigl(0, \nicefrac{1}{\nrm{L}^2}\bigr],
                        \overbrace{1 + \tfrac{1}{2\gamma}\DRSrhocom + \tfrac{1}{2\tau}\DRSrhomon}^{\smash{\Gamma_1\, \dfn}} > 0
                        \ifjota
\else
                            \text{ and }
                        \fi
                        \overbrace{1 + \tfrac{1}{\gamma}\DRSrhocom + \tfrac{1}{\tau}\DRSrhomon + \tfrac{1}{\gamma\tau}\DRSrhocom\DRSrhomon(1-\gamma\tau\sigma^2)}^{\smash{\Gamma_2(\sigma)\, \dfn}} > 0
                    \right\}.
            \end{align*}
            Define $c_1(\sigma, \gamma) \dfn \gamma(1 - \DRSrhocom\DRSrhomon\sigma^2) + \DRSrhocom$ and $c_2(\gamma) \dfn \DRSrhomon(\gamma + \DRSrhocom)$, so that 
            \(
                \gamma\tau\Gamma_2(\sigma) 
                    {}={}
                c_1(\sigma, \gamma)\tau + c_2(\gamma)
            \).
            \begin{proofitemize}
                \item \ref{it:cor:quadratic:CP:existence:pos}: If $\min\{\DRSrhocom,\DRSrhomon\} \geq 0$, then $\Gamma_1 > 0$ and $\Gamma_2(\sigma) > 0$ since $1 - \gamma\tau \sigma^2 \geq 0$.
                Therefore, it follows by algebraic manipulation that $\Gamma(\sigma)$ is nonempty and equal to
                \[
                    \Gamma(\sigma)
                        {}={}
                    \set{(\gamma, \tau) \in \R^2}[
                    \gamma
                        {}\in{}
                    (0, +\infty)
                    \text{ and }
                    \tau
                        {}\in{}
                    \bigl(
                        0, \nicefrac{1}{\gamma\nrm{L}^2}
                    \bigr\rbrack
                    ].
                \]
                \item \ref{it:cor:quadratic:CP:existence:mixed}:
                If $\min\{\DRSrhocom,\DRSrhomon\} < 0$ and
                $\max\{\DRSrhocom,\DRSrhomon\} \geq 0$, then
                either $\tfrac{1}{\gamma}\DRSrhocom + \tfrac{1}{\tau}\DRSrhomon \geq 0$, in which case by definition $\Gamma_1 > 0$, or $\tfrac{1}{\gamma}\DRSrhocom + \tfrac{1}{\tau}\DRSrhomon < 0$, in which case $\Gamma_1 > \Gamma_2(\sigma_d)$ since $1-\gamma\tau\sigma_d^2 \geq 0$. Therefore, regardless of the sign of $\tfrac{1}{\gamma}\DRSrhocom + \tfrac{1}{\tau}\DRSrhomon$ it holds that
                \(
                    \Gamma(\sigma_d)
                \)
                reduces to
                \(
                    \set{(\gamma, \tau) \in \R^2_{++}}[
                        \gamma\tau \in \bigl(0, \nicefrac{1}{\nrm{L}^2}\bigr\rbrack \text{ and } \Gamma_2(\sigma_d) > 0
                    ].
                \)
                \begin{proofitemize}
                    \item $\DRSrhocom \geq 0$, $\DRSrhomon < 0$: Then, by construction $c_1(\sigma_d, \gamma) > 0$ and $c_2(\gamma) < 0$ for all $\gamma \in \R_+$. 
                    Consequently, 
                    $
                        \Gamma_2(\sigma_d) > 0
                            \Longleftrightarrow
                        c_1(\sigma_d, \gamma)\tau + c_2(\gamma) > 0
                            \Longleftrightarrow
                        \tau > \nicefrac{-c_2(\gamma)}{c_1(\sigma_d, \gamma)} = \tau_{\rm min}(\gamma)
                    $
                    and the set $\Gamma(\sigma_d)$ is nonempty if for some $\gamma > 0$ it holds that
                    \begin{align*}
                        \tau_{\rm min}(\gamma)
                            {}<{}
                        \nicefrac{1}{\gamma\nrm{L}^2} 
                            \;\Longleftrightarrow\;
                        c_2(\gamma)\gamma\nrm{L}^2 + c_1(\sigma_d, \gamma)
                            {}={}
                        \DRSrhomon\nrm{L}^2\gamma^2 + 
                        \Bigl(1 + \DRSrhocom\DRSrhomon\bigl(\nrm{L}^2 - \sigma_d^2\bigr)\Bigr)
                        \gamma + \DRSrhocom
                        > 0,
                        \numberthis\label{eq:cor:quadratic:CP:set:plusmin}
                    \end{align*}
                    which is guaranteed by \cref{it:lem:quadratic:existence}. 
                    Therefore, it follows from \Cref{it:lem:quadratic:solutions} that
                    \[
                        \Gamma(\sigma_d)
                            {}={}
                        \left\{(\gamma, \tau) \in \R^2 \; \middle| \;
                        \gamma
                            {}\in{}
                        \Bigl(0, \gamma_{\rm max}\Bigr)
                        \text{ and }
                        \tau
                            {}\in{}
                        \Bigl(
                            \tau_{\rm min}(\gamma), \nicefrac{1}{\gamma\nrm{L}^2}
                        \Bigr]
                        \right\}.
                    \]
                    \item $\DRSrhocom < 0$, $\DRSrhomon \geq 0$:
                    Observe that for all 
                    $
                        \gamma 
                            \leq
                        \nicefrac{-\DRSrhocom}{(1 - \DRSrhocom\DRSrhomon\sigma_d^2)}
                            <
                        -\DRSrhocom   
                    $
                    it holds that $c_1(\sigma_d, \gamma) \leq 0$ and $c_2(\gamma) < 0$.
                    Therefore, if $\Gamma_2(\sigma_d) > 0$ then necessarily 
                    $
                        \gamma > \nicefrac{-\DRSrhocom}{(1 - \DRSrhocom\DRSrhomon\sigma_d^2)}
                    $, implying 
                    $
                        c_1(\sigma_d, \gamma) > 0
                    $, so that     
                    \(
                        \Gamma(\sigma_d)
                    \)
                    reduces to
                    \(
                        \left\{(\gamma, \tau) \in \R^2_{++} \; \middle| \; \begin{array}{l}
                            \gamma\tau \in \bigl(0, \nicefrac{1}{\nrm{L}^2}\bigr],
                            \gamma > \nicefrac{-\DRSrhocom}{(1 - \DRSrhocom\DRSrhomon\sigma_d^2)}
                            \text{ and }
                            \tau > \tau_{\rm min}(\gamma)
                        \end{array}
                            \right\}.
                    \)
                    Ensuring that $\tau_{\rm min}(\gamma) <  \nicefrac{1}{\gamma\nrm{L}^2}$, i.e., solving \eqref{eq:cor:quadratic:CP:set:plusmin} as before, yields $                    
                        \gamma
                            {}\in{}
                        \left(\gamma_{\rm min}, +\infty\right)
                    $.
                    Finally, by observing that $\gamma_{\rm min} > -\DRSrhocom > \nicefrac{-\DRSrhocom}{(1 - \DRSrhocom\DRSrhomon\sigma_d^2)}$, it follows that
                    \[
                        \Gamma(\sigma_d)
                            {}={}
                        \left\{(\gamma, \tau) \in \R^2 \; \middle| \; 
                        \gamma
                            {}\in{}
                        \Bigl(\gamma_{\rm min}, +\infty\Bigr)
                        \text{ and }
                        \tau
                            {}\in{}
                        \Bigl(
                            \tau_{\rm min}(\gamma), \nicefrac{1}{\gamma\nrm{L}^2}
                        \Bigr]
                        \right\}.
                    \]
                \end{proofitemize}
                \item \ref{it:cor:quadratic:CP:existence:neg}:
                If $\max\{\DRSrhocom,\DRSrhomon\} < 0$, then 
                $\Gamma_1 > 0$ for $(\gamma, \tau) \in \R_{++}^2$ if and only if 
                $
                    \gamma > \nicefrac{-\DRSrhocom}2
                $
                and
                $
                    \tau > \nicefrac{-\gamma\DRSrhomon}{(2\gamma+\DRSrhocom)}
                $.
                As a result, the set $\Gamma(\nrm{L})$ is empty when $\DRSrhocom\DRSrhomon \geq \nicefrac{1}{\nrm{L}^2}$, as in this case the inequality
                $\nicefrac{-\gamma\DRSrhomon}{(2\gamma+\DRSrhocom)} < \nicefrac{1}{\gamma \nrm{L}^2}$ does not have a positive solution for $\gamma > \nicefrac{-\DRSrhocom}2$.
                Consider the following cases, assuming that $\DRSrhocom\DRSrhomon < \nicefrac{1}{\nrm{L}^2}$.
                \begin{proofitemize}
                    \item 
                    $\gamma \in
                    \left(
                        \nicefrac{-\DRSrhocom}{2},
                        \nicefrac{-\DRSrhocom}{(1 - \DRSrhocom\DRSrhomon\nrm{L}^2)}
                    \right)
                    $:
                    Then, $c_1(\nrm{L}, \gamma) < 0$ and thus $\Gamma_2(\nrm{L}) > 0$ if and only if
                    \(
                        \tau
                            {}<{}
                        \tau_{\rm min}(\gamma)
                    \).
                    Since it is easy to verify that in this case 
                    $
                        \nicefrac{-\gamma\DRSrhomon}{(2\gamma+\DRSrhocom)} 
                            > 
                        \tau_{\rm min}(\gamma)
                    $,
                    no such $\gamma$ belong to the set $\Gamma(\nrm{L})$.
                    \item $\gamma = \nicefrac{-\DRSrhocom}{(1 - \DRSrhocom\DRSrhomon\nrm{L}^2)}$: Then, $c_1(\nrm{L}, \gamma) = 0$ and $c_2(\gamma) < 0$, so that $\Gamma_2(\nrm{L}) < 0$ for all $\tau \in \R$.
                    \item $\gamma > \nicefrac{-\DRSrhocom}{(1 - \DRSrhocom\DRSrhomon\nrm{L}^2)}$: Then, $c_1(\nrm{L}, \gamma) > 0$ and thus $\Gamma_2(\nrm{L}) > 0$ if and only if
                    $
                        \tau
                            > 
                        \tau_{\rm min}(\gamma)
                    $,
                    which is by construction larger than
                    $
                        \nicefrac{-\gamma\DRSrhomon}{(2\gamma + \DRSrhocom)}.
                    $
                    The set $\Gamma(\sigma_d)$ is nonempty if for some $\gamma > 0$ it holds that
                    $
                        \tau_{\rm min}(\gamma)
                            <
                        \nicefrac{1}{\gamma\nrm{L}^2}
                    $,
                    which holds by \Cref{it:lem:quadratic:existence} if and only if  
                    $            
                        [\DRSrhocom]_-[\DRSrhomon]_- < \nicefrac{1}{4\nrm{L}^2}
                    $.
                    Finally, noting that $\gamma_{\rm min} > \nicefrac{-\DRSrhocom}{2}$,
                    it follows from \Cref{it:lem:quadratic:solutions} that 
                    \[
                        \Gamma(\nrm{L})
                            {}={}
                        \left\{(\gamma, \tau) \in \R^2 \; \middle| \; 
                        \gamma
                            {}\in{}
                        \Bigl(\gamma_{\rm min}, \gamma_{\rm max}\Bigr)
                        \text{ and }
                        \tau
                            {}\in{}
                        \Bigl(
                            \tau_{\rm min}(\gamma), \nicefrac{1}{\gamma\nrm{L}^2}
                        \Bigr]
                        \right\}.
                    \]
                \end{proofitemize}
                Since these sets corresponds to the stepsize ranges provided by \Cref{tab:CP:stepsize:rule}, the proof is completed.\qedhere
            \end{proofitemize}
        \end{proof}
\end{lemma}

\begin{lemma}\label{lem:LMI:solution:cartesian}
    Let $\D \dfn \begin{bsmallmatrix}\I_n & \I_n\end{bsmallmatrix}$, $(Y_1, Y_2) \in \parsumset$ and define $Y = Y_1 \oplus Y_2$.  
    Then, \eqref{eq:LMI:condition} holds and $X^\star$ as defined in 
    \eqref{eq:LMI:solution:best:matrix}
    is equal to $Y_1 \Box Y_2$.
    \begin{proof}
        Let
        \(
            E
                {}={}
            \tfrac12
            \begin{bsmallmatrix}
                1 \\
                -1
            \end{bsmallmatrix}
        \).
        Observe that
        \(
            \D^\dagger 
                {}={}
            \tfrac12 \begin{bsmallmatrix}
                \I_n \\
                \I_n
            \end{bsmallmatrix}
        \)
        and
        \(
            \proj_{\ker{\D}} 
                {}={} 
            \tfrac12 
            \begin{bsmallmatrix}
                \I_n & - \I_n \\ - \I_n & \I_n
            \end{bsmallmatrix}
        \),
        so that
        \ifjota
\else
            \begin{equation*}
                \begin{aligned}
                    \proj_{\ker{\D}} Y
                        {}={}&
                    (E \otimes \I_n)
                    \begin{bmatrix}
                        Y_1 & -Y_2
                    \end{bmatrix},
                    &
                    \proj_{\ker{\D}} Y \proj_{\ker{\D}}
                        {}={}&
                    (E E^\top)
                    \otimes
                    (Y_1 + Y_2),
                    \\
                    \rank(\proj_{\ker{\D}} Y)
                        {}={}&
                    \rank\begin{bmatrix}Y_1 & Y_2\end{bmatrix},
                    &
                    \rank(\proj_{\ker{\D}} Y \proj_{\ker{\D}})
                        {}={}&
                    \rank(Y_1 + Y_2).
                \end{aligned}
            \end{equation*}
        \fi
        Consequently, \eqref{eq:LMI:condition} holds owing to \cite[Thm. 9.2.4]{mitra2010Matrix}, since \(Y_1 + Y_2 \succeq 0\) and \(Y_1\) and \(Y_2\) are parallel summable. The claim for $X^\star$ follows from \Cref{it:DXD:optimal}, since
        \begin{equation*}\label{eq:lem:sum:solution}
            \begin{aligned}            
                X^\star
                    {}={}&
                \tfrac14 (Y_1 + Y_2) - \tfrac1{4}(Y_1 - Y_2)
                (E \otimes \I_n)^\top
                \left(
                (E E^\top)
                \otimes (Y_1 + Y_2)
                \right)^\dagger
                (E \otimes \I_n)
                (Y_1 - Y_2)\\
                    {}={}&
                \tfrac14 (Y_1 + Y_2) - \tfrac14 (Y_1 - Y_2)(Y_1 + Y_2)^\dagger (Y_1 - Y_2)\\
                \ifjota\else
                    {}={}&
                \tfrac14 (Y_1 + Y_2) - \tfrac14 (Y_1 + Y_2 - 2Y_2)(Y_1 + Y_2)^\dagger (Y_1 + Y_2 - 2Y_2)\\
                \fi
                    {}={}&
                \tfrac12 Y_2(Y_1 + Y_2)^\dagger (Y_1 + Y_2)
                + \tfrac12 (Y_1 + Y_2)(Y_1 + Y_2)^\dagger Y_2
                - Y_2 (Y_1 + Y_2)^\dagger Y_2
                    {}={}
                Y_1 \Box Y_2,
            \end{aligned}
        \end{equation*}
        where 
        the second equality holds since for arbitrary matrices $Z_1$, $Z_2$ it holds that $(Z_1 \otimes Z_2)^\dagger = Z_1^\dagger \otimes Z_2^\dagger$ and 
        the final equality holds 
        by definition of the parallel sum and parallel summability.
    \end{proof}
\end{lemma}

\begin{lemma}\label{lem:semi:bounds}
    Suppose that \Cref{ass:SWMVICPsemimonotone} holds and that the sets $\gph(\proj_{\rge{L^\top}} A^{-1})$ and $\gph(\proj_{\rge{L}} B)$ are not singletons. Then, it holds that $[\monA]_+[\comA]_+ \leq \nicefrac{1}{4\sigma_d^2}$ and $[\monB]_+[\comB]_+ \leq \nicefrac{1}{4\sigma_d^2}$.
    \begin{proof}
        Suppose that $\monA, \comA > 0$ and $\monB, \comB > 0$, for otherwise the two claims hold trivially.
        Consider $y_D \in \dom\bigl(A^{-1}\circ (-L^\top)\bigr) \cap \dom\bigl(B^{-1}\bigr) = \dom \Td \neq \emptyset$ and let $y = -L^\top y_D$. By semimonotonicity of $A$ at $(x^\star, - L^\top y^\star)$  it holds for all
        $(x^\star, y^\star) \in \pazocal{S}^\star$,
        $x \in A^{-1}(y) = A^{-1} \circ (- L^\top)(y_D)$ that
        \begin{align*}
            \inner{x - x^\star, y + L^\top y^\star} 
                {}\geq{}
            \qindef{x - x^\star}{\monA L^\top L} + \comA\nrm{y + L^\top y^\star}^2
                {}\geq{}
            \monA \sigma_d^2 \nrm{\proj_{\rge{L^\top}} (x - x^\star)}^2 + \comA\nrm{y + L^\top y^\star}^2.\numberthis\label{eq:cor:CP:semi:gammamax:semimon}
        \end{align*}
        Noting that by definition \((y , \proj_{\rge{L^\top}}x) \in \gph (\proj_{\range{L^\top}} A^{-1})\) and since $\gph(\proj_{\range{L^\top}} A^{-1})$ is not equal to the singleton $\set{(-L^\top y^\star, \proj_{\range{L^\top}} x^\star)}$, both of the above involved norms are nonzero. 

        On the other hand, since $y = -L^\top y_D \in \range{L^\top}$, it holds by the Fenchel-Young inequality with modulus  $2\monA\sigma_d^2 > 0$ that
        \begin{align*}
            \inner{x - x^\star, y + L^\top y^\star} 
                {}={}
            \inner{\proj_{\range{L^\top}}(x - x^\star), y + L^\top y^\star}
                {}\leq{}
            \monA\sigma_d^2\nrm{\proj_{\range{L^\top}}(x - x^\star)}^2 + 
            \tfrac{1}{4\monA\sigma_d^2}\nrm{y + L^\top y^\star}^2.
            \numberthis\label{eq:cor:CP:semi:gammamax:young}
        \end{align*}
        Combining \eqref{eq:cor:CP:semi:gammamax:semimon} and \eqref{eq:cor:CP:semi:gammamax:young}, it follows that $\comA \leq \tfrac1{4\monA\sigma_d^2}$.

        Analogously, consider $x_P \in \dom(\A) \cap \dom(\B \circ L) = \dom \Tp \neq \emptyset$ and let $x = Lx_P$. Then, it holds for all
        $(x^\star, y^\star) \in \pazocal{S}^\star$,
        $y \in B(x) = B(Lx_P)$ by the semimonotonicity assumption of $B$ at $(Lx^\star, y^\star)$ that
        \begin{align*}
            \inner{x - Lx^\star, y - y^\star} 
                {}\geq{}
            \monB\nrm{x - L x^\star}^2 + \qindef{y - y^\star}{\comB L L^\top}
                {}\geq{}
            \monB\nrm{x - L x^\star}^2 + \comB \sigma_d^2 \nrm{\proj_{\range{L}}(y - y^\star)}^2,\numberthis\label{eq:cor:CP:semi:gammamin:semimon}
        \end{align*} 
        where the involved norms are nonzero since \((x, \proj_{\rge{L}}y)\in \gph(\proj_{\range{L}} B)\) and  $\gph(\proj_{\range{L}} B)$ is not equal to the singleton $\set{L x^\star, \proj_{\range{L}} y^\star}$.
        On the other hand, since $x = Lx_P \in \range{L}$, it holds by the Fenchel-Young inequality with modulus $2\monB > 0$ that
        \begin{align*}
            \inner{x - L x^\star, y - y^\star} 
                {}={}
            \inner{x - L x^\star, \proj_{\range{L}}(y - y^\star)}
                {}\leq{}
            \monB\nrm{x - L x^\star}^2 + 
            \tfrac{1}{4\monB}\nrm{\proj_{\range{L}}(y - y^\star)}^2.
            \numberthis\label{eq:cor:CP:semi:gammamin:young}
        \end{align*}
        Finally, combining \eqref{eq:cor:CP:semi:gammamin:semimon} and \eqref{eq:cor:CP:semi:gammamin:young}, it follows that $\comB \sigma_d^2 \leq \tfrac1{4\monB}$, establishing the claim.
    \end{proof}
\end{lemma}

\begin{proposition}[normal cone of a box]\label{prop:box}
    The normal cone operator $N_C: \R^n \rightrightarrows \R^n$ of the $n$\hyp{}dimensional box 
    \(
        C \dfn \set{x\in \R^n}[l_i \leq x_i \leq u_i, \; i = 1, \dots, n]
    \) is 
    $
        \Bigl(\diag\left(\frac{|\tilde v_1|}{u_1 - l_1}, \dots, \frac{|\tilde v_n|}{u_n - l_n}\right), 0\Bigr)
    $-semimonotone at $(\tilde x, \tilde v) \in \gph N_C$.
    \begin{proof}
        By \Cref{lem:calculus:cartesian} it suffices to show that $N_{C_i}$ is $\Bigl(\frac{|\tilde v_i|}{u_i - l_i}, 0\Bigr)$-semimonotone at $(\tilde x_i, \tilde v_i) \in \gph N_{C_i}$. Using the fact that \(|\tilde x_i - x_i| \leq u_i-\ell_i\) and monotonicity of \(N_C\), we have for all $x_i \in C_i$ that
        \(
            \tfrac{|\tilde v_i|}{u_i-\ell_i}|\tilde x_i - x_i|^2 \leq  |\tilde v_i||\tilde x_i - x_i| = \inner{\tilde v_i, \tilde x_i - x_i}.
        \)
    \end{proof}
\end{proposition}

\section{Omitted proofs}\label{sec:omitted}

\begin{appendixproof}{lem:SWMI:primal-dual}
    Note that
    $
        \gph \Tpd
    $
    is equal to the set
    $
        \set{\bigl((x_\A, y_\B),
        (y_\A + L^\top y_\B, 
        x_\B - L x_\A)\bigr)}[(x_\A, y_\A) \in \gph \A, (x_\B, y_\B) \in \gph B].
    $
    Consequently, it holds by assumption for all 
    $(x^\star, y^\star) \in \pazocal{S}^\star$,
    $(x_\A,y_\A) \in \gph \A$ and
    $(x_\B, y_\B) \in \gph B$ that
    \begin{align*}
        \inner*{y_\A + L^\top y_\B, x_\A - x^\star} 
        + \inner*{x_\B - L x_\A, y_\B - y^\star} \geq 
        \qindef{y_\A + L^\top y_\B}{\DRSRho_{\rm P}}
        + \qindef{x_\B - L x_\A}{\DRSRho_{\rm D}},
        \numberthis\label{eq:WMI:primal-dual}
    \end{align*}
    \vspace{-0.4cm}
    \begin{proofitemize}
        \item In \eqref{eq:WMI:primal-dual}, consider $x_\A \in \dom(A) \cap \dom(B \circ L) = \dom \Tp \neq \emptyset$ and let $x_\B = Lx_\A$. Then, it holds for all
        $x^\star \in \pazocal{S}_{\rm P}^\star$,
        $y_\A \in A(x_\A)$
        and
        $y_\B \in B(Lx_\A)$ that 
        \ifjota
\else
            \begin{align*}
                \inner*{y_\A + L^\top y_\B, x_\A - x^\star} 
                    \geq 
                \qindef{y_\A + L^\top y_\B}{\DRSRho_{\rm P}}.
                \numberthis\label{eq:SWMI:primal}
            \end{align*}
        \fi
        Since $(x_\A, y_\A + L^\top y_\B) \in \gph \Tp$ by construction
        and $\pazocal{S}_{\rm P}^\star \subseteq \zer \Tp$ by \Cref{prop:auslender},
        it follows by definition that
        $\Tp$ has $V_{\rm P}$\hyp{}oblique weak Minty solutions at $\pazocal{S}^\star_{\rm P}$.
        \item
        Analogously, consider $y_\B \in \dom\bigl(A^{-1}\circ (-L^\top)\bigr) \cap \dom\bigl(B^{-1}\bigr) = \dom \Td \neq \emptyset$ and let $y_\A = -L^\top y_\B$ in \eqref{eq:WMI:primal-dual}. Then, it holds for all
        $y^\star \in \pazocal{S}_{\rm D}^\star$,
        $x_\A \in A^{-1} \circ (- L^\top)(y_\B)$
        and
        $x_\B \in B^{-1}(y_\B)$ that
        \ifjota
\else
            \begin{align*}
                \inner*{x_\B - L x_\A, y_\B - y^\star}
                    \geq 
                \qindef{x_\B - L x_\A}{\DRSRho_{\rm D}}.
                \numberthis\label{eq:SWMI:dual}
            \end{align*}
        \fi
        Since $(y_\B, x_\B - L x_\A) \in \gph \Td$ by construction
        and $\pazocal{S}_{\rm D}^\star \subseteq \zer \Td$ by \Cref{prop:auslender},
        it follows by definition that $\Td$ has $V_{\rm D}$\hyp{}oblique weak Minty solutions at $\pazocal{S}^\star_{\rm D}$, completing the proof \qedhere
    \end{proofitemize}
\end{appendixproof}

\begin{appendixproof}{ex:CP:example5}[ (saddle point problem)]
    \begin{proofitemize}
        \item \ref{it:CP:example5:1}: By defining
        \(H \dfn \I + \lambda\bigl((\M + \Tpd)^{-1}\M - I\bigr)\) and substituting $\tau = \nicefrac{1}{\gamma \ell^2}$, the update rule for $z^k$ corresponds to the linear dynamical system $z^{k+1} = H z^k$. 
        Global asymptotic stability of this system is achieved if and only if the spectral radius of the matrix \(H\) is strictly less than one, which holds iff $\lambda \in \bigl(0, \min\bigl\{2, \bar{\lambda}\bigr\}\bigr)$. Analogously, the convergence result for $\seq{\projQ z^k}$ can be obtained by analyzing the spectral radius of $\projQ H$.
        \item \ref{it:CP:example5:2}:
        The primal-dual operator and its inverse are given by 
        \begin{align*}
            \Tpd 
                {}={}
            \begin{bmatrix}
                A & \tp L \\ -L & B^{-1}
            \end{bmatrix}
            \quad\text{and}\quad
            \Tpd^{-1}
                {}={}
            \begin{bmatrix}
                (A + \tp LBL)^{-1} & - (A + \tp LBL)^{-1}\tp L B\\
                B L (A + \tp LBL)^{-1} & B - B L (A + \tp LBL)^{-1}\tp L B
            \end{bmatrix}
            \numberthis\label{eq:example:linear:Tpd}
        \end{align*}
        owing to the Schur complement lemma.
        Therefore, $\Tpd$ has a 
        \(
            \DRSRho
        \)\hyp{}oblique weak Minty solution at zero on
        \(\pazocal{U} = \R^5 \times \range{\M}\)
        if and only if
        \begin{equation}
            \tp z \left(\tfrac{\Tpd + \tp \Tpd}{2} - \tp \Tpd \DRSRho \Tpd\right) z \geq 0,\qquad \text{for all $z \in \R^n : z \in \Tpd^{-1}\range{\M}$}.\label{eq:example:linear:WMVI}
        \end{equation}
        Using that 
        $
            L
                {}={}
            \UL \SL \VL^\top
        $,
        where $\UL = \I_3$, 
        $
            \SL 
                = 
            \begin{bsmallmatrix}
                |\ell| \I_2\\
                0
            \end{bsmallmatrix}
        $
        and $\VL = \sign(\ell)\I_2$, it follows by plugging in \eqref{eq:proof:CP:indef:basisP1} into \eqref{eq:proof:CP:basisP} that
        \begin{align*}
            \basisP = 
            \begin{bsmallmatrix} 
                \tfrac{1}{\sqrt{1+\gamma^2\ell^2}}\I_2 & 0 \\
                -\tfrac{\gamma \ell}{\sqrt{1+\gamma^2\ell^2}}\I_2 & 0\\
                0 & 1
            \end{bsmallmatrix}
            \quad\text{and}\quad
            \DRSRho
                {}={}
            \tfrac{b\ell^2}{a^2 + b^2\ell^4} \I_2 \oplus \tfrac{b a^2}{a^2 + b^2\ell^4} \I_2 \oplus c \I_1,
            \numberthis\label{eq:example:linear:UV}
        \end{align*}
        where $\basisP$ is an orthonormal basis for $\range{\M}$. As a result, \eqref{eq:example:linear:WMVI} is satisfied if and only if
        \begin{align*}
            \bigl(\Tpd^{-1} \basisP\bigr)^\top\left(\tfrac{\Tpd + \tp \Tpd}{2} - \tp \Tpd \DRSRho \Tpd\right) \Tpd^{-1} \basisP
            \succeq 0.
        \end{align*}
        Noting that the left-hand side becomes zero when plugging in \eqref{eq:example:linear:Tpd} and \eqref{eq:example:linear:UV}, the claim follows. 
        
        \item
        \ref{it:CP:example5:3}:
        Since condition \eqref{eq:existence:condition} holds for
        $\DRSrhocom = \tfrac{b\ell^2}{a^2 + b^2\ell^4}$,
        $\DRSrhomon = \tfrac{b a^2}{a^2 + b^2\ell^4}$
        and
        $\DRSrhomon^\prime = c$, it follows from \Cref{it:CP:example5:2} that \Cref{ass:SWMVICP} is satisfied.
        Moreover, it can be verified that for these parameters $\delta$ is equal to one, that \Cref{ass:CP:stepsize:rule} reduces to
        \eqref{eq:example5:gamma}
        and that 
        \Cref{ass:CP:relaxation:rule}
        matches \eqref{eq:example5:lambda}.
        Operators $A$ and $B$ are continuous and hence outer semicontinuous, verifying \Cref{ass:CP:1}.
        It only remains to show that $J_{\gamma A}$ and $J_{\tau B^{-1}}$ have full domain and are single-valued continuous.
        For $J_{\gamma A}$ this holds since $\I + \gamma A$ is invertible for any $\gamma \in \R$. For $J_{\tau B^{-1}}$ with $\tau = \nicefrac{1}{\gamma\ell^2}$, invertability of $\I + \tau B^{-1}$ is ensured provided that $-\gamma b\ell^2 \neq 1$ and $-\gamma c\ell^2 \neq 1$. Since these conditions are guaranteed by \Cref{ass:CP:stepsize:rule} and the assumption $c \geq 0$, the claim is established.

        \item \ref{it:CP:example5:4}: 
        For this particular instance 
        \(
            \trace\Tp = -2
        \)
        and 
        \(
            \trace\Td = \trace\Tpd = -8
        \).
        Since the trace of a matrix equals the sum of its eigenvalues, the proof is completed.
        \qedhere
    \end{proofitemize}
\end{appendixproof}

\begin{appendixproof}{ex:CP:second}[ (influence of singular values)]
    \begin{proofitemize}
        \item \ref{ex:CP:second:Tpd}: 
        By \cite[Prop. 5.1(ii)]{bauschke2021Generalized} and using that $A$, $B$ and $L$ are symmetric, $\Tpd$ is $\tfrac12$\hyp{}comonotone if and only if 
        \begin{align*}
            \tfrac{\Tpd + \Tpd^\top}{2} - \tfrac12 \Tpd^\top \Tpd
                {}={}&
            \begin{bmatrix}
                A - \tfrac12 (A^\top A + L^\top L)& -\tfrac12 (AL - LB^{-1})\\
                -\tfrac12 (LA - B^{-1}L) & B^{-1} - \tfrac12 ({B^{-1}}^\top B^{-1} + LL^\top)
            \end{bmatrix}
                \succeq
            0.\numberthis\label{eq:ex:second:comonotone}
        \end{align*}
        Using that $B^{-1} = A$, that $L$ is symmetric and that $A$ and $L$ commute, i.e., $A L = L A$, this condition reduces to $A - \tfrac12 (A^\top A + L^\top L) \succeq 0$, which holds by definition of $A$ and $L$. Noting that $\zer \Tpd = (0_n, 0_n)$, the claim is established.
        \item \ref{ex:CP:second:theorem}:
        Follows from \Cref{thm:CP:full} and \Cref{rem:PDHG:sequences}, using that $\nrm{L} = 1$ and $\gamma \tau = \tfrac{1}{\nrm{L}^2}$. 
        \item \ref{ex:CP:second:spectral}: 
        Analogous to the setting of \Cref{it:CP:example5:1}, the update rule for $\projQ z^k$ can be expressed as the linear dynamical system $\projQ z^{k+1} = \projQ H \projQ z^k$, where \(H \dfn \I + \lambda\left((\M + \Tpd)^{-1}\M - I\right)\) and \(\M\) is defined as in \Cref{eq:CP:P}.
        This system is globally asymptotically stable if and only if the spectral radius of \(\projQ H\) is strictly less than one, i.e., if and only if $\lambda \in (0, \bar \lambda_{\rm spectral})$, where
        \begin{align*}
            \bar \lambda_{\rm spectral}
                \in
            \argmax_\lambda 
            \;\,
            \lambda 
            \;\,
            \stt 
            \;\,
            \nrm{\projQ H(\lambda)}_2 < 1.
            \numberthis\label{eq::CP:second:spectral:lambda}
        \end{align*}
        The values for $\bar \lambda_{\rm spectral}$ reported in \Cref{fig:ex:second:bound} are obtained by solving this problem using SymPy. \qedhere
    \end{proofitemize}
\end{appendixproof}

\begin{appendixproof}{prop:semi:young}
        By the Fenchel-Young inequality, it holds for any
        $
            \Com
                \prec
            0
        $
        that
        \ifjota
\else
            \begin{align*}
                \inner{x-\other{x},y-\other{y}} \geq \tfrac{1}{4}\qindef{x - \other{x}}{\Com^{-1}} - \qindef{y - \other{y}}{\Com}, \qquad \text{for all }(x, y),(\other{x}, \other{y}) \in \gph A.
            \end{align*}
        \fi
        Therefore, \eqref{eq:obliquequasisemimonotonicity:matrix} is satisfied for all $\Mon \preceq \tfrac14 \Com^{-1}$.
\end{appendixproof}

\begin{appendixproof}{prop:semi:maximality}
    If \(T\) is maximally \((\Mon', \Com')\)\hyp{}semimonotone, then by definition
    there exists no \((\Mon', \Com')\)\hyp{}semimonotone operator $\tilde T$ such that \(\gph T \subset \gph \tilde T\).
    Since the class of \((\Mon', \Com')\)\hyp{}semimonotone operator encompasses the class of \((\Mon, \Com)\)\hyp{}semimonotone operators whenever $\Mon' \preceq \Mon$ and $\Com' \preceq \Com$,
    this implies that there also exists no \((\Mon, \Com)\)\hyp{}semimonotone operator $\tilde T$ such that \(\gph T \subset \gph \tilde T\), establishing the claim by definition of maximality.
\end{appendixproof}

\begin{appendixproof}{lem:calculus}
    \begin{proofitemize}
        \item \ref{lem:calculus:scaling}:
        First, consider the assertion where $A$ is semimonotone only at $(\other{x}_A,\other{y}_A)$. Define $\other{s} = \other{x}_A - u$ and $\other{t} = y + \alpha \other{y}_A$, such that $(\other{s},\other{t}) \in \gph{T}$. Then, it holds for all $(s,t) \in \gph T$ that
        \begin{align*}
            \inner*{s - \other{s}, t - \other{t}} 
                {}={}&
            \alpha \inner*{(s+u) - (\other{s}+u), \alpha^{-1}(t-y) - \alpha^{-1}(\other{t}-y)}\\
                \dueto{(semimonotonicity of $A$ at $(\other{x}_A,\other{y}_A)$)}{}\geq{}&
            \alpha \qindef{(s+u) - (\other{s}+u)}{\Mon} + \alpha \qindef{\alpha^{-1}(t-y) - \alpha^{-1}(\other{t}-y)}{\Com}\\
                {}={}&
            \qindef{s - \other{s}}{\alpha\Mon} + \qindef{t - \other{t}}{\alpha^{-1}\Com},
        \end{align*}
        where we used that $(\other{s}+u,\alpha^{-1}(\other{t}-y)) = (\other{x}_A, \other{y}_A)$ and $(s+u,\alpha^{-1}(t-y)) \in \gph{A}$. Hence, it follows that $T$ is $(\alpha\MonA,\alpha^{-1}\ComA)$\hyp{}semimonotone at $(\other{s}, \other{t})$.

        If $A$ is $(\MonA,\ComA)$\hyp{}semimonotone at all $(\other{x}_A,\other{y}_A) \in \gph A$, we then know that $T$ is $(\alpha\MonA,\alpha^{-1}\ComA)$\hyp{}semimonotone at all points in the set $\set*{(\other{x}_A-u,y + \alpha\other{y}_A)}[(\other{x}_A,\other{y}_A) \in \gph A]$. Since this set is equal to $\gph T$, it follows that $T$ is $(\alpha\MonA,\alpha^{-1}\ComA)$\hyp{}semimonotone (everywhere). 

        \item \ref{lem:calculus:cartesian}:
        Let $A$ and $B$ be semimonotone at respectively $(\other{x}_A, \other{y}_A)$ and $(\other{x}_B, \other{y}_B)$. Since $\gph T$ is equal to the set $\set{((x_\A,x_\B),(y_\A,y_\B))}[y_\A \in \A x_\A, y_\B \in \B x_\B]$, it holds for all $(x, y) \in \gph T$ that
        \begin{align*}
            \inner*{x - \other{x}, y - \other{y}}
                {}={}&
            \inner*{x_\A - \other{x}_\A, y_\A - \other{y}_\A} + \inner*{x_\B - \other{x}_\B, y_\B - \other{y}_\B}\\
                \dueto{(semimonotonicity of $A$ and $B$)}{}\geq{}&
            \qindef{x_\A - \other{x}_\A}{\MonA} + \qindef{y_\A - \other{y}_\A}{\ComA} +  \qindef{x_\B - \other{x}_\B}{\MonB} + \qindef{y_\B - \other{y}_\B}{\ComB}\\
                {}={}&
            \qindef{x - \other{x}}{\MonA \oplus \MonB} + \qindef{y - \other{y}}{\ComA \oplus \ComB},
        \end{align*}
        and thus $T$ is $\bigl(\MonA \oplus \MonB, \ComA \oplus \ComB\bigr)$\hyp{}semimonotone at $(\other{x}, \other{y}) \in \gph T$. 

        If $A$ and $B$ are semimonotone at all points in their graph, then $T$ is $\bigl(\MonA \oplus \MonB, \ComA \oplus \ComB\bigr)$\hyp{}semimonotone at all points in $\gph T$, which completes the proof.
        \qedhere
    \end{proofitemize}
\end{appendixproof}

\begin{appendixproof}{prop:semi:lin}
    Owing to the linearity of $D$, \((\Mon, \Com)\)\hyp{}semimonotonicity corresponds to having $\inner*{x, Dx} \geq \qindef{x}{\Mon} + \qindef{Dx}{\Com}$ for all $x \in \R^n$, which is equivalent to the claimed LMI.
\end{appendixproof}

\begin{appendixproof}{lem:LMI:solution}
    \begin{proofitemize}
        \item \ref{it:DXD:existence}: First, note that the problem of finding an \(X \in \sym{n}\) such that \(\tp \D X \D \preceq Y\) is equivalent to the problem of finding a pair \((X, Z) \in \sym{n} \times \sym{m}\) such that
        \ifjota
\else
            \begin{align*}
                \tp \D X \D = Y - Z\numberthis\label{eq:DXD:equivalent}
            \end{align*}
        \fi
        and \(Z \succeq 0\).
        Second, observe that by 
        \cite[Thm. 2]{penrose1955generalized}, \cite[Prop. 1]{grosz1998note}
        the involved linear matrix equality is solvable for $X \in \sym{n}$ if and only if
        \(
            \range{Y - Z} \subseteq \range{\tp \D}
        \),
        i.e.,
        \(
            \proj_{\ker{\D}}Z = \proj_{\ker{\D}} Y.
        \)
        By \cite[Thm. 2.2]{khatri1976hermitian}, a matrix $Z \succeq 0$ satisfying this condition exists 
        if and only if \eqref{eq:LMI:condition} holds, and the general solution is given by
        \begin{align*}
            Z = Y \proj_{\ker{\D}}(\proj_{\ker{\D}} Y \proj_{\ker{\D}})^\dagger \proj_{\ker{\D}} Y + \tp \D (\D^\dagger)^\top G \tp \D (\D^\dagger)^\top,\numberthis\label{eq:DXD:equivalent:ness+suff:solution}
        \end{align*}
        where \(G \in \sym{n}\) is an arbitrary symmetric positive semidefinite matrix.
        \item \ref{it:DXD:optimal}: 
        \ifjota
\else
            Substituting \eqref{eq:DXD:equivalent:ness+suff:solution} into \eqref{eq:DXD:equivalent} yields
            \begin{align*}
                \tp \D X \D = Y - Y \proj_{\ker{\D}}(\proj_{\ker{\D}} Y \proj_{\ker{\D}})^\dagger \proj_{\ker{\D}} Y - \tp \D (\D^\dagger)^\top G \tp \D (\D^\dagger)^\top,
            \end{align*}
            of which the general solution is given by \cite[Prop. 1]{grosz1998note}
        \fi
        \begin{align*}
            X
            {}={}
            \underbracket{
                (\D^\dagger)^\top (Y
                - Y \proj_{\ker{\D}}(\proj_{\ker{\D}} Y \proj_{\ker{\D}})^\dagger \proj_{\ker{\D}} Y - G)\D^\dagger
            }_{=\, X^\star - (\D^\dagger)^\top G \D^\dagger}
            + H - \proj_{\rge{\D}} H \proj_{\rge{\D}},
            \numberthis\label{eq:DXD:solution:general}
        \end{align*}
        where \(H \in \sym{n}\) is an arbitrary matrix.
        Substituting \eqref{eq:DXD:solution:general} into \(\tp \D X \D \preceq Y\) shows that
        \ifjota
\else
            \begin{align*}
                Y - \tp \D X \D = Y - \tp \D X^\star \D + \proj_{\range{\tp \D}} G \proj_{\range{\tp \D}} \succeq Y - \tp \D X^\star \D \succeq 0.
            \end{align*}
        \fi
        Finally, the alternative expression for $X^\star$ given in \eqref{eq:LMI:solution:best:matrix} follows directly from \cite[Lem. 3]{yonglin1990g}.
        \item \ref{it:DXD:consistent}: \cite[Fact 6.4.38]{bernstein2009matrix}
        \qedhere
    \end{proofitemize}
\end{appendixproof}

\begin{appendixproof}{lem:calculus:composition:1}
    If \eqref{eq:LMI:condition} holds, then it follows from \Cref{lem:LMI:solution} that \(\D^\top X^\star \D \preceq Y\). Therefore, it only remains to be to shown that this implies \((\D \Mon \D^\top,X^\star)\)\hyp{}semimonotonicity of \(\D \F \D^\top\) \optional{at \((\other{x}, \D \other{y})\)}.
    \ifjota\else

    \fi
    First, consider the case where \(\F\) is semimonotone only at a single point \((\D^\top\other{x},\other{y})\). Let \((\D^\top x,y) \in \gph \F\) and denote \(u = \D y\) and \(\other{u} = \D \other{y}\). Then, \((x,u),(\other{x},\other{u}) \in \gph \D \F \D^\top\) and it holds that
    \begin{align*}
        \inner{x-\other{x},u-\other{u}} 
            {}={}&
        \inner{\D^\top(x-\other{x}),y-\other{y}}\\ 
            \dueto{(semi. of \(\F\) at $(\D^\top\other{x}, \other{y})$)}
            {}\geq{}&
        \qindef{\D^\top(x-\other{x})}{\Mon} + \qindef{y-\other{y}}{Y}\\
            {}\geq{}&
        \qindef{x-\other{x}}{\D \Mon \D^\top} + \qindef{\D(y-\other{y})}{X^\star}
            {}={}
        \qindef{x-\other{x}}{\D \Mon \D^\top} + \qindef{u-\other{u}}{X^\star},
    \end{align*}
    where \(\D^\top X^\star \D \preceq Y\) was used in the second inequality, showing that \(\D \F \D^\top\) is \((\D \Mon \D^\top,X^\star)\)\hyp{}semimonotone at \((\other{x}, \D \other{y})\).
    Finally, if \(\F\) is \((\Mon,Y)\)\hyp{}semimonotone at all \((\D^\top\other{x},\other{y}) \in \gph \F\), then \(\D \F \D^\top\) is \((\D \Mon \D^\top,X^\star)\)\hyp{}semimonotone at all points in \(\set{(\other{x}, \D \other{y})}[(\D^\top\other{x},\other{y}) \in \gph \F]\), which equals \(\gph \D \F \D^\top\).
\end{appendixproof}

\begin{appendixproof}{prop:calculus:sum+parsum}
    Let $\D = \left[\begin{smallmatrix}\I_n & \I_n\end{smallmatrix}\right]$. Then, \(A + B\) is equal to \(\D T \D^\top\), where $
    T \dfn A \times B$.
    By \Cref{lem:calculus:cartesian}, operator
    \(T\)
    is \((\Mon, \Com) = \bigl(\MonA \oplus \MonB, \ComA \oplus \ComB\bigr)\)\hyp{}semimonotone \optional{at $((\other{x}, \other{x}), (\other{y}_\A, \other{y}_\B)) \in \gph T$}. 
    Consequently, it follows from \Cref{lem:calculus:composition:1,lem:LMI:solution:cartesian} that $\D T \D^\top = \A + \B$ is
    \((\MonA + \MonB, \ComA \Box \ComB)\)\hyp{}semimonotone \optional{at $(\other{x}, \other{y}_\A + \other{y}_\B)$}. Finally, the claim for the parallel sum follows directly from those for the sum and \Cref{prop:calculus:inverse}, since $A \Box B \dfn (A^{-1} + B^{-1})^{-1}$. 
\end{appendixproof}

\begin{appendixproof}{lem:calculus:sum:linear}
    First, consider the assertion where \(\F\) is semimonotone at \((\other{x}, \other{y})\). Let \((x,y) \in \gph{\F}\). Then, \((x, y + Dx), (\other{x}, \other{y} + D \other{x}) \in \gph{(T+D)}\) and
    \ifjota
\else
        \begin{align*}
            \inner{x - \other{x}, y - \other{y} + \D(x  - \other{x})}
                {}\geq{}&
            \qindef{x - \other{x}}{D^\top \Mon D}
            +
            \qindef{y - \other{y}}{\Com + \Com^\prime}
            +
            \inner{x - \other{x}, \D (x - \other{x})}.
        \end{align*}
    \fi
    By \Cref{prop:semi:lin:optimal} and skew-symmetry of $\D$, it follows that \(D\) is \((-D^\top \Mon D, \proj_{\range{D}}\Mon\proj_{\range{D}})\)\hyp{}semimonotone. Consequently, 
    \begin{align*}
    \inner{x - \other{x}, y - \other{y} + D (x - \other{x})}
    {}\geq{}&
    \qindef{y - \other{y}}{\Com + \Com^\prime}
    +
    \qindef{D x - D\other{x}}{\proj_{\range{D}}\Mon\proj_{\range{D}}}\\
    {}={}&
    \qindef{y - \other{y} + \D (x - \other{x})}{\Com^\prime}
    +
    \qindef{y - \other{y}}{\Com}
    +
    \qindef{D(x - \other{x})}{\Mon}\\
        {}={}&
    \qindef{y - \other{y} + \D (x - \other{x})}{\Com^\prime}
    +
    q_{\Mon \oplus \Com}
    \left(
        (D(x - \other{x}),
        y - \other{y})
    \right)
    \numberthis\label{eq:proof:lem:calculus:sum:linear:auxiliary}
    \end{align*}
    where the fact that \(\range{\Com'} \subseteq \ker{\D}\) was used in the first equality.
    Let $E \coloneqq \left[\begin{smallmatrix}\I_n & \I_n\end{smallmatrix}\right]$.
    Since \((\Mon, \Com) \in \parsumset\), it holds by \Cref{lem:LMI:solution,lem:LMI:solution:cartesian}
    that 
    \(
        E^\top (\Mon \Box \Com) E \preceq \Mon \oplus \Com
    \).
    Consequently, it follows from \eqref{eq:proof:lem:calculus:sum:linear:auxiliary} that
    \begin{align*}
        \inner{x - \other{x}, y - \other{y} + D (x - \other{x})}
            {}\geq{}&
        \qindef{y - \other{y} + \D(x  - \other{x})}{\Com^\prime + \Mon \Box \Com},
    \end{align*}
    proving the claimed semimonotonicity of $T + D$
    at
    $(\other{x}, \other{y} + D \other{x})$.
    Finally, when \(T\) is semimonotone everywhere,
    the claim follows analogously by considering all $(\other{x}, \other{y}) \in \gph T$, completing the proof.
\end{appendixproof}
\begin{appendixproof}{prop:constrainedQP}
    \begin{proofitemize}
        \item \ref{it:prop:constrainedQP:A}: 
        Let $\Mon \dfn L^\top \MonA L = \proj_{\rge{L^\top}} Q \proj_{\rge{L^\top}}$ and observe that
        \begin{align*}
            \proj_{\ker{Q}} \Mon
                {}={}
            \proj_{\ker{Q}} \proj_{\rge{L^\top}} Q \proj_{\rge{L^\top}}
                {}={}
            \proj_{\ker{Q}} (\I_n - \proj_{\ker{L}}) Q \proj_{\rge{L^\top}}
                {}={}
            0,
            \numberthis\label{eq:proof:ex:constrainedQP:projQ-M}
        \end{align*}
        where the final equality holds since $\proj_{\range{L^\top}} Q \proj_{\ker{L}} = 0$.
        Therefore, \eqref{eq:LMI:condition:lin} is satisfied for $D = Q$ and $Q$ is $(L^\top \MonA L, \Com^\star)$\hyp{}semimonotone by \Cref{prop:semi:lin:optimal}, where
        \(
            \Com^\star 
                {}={}
            Q^\dagger - Q^\dagger \Mon Q^\dagger
                {}={}
            Q^\dagger - Q^\dagger \proj_{\rge{L^\top}}Q \proj_{\rge{L^\top}} Q^\dagger
        \)
        due to \eqref{eq:proof:ex:constrainedQP:projQ-M} and symmetry of $Q$.
        Moreover, since
        $
            \proj_{\range{L^\top}} Q \proj_{\ker{L}} = 0,
        $
        it holds by \cite[Fact 6.4.34]{bernstein2009matrix} that
        $
            Q^\dagger = (\proj_{\range{L^\top}} Q \proj_{\range{L^\top}})^\dagger + (\proj_{\ker{L}} Q \proj_{\ker{L}})^\dagger
        $,
        so that 
        $
            \Com^\star
                {}={}
            \proj_{\ker{L}} Q^\dagger \proj_{\ker{L}}
                {}={}
            \ComA^\prime
        $.
        The claimed result for $A : x \mapsto Qx + q$ then follows from \Cref{lem:calculus:scaling}.
        \item \ref{it:prop:constrainedQP:B}: By \Cref{prop:auslender}, it holds that $(x^\star, y^\star) \in \zer \Tpd$ if and only if $(x^\star, -L^\top y^\star) \in \gph A$ and $(L x^\star, y^\star) \in \gph B$. The claimed result then follows directly from \Cref{prop:box}.
        \item \ref{it:prop:constrainedQP:Tpd}: 
        By \Cref{thm:calculus:primaldual}, $\Tpd$
        has
        $((\proj_{\ker{L}} Q^\dagger \proj_{\ker{L}}) \oplus (\MonA \Box \MonB))$\hyp{}oblique weak Minty solutions at $\pazocal{S}^\star = \set{(x^\star, y^\star)}$.
        Therefore, recalling that $\proj_{\ker{L}} = \VL^\prime {\VL^\prime}^\top$, it follows that the parameters $\DRSrhocom, \DRSrhomon, \DRSrhocom'$ and $\DRSrhomon'$ given in \ref{it:prop:constrainedQP:Tpd} match those from \eqref{eq:calculus:primaldual:params}.
        As remarked below \Cref{thm:calculus:primaldual}, \Cref{ass:SWMVICP} holds since these parameters satisfy condition \eqref{eq:existence:condition}, completing the proof.
        \item \ref{it:prop:constrainedQP:convergence}:
        Note that both $A$ and $B$ are outer semicontinuous since $A$ is a linear mapping and $B$ is a normal cone \cite[Prop. 6.6]{rockafellar2009Variational}, verifying \Cref{ass:CP:1}.
        Since the parameters 
        \(
            \DRSrhocom,\DRSrhomon,\DRSrhocom^\prime
        \)
        and
        \(
            \DRSrhomon^\prime
        \)
        from \Cref{it:QP:Tpd:full:Semi}
        satisfy condition \eqref{eq:existence:condition}, it follows that also \Cref{ass:SWMVICP} holds.
        Moreover, it is straightforward to verify that for these parameter values \Cref{ass:CP:stepsize:rule,ass:CP:relaxation:rule} correspond to \eqref{eq:QP:convergence:stepsizes} (excluding the case $\tau = \nicefrac{1}{\gamma\nrm{L}^2}$ for simplicity).
        Therefore, it only remains to show that the resolvents $J_{\gamma A}$ and $J_{\tau^{-1} B}$ have full domain and are (single-valued) continuous for any $\gamma$ and $\tau$ complying with \eqref{eq:QP:convergence:stepsizes}.
        Note that for any $\tau > 0$ it holds that
        $
            J_{\tau B^{-1}}
                =
            \id - \tau J_{\tau^{-1} B} \circ (\tau^{-1}\id).
        $
        Since $J_{\tau^{-1} B} = \proj_C$ is (single-valued) continuous and has full domain, this implies that so does $J_{\tau B^{-1}}$.
        On the other hand, since $A(x) = Qx + q$ it holds that
        \(
            J_{\gamma A}(\cdot) = J_{\gamma Q}(\cdot - \gamma q)
        \).
        Therefore, $J_{\gamma A}$ has full domain and is (single-valued) continuous if and only if so is $J_{\gamma Q}$, i.e., if and only if the matrix $\id + \gamma Q$ is invertible.
        We distinguish between two cases.
        \begin{itemize}
            \item $\lambda_{\rm min}(Q) \geq 0$: Then, $\id + \gamma Q$ is invertible for any $\gamma > 0$ and the claim is immediate. 
            \item $\lambda_{\rm min}(Q) < 0$:
            Let $\VL, \VL', \UL, \UL'$ be defined as in \eqref{eq:L:svd}, so that $L = \UL \SL \VL^\top$.
            We begin by showing the following intermediate claims.
            \begin{claims}
                \item\label{claim:QP:1}
                    \textit{
                        The smallest eigenvalue of $Q$ is equal to the smallest eigenvalue of 
                        \(
                            \VL^\top Q \VL
                        \).
                    }\newline
                    Since by assumption 
                    $
                        \proj_{\range{L^\top}} Q \proj_{\ker{L}} = 0
                    $
                    where 
                    \(
                        \proj_{\range{L^\top}} = \VL \VL^\top 
                    \)
                    and 
                    \(
                        \proj_{\ker{L}} = \VL^\prime {\VL^\prime}^\top
                    \),
                    we have
                    \[
                        Q =
                        \begin{bmatrix}
                            \VL & \VL^\prime
                        \end{bmatrix}
                        \underbrace{
                        \begin{bmatrix}
                            \VL^\top Q \VL & 0\\
                            0 & {\VL^\prime}^\top Q \VL^\prime
                        \end{bmatrix}}_{\eqqcolon \bar Q}
                        \begin{bmatrix}
                            \VL^\top \\ {\VL^\prime}^\top
                        \end{bmatrix}
                    \]
                    Consequently, it follows by the similarity transformation with similarity matrix
                    \(
                        \begin{bmatrix}
                            \VL & \VL^\prime
                        \end{bmatrix}
                    \)
                    \cite[Cor. 1.3.4]{horn2012matrix}
                    that
                    \(Q\) is similar to \(\bar Q\), and therefore
                    \(
                        \lambda_{\rm min}(Q)
                            {}={}
                        \min\{\lambda_{\rm min}(\VL^\top Q \VL), \lambda_{\rm min}({\VL^\prime}^\top Q \VL^\prime)\}
                    \),
                    where we used the block diagonal structure of $\bar Q$. 
                    Moreover, since we assumed that
                    $
                        \proj_{\ker{L}} Q \proj_{\ker{L}} \succeq 0
                    $
                    it follows that also 
                    $
                        {\VL^\prime}^\top Q \VL^\prime
                            {}={}
                        {\VL^\prime}^\top \proj_{\ker{L}} Q \proj_{\ker{L}} \VL^\prime
                            \succeq
                        0
                    $
                    and since we assumed that 
                    \(
                        \lambda_{\rm min}(Q) < 0
                    \)
                    necessarily 
                    \(
                        \lambda_{\rm min}(Q)
                            {}={}
                        \lambda_{\rm min}(\VL^\top Q \VL)
                    \).
                \item \label{claim:QP:2}
                \textit{
                    The parameter
                    \(
                        \DRSrhomon
                    \)
                    is upper bounded by 
                    \(
                        \tfrac1{\nrm{L}^2}\lambda_{\rm min}(Q)
                    \).
                }\newline
                Since $(\MonA, \MonB(y^\star)) \in \parsumset$, 
                we have that
                \(
                    \MonA + \MonB(y^\star) \succeq 0
                \)
                and consequently it follows from \eqref{eq:def:parsum:prop} that
                \(
                    \MonA
                    \Box 
                    \MonB(y^\star)
                        {}\preceq{}
                    \MonA
                \), where
                \(
                    \MonA
                        {}={}
                    {L^\dagger}^\top Q L^\dagger
                \). 
                Since
                \(
                    L^\dagger
                        {}={}
                    \VL \SL^{-1} \UL^\top
                \)
                and 
                $\UL^\top \UL = \I$,
                this implies that
                \begin{align*}
                    \DRSrhomon
                        {}={}&
                    \lambda_{\rm min}\left(\UL^\top\left(
                        \MonA
                        \Box 
                        \MonB(y^\star)
                    \right)\UL\right)\\
                        {}\leq{}&
                    \lambda_{\rm min}\left(
                        \SL^{-1}
                        \VL^\top 
                        Q
                        \VL
                        \SL^{-1}
                    \right)
                        {}={}
                    \min_{s \in \R^r \backslash \{0\}}
                    \frac{
                        s^\top                         
                        \SL^{-1}
                        \VL^\top 
                        Q
                        \VL
                        \SL^{-1}
                        s
                    }{s^\top s}
                        {}={}
                    \min_{t \in \R^r \backslash \{0\}}
                    \frac{
                        t^\top
                        \VL^\top 
                        Q
                        \VL
                        t
                    }{t^\top \SL^2 t},
                \end{align*}
                where we defined $t \coloneqq \SL^{-1}s$ in the final step.
                Since 
                \(
                    \lambda_{\rm min}(                                
                        \VL^\top 
                        Q
                        \VL
                    ) 
                        <
                    0
                \)
                due to \Cref{claim:QP:1}
                and since $\SL^2 \preceq \|L\|^2\I$, it follows that
                \[
                    \DRSrhomon
                        {}\leq{}
                    \tfrac{1}{\|L\|^2}
                    \min_{t \in \R^r \backslash \{0\}}
                    \frac{
                        t^\top
                        \VL^\top 
                        Q
                        \VL
                        t
                    }{t^\top t}
                        {}={}
                    \tfrac{1}{\|L\|^2}
                    \lambda_{\rm min}\left(
                        \VL^\top 
                        Q
                        \VL
                    \right)
                    {}={}
                    \tfrac{1}{\|L\|^2}
                    \lambda_{\rm min}\left(
                        Q
                    \right).
                \]
            \end{claims}
            Finally, it follows from \Cref{claim:QP:2}
            that 
            \(
                \I + \gamma Q
                    \succeq
                \bigl(1 + \gamma \lambda_{\rm min}(Q)\bigr) \I
                    \succeq
                \bigl(1 + \gamma \DRSrhomon \|L\|^2\bigr) \I
            \),
            showing that $\I + \gamma Q$ is invertible 
            since $\gamma < \nicefrac{1}{[\DRSrhomon]_-\|L\|^2}$ as in \eqref{eq:QP:convergence:stepsizes}.
            \qedhere
        \end{itemize}
    \end{proofitemize}
\end{appendixproof}

\begin{appendixproof}{ex:constrainedQP:Semi}[ ]
    The projections are given by 
    $
        \proj_{\rge{L^\top}} = \diag(1, 1, 0)
    $
    and
    $
        \proj_{\ker{L}} = \diag(0, 0, 1)
    $, leading to 
    $
        \proj_{\range{L^\top}} Q \proj_{\ker{L}} = 0
    $
    and
    $
        \proj_{\ker{L}} Q \proj_{\ker{L}} = \diag(0,0,2) \succeq 0
    $.
    Therefore, we can invoke \Cref{prop:constrainedQP} to establish the assertions.
    In doing so, the following intermediate computations arise:
    $
        L^\dagger =         
        \begin{bsmallmatrix}
            1 & 0 & 0\\
            -\frac14 & 1 & 0
        \end{bsmallmatrix}^\top
    $, 
    $
        \|L\|
            =
        \tfrac18(1 + \sqrt{65})
    $,
    $
        \smash{
        \MonA
            {}={} 
        {L^\dagger}^\top Q L^\dagger
            {}={}
        \tfrac{1}{16}
        \begin{bsmallmatrix}
            16 & -4\\
            -4 & -15
        \end{bsmallmatrix}
        }
    $,
    $
        L^\top \MonA L = \diag(1, -1, 0) 
    $,
    $
        \ComA^\prime
            =
        \proj_{\ker{L}} Q^{-1} \proj_{\ker{L}} 
            =
        \diag(0,0,\nicefrac12)
    $,
    $
        \MonB
            =
        \diag(0, \nicefrac32)
    $,
    $
        X^\prime
            =
        \begin{bsmallmatrix}
            0&
            0&
            1
        \end{bsmallmatrix}^\top
    $,
    $
        {\VL'}^\top
        Q^{-1}
        \VL'
            {}={}
        \nicefrac12
    $,
    $
        \MonA \Box \MonB = \diag(0, -3)
    $
    and that
    $
        \DRSrhomon
            =
        \lambda_{\rm min}
        \left(\UL^\top\left( 
            \diag(0, -3)
        \right)\UL\right)
            =
        -3
    $
    since $Y$ is an orthonormal basis with full rank.
\end{appendixproof}

\begin{appendixproof}{ex:constrainedQP:semi}[ ]
    The proof is analogous to the one of \Cref{ex:constrainedQP:Semi}. It holds that
    $
        \proj_{\rge{L^\top}} = \diag(1, 1, 0)
    $
    and
    $
        \proj_{\ker{L}} = \diag(0, 0, 1)
    $.
    Therefore, 
    $
        \proj_{\range{L^\top}} Q \proj_{\ker{L}} = 0
    $
    and
    $
        \proj_{\ker{L}} Q \proj_{\ker{L}} = \diag(0,0,1) \succeq 0
    $.
    Hence, we can invoke \Cref{prop:constrainedQP} to establish the assertions.
    In doing so, the following intermediate computations arise:
    $
        L^\dagger = 
        \tfrac16         
        \begin{bsmallmatrix}
            2 & 2 & 2\\
            0 & 3 & -3\\
            0 & 0 & 0
        \end{bsmallmatrix}
    $,
    $
        \|L\| = \sqrt{3}
    $,
    $
        \MonA
            {}={} 
        {L^\dagger}^\top Q L^\dagger
            {}={}
        -\I_3
    $,
    $
        \ComA^\prime
            =
        \proj_{\ker{L}} Q^{-1} \proj_{\ker{L}} 
            =
        \diag(0,0,1)
    $,
    $
        \MonB
            =
        2\I_3
    $,
    $
        \VL^\prime
            =
        \begin{bsmallmatrix}
            0&
            0&
            1
        \end{bsmallmatrix}^\top
    $,
    $
        Y
            {}={}
        \begin{bsmallmatrix}
            \nicefrac1{\sqrt{3}} & 0\\
            \nicefrac1{\sqrt{3}} & \nicefrac1{\sqrt{2}}\\
            \nicefrac1{\sqrt{3}} & -\nicefrac1{\sqrt{2}}
        \end{bsmallmatrix}
    $,
    $
        {\VL'}^\top
        Q^{-1}
        \VL'
            {}={}
        1
    $
    and
    $
        \MonA \Box \MonB = -2\I_3
    $.
\end{appendixproof}

\begin{appendixproof}{cor:CP:semi}
    To show that the claims from \Cref{thm:CP:full} hold, it suffices to verify that \Cref{ass:SWMVICP,ass:CP:stepsize:rule,ass:CP:relaxation:rule} hold. 
    By \Cref{prop:semi:WMVI},
    the primal-dual operator $\Tpd$ has $\DRSRho$\hyp{}oblique weak Minty solutions at \(\pazocal{S}^\star\), where $\DRSRho$ is defined as in \eqref{eq:WMSDRS} with parameters $\DRSrhocom, \DRSrhocom^\prime, \DRSrhomon$ and $\DRSrhomon^\prime$ given by
    \begin{align*}
        \DRSrhocom = \comA \Box \comB,
        \quad
        \DRSrhomon = \monA \Box \monB,
        \quad
        \DRSrhocom^\prime = 
        \begin{cases}
            \infty,  &\quad \text{if } \rank L = n,\\
            \comA, &\quad \text{if } \rank L < n,
        \end{cases}, 
        \quad 
        \DRSrhomon^\prime = 
        \begin{cases}
            \infty,  &\quad \text{if } \rank L = m,\\
            \monB, &\quad \text{if } \rank L < m.
        \end{cases}
        \numberthis\label{eq:DRSRho:semimonotonicity}
    \end{align*}
    To show that \Cref{ass:SWMVICP} holds, we need to verify that these parameters satisfy \eqref{eq:existence:condition}, i.e., that
    \(
        [\DRSrhocom]_-[\DRSrhomon]_- < \nicefrac{1}{4\nrm{L}^2}
    \)
    and that\\[-0.35cm]
    \noindent\begin{minipage}[t]{.475\linewidth}
        \begin{equation}
            \DRSrhocom^\prime
                {}\geq{}
            \begin{cases}
                0, &\;\; \text{if } \DRSrhocom \geq 0,\\
                {\frac{\DRSrhocom}{1 - \DRSrhocom\DRSrhomon\nrm{L}^2}}, &\;\; \text{if } \DRSrhocom < 0, \DRSrhomon < 0,\\[0.15cm]
                {\frac{\DRSrhocom}{1 - \DRSrhocom\DRSrhomon\sigma_d^2}}, &\;\; \text{if } \DRSrhocom < 0, \DRSrhomon \geq 0,
            \end{cases}
            \label{eq:proof:existence:condition:1}
        \end{equation}
    \end{minipage}\hspace{0.05\linewidth}%
    \begin{minipage}[t]{.475\linewidth}
        \begin{equation}
            \DRSrhomon^\prime
                {}\geq{}
            \begin{cases}
                0, &\;\; \text{if } \DRSrhomon \geq 0,\\
                {\frac{\DRSrhomon}{1 - \DRSrhocom\DRSrhomon\nrm{L}^2}}, &\;\; \text{if } \DRSrhomon < 0,\DRSrhocom < 0,\\[0.15cm]
                {\frac{\DRSrhomon}{1 - \DRSrhocom\DRSrhomon\sigma_d^2}}, &\;\; \text{if } \DRSrhomon < 0, \DRSrhocom \geq 0.
            \end{cases}
            \label{eq:proof:existence:condition:2}
        \end{equation}
    \end{minipage}
    The first inequality holds owing to \cref{eq:DRSRho:semimonotonicity,ass:CP:semi:case:nonmon}.
    Before verifying \eqref{eq:proof:existence:condition:1} and \eqref{eq:proof:existence:condition:2}, we first show three intermediate results.
    \begin{claims}
        \item \label{claim:CP:cond-positive}
        \textit{%
            Condition \eqref{eq:proof:existence:condition:1} holds for all $\DRSrhocom^\prime \geq 0$ and \eqref{eq:proof:existence:condition:2} holds for all $\DRSrhomon^\prime \geq 0$.
        }
        \newline
        We already established that 
        \(
            [\DRSrhocom]_-[\DRSrhomon]_- < \nicefrac{1}{4\nrm{L}^2}
        \)
        and therefore       
        \(
            1 - \DRSrhocom\DRSrhomon\nrm{L}^2 > \nicefrac{3}{4\nrm{L}^2}
        \)
        when
        $\DRSrhocom < 0$ and $\DRSrhomon < 0$.
        Moreover, it holds by construction that
        \(
            1 - \DRSrhocom\DRSrhomon\sigma_d^2 \geq 1
        \)
        when $\DRSrhocom\DRSrhomon \leq 0$.
        Consequently, the fractions appearing in the cases of \eqref{eq:proof:existence:condition:1} and \eqref{eq:proof:existence:condition:2} are strictly negative as their numerators are strictly positive and their denominators strictly negative, proving the claim.

        \item \label{claim:CP:parsum-ineq}
        \textit{It holds that }
        \(
            \DRSrhocom = \comA \Box \comB \leq \comA
        \)
        \textit{and}
        \(
            \DRSrhomon = \monA \Box \monB \leq \monB.
        \)
        \newline
        This claim follows by definition of the parallel sum,
        since either $\comA = \comB = 0$ or $\comA + \comB > 0$ and either $\monA = \monB = 0$ or $\monA + \monB > 0$ owing to \cref{ass:SWMVICPsemimonotone}.

        \item \label{claim:CP:singleton}
        \textit{The graphs $\gph(\proj_{\rge{L^\top}} \circ A^{-1})$ and $\gph(\proj_{\rge{L}} \circ B)$ are not singletons.}\newline
        By \cref{ass:CP:2} the resolvents $J_{\gamma A}$ and $J_{\tau B^{-1}}$ have full domain. Consequently, $\gamma\id + A^{-1}$ and $\tau\id + B$ have full range \cite[Prop. 2.11(ii)]{bauschke2021Generalized},
        which implies that
        $
            \rge{\proj_{\rge{L^\top}}\circ (\gamma\id + A^{-1})} = \rge{L^\top}
        $
        and
        $
            \rge{\proj_{\rge{L}}\circ (\tau\id + B)} = \rge{L}
        $.
        Since $\rge{\proj_{\rge{L^\top}}\circ (\gamma\id + A^{-1})}$
        is a singleton whenever $\gph(\proj_{\rge{L^\top}} \circ A^{-1})$
        is a singleton
        and $\rge{\proj_{\rge{L}}\circ (\tau\id + B)}$
        is a singleton whenever $\gph(\proj_{\rge{L}} \circ B)$
        is a singleton, the claim is established.
    \end{claims}
    Using these claims, we will now show that condition \eqref{eq:proof:existence:condition:1} holds.
    For any $\DRSrhocom^\prime \geq 0$ this follows directly from \Cref{claim:CP:cond-positive}.
    Consider the case where $\DRSrhocom^\prime < 0$, which by \eqref{eq:DRSRho:semimonotonicity} implies that $\rank{L} < n$, $\DRSrhocom^\prime = \comA$ and $\comA < 0$.
    Then, it follows from \Cref{ass:SWMVICPsemimonotone} that (i) $\comB > -\comA > 0$, which implies by definition of $\DRSrhocom$ from \eqref{eq:DRSRho:semimonotonicity} that $\DRSrhocom = \comA \Box \comB < 0$,
    and that (ii) either $\monA = \monB = 0$ or $\monA + \monB > 0$.
    \begin{proofitemize}
        \item 
        $\monA = \monB = 0$: By definition of $\DRSrhomon$ from \eqref{eq:DRSRho:semimonotonicity} it then follows that $\DRSrhomon = \monA \Box \monB = 0$
        and that condition \eqref{eq:proof:existence:condition:1} reduces to
        $
            \comA
                {}\geq{}
            \DRSrhocom
        $,
        which holds by \cref{claim:CP:parsum-ineq}.
        
        \item 
        $\monA + \monB > 0$: Then, it holds that $\DRSrhomon \neq 0$, and since $\DRSrhocom < 0$ condition \eqref{eq:proof:existence:condition:1} reduces to
        \begin{align*}
            \comA \geq \frac{\DRSrhocom}{1 - \DRSrhocom([\DRSrhomon]_+\sigma_d^2 - [\DRSrhomon]_-\|L\|^2)}.
        \end{align*}
        Plugging in the definitions from \eqref{eq:DRSRho:semimonotonicity}, this can be equivalently stated as
        \begin{align*}
            \comA
                {}\geq{}
            \frac{\comA \Box \comB}{1 - (\comA \Box \comB)\Bigl([\monA \Box \monB]_+\sigma_{d}^2 - [\monA \Box \monB]_-\nrm{L}^2\Bigr)}
                {}={}
            \frac{\comA\comB(\monA + \monB)}{\comA + \comB - \comA \comB\Bigl([\monA\monB]_+\sigma_{d}^2 - [\monA\monB]_-\nrm{L}^2\Bigr)},
        \end{align*}
        where we multiplied the numerator and denominator by $(\monA + \monB)(\comA + \comB) > 0$ in the equality.
        Multiplying both sides of the inequality by the strictly positive denominator (see the proof of \cref{claim:CP:cond-positive}) and reordering, we obtain
        \(
            \comA^2\left(
            \monA + \monB - [\monA\monB]_+ \comB \sigma_d^2 + [\monA\monB]_-\comB \nrm{L}^2\right) \geq 0
        \),
        which is vacuously satisfied when $[\monA\monB]_+ = 0$. Otherwise, it holds that $\monA > 0$ and $\monB > 0$, so that it becomes
        \(
            \monB + \monA(1 - \monB\comB\sigma_d^2) \geq 0
        \), which
        is satisfied owing to \Cref{lem:semi:bounds} and \cref{claim:CP:singleton}.
    \end{proofitemize}
    The argument for condition \eqref{eq:proof:existence:condition:2}
    is analogous, establishing that \Cref{ass:SWMVICP} holds.

    It is easy to verify that the stepsize intervals provided in \Cref{tab:CP:stepsize:rule:semi} match those from \Cref{ass:CP:stepsize:rule} when $\DRSrhocom$ and $\DRSrhomon$ are given by \eqref{eq:DRSRho:semimonotonicity}.
    Therefore, it only remains to verify that 
    that \Cref{ass:CP:relaxation:rule:semi} is equivalent to \Cref{ass:CP:relaxation:rule}. 
    First, observe that the definition of $\eta^\prime$ from \Cref{tab:CP:degen:semi} is obtained by plugging in $\DRSrhocom^\prime$ and $\DRSrhomon^\prime$ from \eqref{eq:DRSRho:semimonotonicity} into \Cref{tab:CP:degen}.
    Moreover, as a consequence of \cref{claim:CP:parsum-ineq},
    $
        \min\bigl\{\tfrac1\gamma \DRSrhocom, \tfrac1\tau \DRSrhomon\bigr\} \leq \min\bigl\{\tfrac1\gamma \comA, \tfrac1\tau \monB\bigr\} \leq \eta^\prime
    $
    and the following assertions hold.
    \begin{enumerate}
        \item If $\max\{\DRSrhocom, \DRSrhomon\} \leq 0$, then $\tfrac1\gamma \DRSrhocom + \tfrac1\tau \DRSrhomon \leq \eta^\prime$.
        \item if $\DRSrhocom\DRSrhomon \geq 0$, then $\tfrac1{2\gamma} \DRSrhocom + \tfrac1{2\tau} \DRSrhomon - \theta_{\gamma\tau}(\sigma) \leq \eta^\prime$ for any $\sigma \in (0, \nrm{L}]$, since
        \begin{align*}
            \theta_{\gamma \tau} (\sigma) 
                {}={}
            \sqrt{\left(\tfrac{1}{2\gamma}\DRSrhocom + \tfrac{1}{2\tau}\DRSrhomon\right)^2 - \tfrac{1}{\gamma\tau}\DRSrhocom\DRSrhomon(1 - \gamma\tau\sigma^2)}
                {}\geq{}
            \left|\tfrac{1}{2\gamma}\DRSrhocom + \tfrac{1}{2\tau}\DRSrhomon\right|.
        \end{align*}
    \end{enumerate}
    The claimed equivalence follows immediately from these two assertions, completing the proof.
\end{appendixproof}
\end{appendix}

\ifjota
\else
  \bibliographystyle{siamplain}%
  \bibliography{TeX/Minty-DRS.bib}%
\fi

\end{document}